\documentclass[11pt,a4paper]{article}
\title{\bf On a particle approximation to the Dean-Kawasaki type 
	equation with logarithmic interactions}
\author{Hao DING\footnote{Email: dinghao16@mails.ucas.ac.cn}
}

\usepackage[utf8]{inputenc}
\usepackage{mathtools}
\usepackage{hyperref}
\hypersetup{hypertex=true,
	colorlinks=true,
	linkcolor=blue,
	anchorcolor=blue,
	citecolor=blue}
\usepackage[nomsgs]{refcheck}
\usepackage{amssymb,amsmath,amsfonts,amsthm,array,bm,color}
\usepackage{mathrsfs}
\usepackage{bbm}
\def\N{\mathbb{N}}
\def\R{\mathbb{R}}
\def\E{\mathbb{E}}
\def\P{\mathbb{P}}
\def\T{\mathbb{T}}

\def\Z{\mathbb{Z}}
\def\G{\mathcal{G}}

\def\F{\mathcal{F}}
\def\C{\mathcal{C}}
\def\H{\textup{H}}

\def\L{\mathcal L}
\def\S{\mathbb{S}}

\def\d{\textup{d}}

\def\eps{\varepsilon}
\def\<{\langle}
\def\>{\rangle}

\def\W{\mathbb{W}}

\let\ra=\rightarrow
\def\one{\uppercase\expandafter{\romannumeral1}}
\def\two{\uppercase\expandafter{\romannumeral2}}
\def\three{\uppercase\expandafter{\romannumeral3}}
\def\four{\uppercase\expandafter{\romannumeral4}}

\newtheorem{theorem}{Theorem}
\newtheorem{lemma}[theorem]{Lemma}       
\newtheorem{corollary}[theorem]{Corollary}
\newtheorem{proposition}[theorem]{Proposition}
\newtheorem{remark}{Remark}

\newtheorem{definition}{Definition}
\newtheorem*{thm}{Theorem}

\begin{document}
	
	\maketitle
	\makeatletter 
	\renewcommand\theequation{\thesection.\arabic{equation}}
	\@addtoreset{equation}{section}
	\makeatother 
	
	\vspace{-8mm}
	\begin{abstract} 
		We consider a class of Dean-Kawasaki type equations on $\T$ with logarithmic repulsive interactions depending on the inverse temperature $\beta$ and a new spectral approximation to the noise part, which approximately features Otto's metric in $\P(\T)$. Following the idea of intrinsic constructions of Brownian motions on the Wasserstein space, we construct a class of particle models whose fluctuating hydrodynamic limits, denoted as $p_t^\beta$, are solutions to the martingale problems of this class of equations. Specifically, we give a quantitative convergence rate of the particle approximation, which allows us to identify a unique limit distribution depending on $\beta$. 
		
		As the inverse temperature rises, the regularizing effect of repulsive interactions becomes stronger. We prove that there exists three thresholds $0<\lambda_0\leq\lambda_1<\lambda_2$ depending on the noise such that, when $\beta>\lambda_0$, $p_t^\beta$ is a non-atomic measure process in $\P(\T)$; when $\beta>\lambda_1$, $p_t^\beta$ is absolutely continuous with respect to Lebesgue measure almost surely; when $\beta>\lambda_2$, the expectation of the R\'enyi entropy of $p_t^\beta$ satisfies an exponential decay estimate.
	\end{abstract}
	
	
	
	\tableofcontents
	\section{Introduction}
	\subsection{Background}
	\quad Dean-Kawasaki type equations are a class of nonlinear stochastic partial differential equations (SPDEs) arising in fluctuating hydrodynamics theory \cite{eyink1990,dean1996,kawasaki1998}. In general, one may consider the following Dean-Kawasaki type equation on $\T^d$
	\begin{equation}\label{onex1}
		\begin{split}
			\partial_t \mu &=K\Delta\mu-\nabla.(\mu\Gamma(\mu))-\nabla.(\sqrt{\mu} \dot{\xi}) , 
		\end{split}
	\end{equation}
	
	for space-time white noise $\dot{\xi}$ and $\alpha >0$, $\Gamma(\mu)$ is a measure dependent vector field on $\T^d$. We say a continuous measure-valued process $\mu_t$ is a solution to the martingale problem $(MP)_{\mu_0}^{\alpha}$ of $\eqref{onex1}$ with initial condition $\mu_0$ if there exists a filtered probability space $(\Omega, \mathcal{F}, \{\mathcal{F}_t\}_{t\geq 0}, \P)$ such that for all $\phi \in \C^{2}(\T)$ ,
	
	\begin{equation*}
		M_t(\phi):=\langle \mu_t,\phi\rangle-\langle\mu_0,\phi\rangle -\int_{0}^{t}\langle\Gamma(\mu_s) , \nabla \phi \cdot  \mu_s\rangle d s -\alpha \int_0^t \langle\mu_s,\Delta \phi\rangle d s
	\end{equation*}
	
	is a $\F_t-$adapted martingale, whose quadratic variation is given by 
	
	\begin{equation}\label{charac}
		< M_t(\phi)>=\int_0^t ||\nabla \phi||^2_{\L^2(\mu_s)}d s .
	\end{equation}
	We only consider the case when d=1 in this paper and always realize $\T=\R/\Z$ as $[0,1]$ in this paper. We also denote $Leb$ as the Lebesgue measure on $\T$.
	
	
	
	\quad Dean-Kawasaki type equations are closely related to many fields. For instance, the Large deviation principle for Dean-Kawasaki equations \cite{mariani2010,ZR22,FG23} can help predicting exponential decay of rare events in the interacting particle systems, so that it acts as a bridge between the macroscopic fluctuation theory \cite{RevModPhys.87.593} and fluctuating hydrodynamic theory. Also, Dean-Kawasaki equations are effective tools in
	accurately simulating the density fluctuations in systems of diffusing particles \cite{CF23}. As one of the primary motivations of this paper, Dean-Kawasaki type equations can be recognized as a class of SPDEs that characterize the infinite dimensional analogs of Brownian motions on the Wasserstein space, called as Wasserstein diffusions, since the quadratic variation \eqref{charac} coincides with Otto's formal Riemannian metric \cite{von2009,Otto01}. In particular, we are more interested in non-trivial Wasserstein diffusions $\mu_t$, i.e. measure without atoms, mainly from the following two perspectives. First, for a conservative stochastic particle system whose hydrodynamic limit is governed by a deterministic Fokker-Planck equation 
	$$\partial_t \mu= \nabla.(\mu\Gamma(\mu)),$$
	one would expect a regular solution in the continuum model in general cases. Similarly, we also expect to observe non-trivial macroscopic fluctuations around $(\mu_t)_{t\geq 0}$, when the microscopic noises can not be averaging out. The second perspective concerns the regularization effects of noise \cite{F2015,VM21}. Analogous to Brownian motions in Euclidean spaces, one would like to construct a Wasserstein diffusion that has similar regularization effects on some conservative SPDEs, and the necessary condition is the noise itself should support on measures without atoms.


	\quad However, as pointed out in \cite{konarovskyi2019,cornalba2019,fehrman2021}, well-posedness of \eqref{onex1} is challenging due to the conservative noise with the coefficient $\sqrt{\mu}$ which causes irregularity, nonlinearity and possible lack of Lipschitz continuity. It is not renormalisable by the theory of
	regularity structures \cite{hairer2014}. Even worse, establishing the well-posedness of the pure diffusive Dean-Kawasaki equation results in trivial solutions. In detail, when $\Gamma(\mu)=0$, it has been proven in \cite{konarovskyi2019} that there exists a unique measure-valued solution to the martingale problem of \eqref{onex1} if and only if
	$2\alpha \in \N^+$, and in this case, the solution is trivial (atomic), i.e. $\mu_t=\frac{1}{N}\sum_{i=1}^{N}\delta_{W_{\sqrt{N}t}^i}$ , where $\{W_t^i\}_{i=1,\ldots, N}$ are N independent standard Brownian particles. For this solution, the diffusive term $\frac{N}{2}\Delta \mu$ could be understood as an Ito's correction term induced by the noise, and the unique solution satisfies the following equation in stratonovich's form:
	\begin{equation}\label{onex2}
		\partial_t \mu =-\nabla.(\sqrt{\mu} \circ\dot{\xi}) .
	\end{equation}
	Moreover, it has been demonstrated \cite{Konarovskyi2018OnDD} that adding a "good" drift in form of gradient term $\nabla.(\mu\nabla \frac{\delta \mathcal{F}}{\delta \mu }(\mu))$ for some smooth functional $\mathcal{F}$ on the Wasserstein space to \eqref{onex2} can not change the nature of triviality. The triviality of solutions to the pure diffusive Dean-Kawasaki equations in turn suggests that adding  an irregular drift or renormalizing the noise is essential for constructing an analytically non-trivial diffusions on the Wasserstein space. An interesting example from \cite{fehrman2021} supports this notion. The well-posedness results in \cite{fehrman2021} apply to the following equation in stratonovich form on $\T$:
	\begin{equation*}\label{onex3}
		\partial_t \mu=\nu \partial_{xx} \mu -\partial_x(\sqrt{\mu}\circ \dot{\xi^\gamma}),
	\end{equation*}
	where $\xi^\gamma=\sum \frac{1}{|k|^\gamma}e_k B_t^k$ and 
	\begin{equation*}
		\begin{split}
			&e_k=\sqrt{2}\sin(2k\pi x),  \quad \ k=1,2...;\\
			&e_0=1;\\
			&e_k=\sqrt{2}\cos(2k\pi x), \quad \ k=-1,-2,... \ ,
		\end{split}
	\end{equation*}
	are canonical basis in $\L^2(\T, Leb)$. It reveals that adding an additional viscosity term and adopting a spectral approximation of $\dot{\xi}$ lead to a regular solution to the Dean-Kawasaki equation in sense of stochastic kinetic solution. The viscosity term can be interpreted as a negative gradient of the Boltzmann entropy under Otto's metric \cite{Otto01}. 
	
	\quad In this paper, by a particle approximation, we construct a class of non-trivial Langevin diffusions $p_t^\beta$ on the Wasserstein space satisfying Dean-Kawasaki type equations with a singular drift in form of the negative gradient of the logarithmic interaction energy and a new spectral approximation of the noise term $"\partial_x(\sqrt{\mu} \dot{\xi})"$, which can be understood as a realization of the intrinsic construction of a $Q-$Brownian motion $\textbf{B}^Q_t$ on the Wasserstein space. We also investigate the pathwise regularity of $p_t^\beta$ as the inverse temperature $\beta$ varies. In particular, decay estimates on the expectation of Boltzmann and R\'enyi entropy of $p_t^\beta$ are obtained.

	\subsection{Main ideas of the particle model}\label{subsec1.1}
	\quad In this subsection, we introduce the motivation of the particle model approximation. Following the idea of intrinsic constructions of Brownian motions on Riemannian manifolds \cite{eells1976stochastic,malliavin1974formules,eells2011wiener}, we aim to give a class of non-trivial diffusions on $\P(\T)$ satisfying Dean-Kawasaki type equations whose noise part approximately features Otto's metric. Specifically, we firstly adopt Lagrangian's viewpoint to interpret an approximation Wasserstein diffusion $\mu^\gamma_t$ as a pushforward measure process of a certain stochastic flow $X^\gamma_t$, i.e. $\mu^\gamma_t=(X^\gamma_t)_{\#}\mu_0$, starting from an initial measure $\mu_0$. Subsequently, we want to construct $X^\gamma_t$ as a stochastic flow on $\T$ generated by a so called stochastic moving frame $\{v_k(\mu^\gamma_t)\, \}_{k\in \Z}$ along $\mu^\gamma_t$, i.e. 
	\begin{equation}\label{xbeta}
		\left\{
		\begin{aligned}
			&dX^\gamma_t=\sum_{k=-\infty}^{+\infty}\frac{1}{|k|^\gamma}v_k(\mu^\gamma_t,X^\gamma_t)\circ dW_t^k,\\
			&X_0^\gamma(x)=x, \ \mu_0^\gamma=\mu_0, \ v_k(\mu_0^\gamma)=v_k(0), \ \ \ \forall \, k\in \N,
		\end{aligned}
		\right.
	\end{equation}
	where $\{v_k(\mu^\gamma_t)\, \}_{k\in \Z}$ are stochastically parallel translated in Ito's sense along $\mu^\gamma_t$ on $\P(\T)$ under Lott's Levi-civita connection \cite{lot2006som}. Without confusions on terminology, we will not check in this paper whether $X_t^\gamma$ satisfies a flow property, since we only use $X_t^\gamma$ to explain the motivation of our particle model. According to \cite{fan2021}, the well-posedness of the stochastic parallel translation along a pushforward measure process requires the driving flow to be diffeomorphic; also, a close form expression of $v_k(\mu^\gamma_t)$ seems impossible. This result, combined with non-triviality requirement on $\mu_t^\gamma$, shows that the construction of $X_t^\gamma$ is analytically difficult and intricate. In order to realize the main idea introduced above, we need to modify the original model \eqref{xbeta} so that it
	\begin{enumerate}
		\item does not exhibit mass concentration, which is a necessary condition for diffeomorphism, 
		\item satisfies a closed form equation, which means we do not need additional equations to determine the stochastic moving frames along with \eqref{xbeta}.
	\end{enumerate} 
Our strategy is to construct a particle model with some additional microscopic informations so that it can realize a macroscopic flow $X_t^\gamma$  satisfying these two requirements. 
	
	\quad In order to satisfy the second requirement, we adopt a similar construction of stochastic parallel translations on $\P_2(\R)$ instead of $\P(\T)$. We observed that, due to flatness of $\P_2(\R)$ in sense of Otto's Riemannian structure, stochastic parallel translations of an initial vector $v_k(0)=e_k(F_{\mu_0})$ along a pushforward measure process $\mu_t$ on $\P_2(\R)$ induced by a stochastic flow should be uniquely determined as $v_k(t)=e_k(F_{\mu_t})$, where $F_\mu$ is the distribution function of $\mu$ on $\R$. Therefore, by using this form of stochastic moving frames, we want to construct the noise part of $X_t^\gamma$ on $\T$ as 
	\begin{equation}\label{xbeta2}
		\sum_{k=-\infty}^{+\infty}\frac{1}{|k|^\gamma}e_k\big(F_{\mu^\gamma_t}(X^\gamma_t)\big)\circ dW_t^k,
	\end{equation} 
	where $F_{\mu}$ is defined as 
	\begin{equation*}
		F_{\mu}(x):=\int_{0}^x\mu(dy) , \quad \quad \forall \, x\in [0,1).
	\end{equation*} 
	Under certain conditions, a formal application of Ito's formula suggests that $\mu_t^\gamma$ does share an identical fluctuation with the solution to the martingale problem of \eqref{onex1} on $\T$ for $\gamma= 0$. Although it is impossible to rigorously achieve this consistency, it still demonstrates that the construction on the noise part of $\mu_t^\gamma$ indeed approximately features Otto's metric. Details about this construction of stochastic moving frames and formal consistency will be presented in Section \ref{sec1.2}.
	
	\quad The left requirement is to avoid mass concentrations. To fulfill this macroscopic property, we introduce an additional repulsive interaction into the microscopic particle model to prevent particle colliding within a finite time. In the prior work \cite{rog1993, li2020}, they consider the generalised Dyson Brownian motions $\lambda_N(t)=(\lambda_N^i(t))_{i=1}^N\in \R^N$, given by
	\begin{equation*}
		d\lambda_N^i(i)=\frac{1}{\sqrt{N}}dW_t^i+\frac{1}{N}\sum_{j\neq i, j=1}^{N}\frac{1}{\lambda_N^i(t)-\lambda_N^j(t)} dt-\frac{1}{2}V'(\lambda_N^i(t)) dt, \ \ \ i=1,\ldots,N,
	\end{equation*}
	with some initial value $(\lambda_N^i(0))_{i=1}^N$. For Lyapunov function $f(x_1,\ldots,x_N)=\frac{1}{N}\sum_{i=1}^{N}V(x_i)-\frac{1}{N^2}\sum_{i\neq j}\ln|x_i-x_j|$, it is proved in \cite{rog1993,li2020} that $\E[f(\lambda_N(t))]\leq C_0+C t$, where $C$ is independent of the particle number $N$. This uniform estimate indicates the absence of collisions almost surely within finite time even as $N\to \infty$. Motivated by this observation, we also add a logarithmic repulsion on $\T$ into the particle approximation.

	\subsection{Statements of main results}
	\quad Following the ideas introduced above, we construct a particle model $X_N(t)=(X_N^i(t))_{i=1}^N\in \R^N$ approximating $X_t^\gamma$. Analogous to the form of the noise part \eqref{xbeta2}, we construct the stochastic moving frame as $\displaystyle \{e_k(F_{\bar{L}_N(t)}))\}_{k\in \Z}$, where $\bar{L}_N(t)$ is the associated empirical measure of $X_N(t)$ on $\R$. We expect $\{X^i_N(t)\}_{i=1}^N$ keeps order so that
	$$e_k(F_{\bar{L}_N(t)}(X_N^i(t)))=e_k(F_{\bar{L}_N(0)}(X_N^i(0)))=e_k(\frac{i}{N}), \ \ \ \forall \, i=1,2,\ldots,N.$$  
	Let $\{W^k_t\}_{k\in \Z}$ be a sequence of independent standard Brownian motions defined
on a filtered probability space $(\Omega,\mathcal{F}, \{\mathcal{F}_t\}_{t\geq 0},\P)$. Given the particle number $N$ and $\gamma>\frac{3}{2}$, we construct the model as the solution to the following SDEs
	\begin{equation}\label{model1}
		dX_N^i(t)=\frac{\beta}{N}\sum_{j\neq i, j=1}^{N}\cot(\pi(X_N^i(t)-X_N^j(t)))dt+\sum_{k=-N}^{N}\frac{1}{|k|^\gamma}e_k(\frac{i}{N}) \circ  dW_t^k,  \ \ i=1,2,\ldots,N
	\end{equation}
	with the initial value $X_N(0)$ and the inverse temperature $\beta>0$. Let the integer function $\displaystyle [\,\cdot \,]:\, \R\ra \N$ be
	\begin{equation*}
		[x] =
		\left\{
		\begin{aligned}
			&x-1,   \quad \quad x\in \N;\\
			&\max\{n\in \N |\, n<x\}, \ \mathrm{otherwise},
		\end{aligned}
		\right.
	\end{equation*}
	and $\{x\}:=x-[x]$. We define the empirical measure of the projection of $X_N(t)$ on $\T$ 
	\begin{equation*}
		L_N(t):=\frac{1}{N}\sum_{i=1}^{N}\delta_{\{X_N^i(t)\}} ,
	\end{equation*}
	
	 and introduce several constants
	$$K_1^\gamma=\sum_{i=1}^{\infty}\frac{4\pi^2}{i^{2\gamma-2}}; \ \ \  K_2^\gamma=\sum_{i=1}^{\infty}\frac{1}{i^{2\gamma}}; \ \  \ \lambda_0^\gamma=(8\pi \sqrt{K_1^\gamma K_2^\gamma})\wedge( \frac{1}{2}\pi K_1^\gamma).$$
In this paper, for a measure $\mu=\rho dx \in \P_{ac}(\T)$, the Boltzmann entropy is introduced as $\displaystyle\textbf{S}(\mu)=\int_0^1 \ln \rho d \mu$ and the R\'enyi entropy is introduced as $\displaystyle\textbf{H}_{p}(\mu)=\int_0^1\rho^{p-1} d\mu$ for $p>1$. The interaction energy is 
\begin{equation}\label{LIE}
	\textbf{W}(\mu)=\int_{0}^1\int_{0}^1\ln|e^{2\pi i x}-e^{2\pi i y}| d\mu(x) d\mu(y).
\end{equation}
We also define the mollified Boltzmann entropy $S^N(L_N(t))$ by
$$S^N(L_N(t))=\sum_{j=1}^{N}\frac{1}{N}\ln\frac{2\pi}{N|e^{2\pi i X_N^{j}(t)}-e^{2\pi i X_N^{j-1}(t)}|},$$ 
where we just denote $X_N^0$ as $X_N^N-1$ for simplicity of notations. 

\quad Our first result is about the well-posedness of the particle model. Let 
	$$\Delta_N=\{(x_i)_{1\leq i\leq N} \in \R^N:x_1<x_2<...<x_N< x_1+1\}.$$

	\begin{theorem}\label{introthm1}
		Suppose that $ \gamma>\frac{3}{2}$ and $\beta>\lambda_0^\gamma$. If  the initial value satisfies $X_N(0)\in \Delta_N$, then there exists a unique strong solution $X_N(t)$, which takes value in $\Delta_N$, to SDE \eqref{model1} with infinite lifetime. Moreover, if $\beta>8\pi\sqrt{K_1^\gamma K_2^\gamma}$, then $\textup{\textbf{W}}(L_N(t))-C_{\gamma,\beta}t$ is a supermartingale for some constant $C_{\gamma,\beta}$ depending on $\beta$ and $\gamma$; if $\beta>\frac{1}{2}\pi K_1^\gamma$, then $S^N(L_N(t))-2\beta \pi t$ is a supermartingale and for $N$ large enough, it holds
		\begin{equation}\label{ugr'}
			\E[S^N(L_N(t))]\leq e^{-2 C_1^\beta t}S^N(L_N(0)) +(1-e^{-2C_1^\beta t})\frac{\beta \pi}{C_1^\beta},  \ \ \ \ \forall \, t\geq0,
		\end{equation}
		where $C_1^\beta=\frac{\beta-\frac{1}{2}\pi K_1^\gamma}{\pi}$.
	\end{theorem}

	\quad The first part of this result is about the well-posedness of \eqref{model1}. It shows that the system $(X_N^i)_{i=1}^{ N}$ is actually constrained to the torus since $|X_N^1(t)-X^N_N(t)|<1$, which means our construction of the stochastic moving frame via stochastic parallel translations on $\P(\R)$ is adapted to the projection model $(\{X^i_N(t)\})_{i=1}^N$ on $\T$. The proof of well-posedness is based on the Lyapunov function method \cite{rog1993,Anderson_Guionnet_Zeitouni_2009,li2020}. For the system \eqref{model1}, we choose two Lyapunov functions, whose macroscopic counterparts are the interaction energy $\textup{\textbf{W}}(\mu)$ and the Boltzmann entropy $\textup{\textbf{S}}(\mu)$. Both of them are able to indicate not only the absence of microscopic collisions between particles but also the absence of macroscopic mass concentrations, once we get a growth control on them (in sense of conditional expectation) independent of the particle number. The second part of this result gives the growth control in detail, which is also for preparation of the results on the regularity of $p_t^\beta$. We will discuss about this part of results after establishing the convergence of the empirical measures.

	\quad Denote the 2-Wasserstein distance in $\P_2(\R)$ and $\P(\T)$ respectively as $\bar{W}_2$ and $W_2$. For the canonical space $\Omega=\C([0,\infty),\P(\T))$  endowed with the metric $$d(p^1,p^2):=\sum_{n=1}^{\infty}\frac{1}{2^n}(\max_{t\in[0,n]}W_2(p^1_t,p^2_t)\wedge 1),$$
	and $P_1,P_2\in \P(\C([0,\infty),\P(\T)))$, we define the 2-Wasserstein distance $\mathbb{W}_2(P_1,P_2)$ by coupling as
	\begin{equation}\label{defW}
		\mathbb{W}_2(P_1,P_2):=(\E_{\Omega}[d^2(p_t^1,p_t^2)])^{\frac{1}{2}},
	\end{equation}
	for any two variables $p^1 , \, p^2 : \ (\Omega, \mathcal{F},P)\to C([0,\infty),\P(\T))$, whose laws are $P_1,\, P_2$ respectively. We then prove the convergence (in distribution) of the empirical measure process $L_N(t)$ in $\C([0,\infty), \P(\T))$ and give the rate of convergence.
	\begin{theorem}\label{ratecon}
		Suppose that $\gamma >\frac{3}{2}$ and $\mu_0\in \P(\T)$. If $\beta>\lambda_0^\gamma$ and the initial empirical measures of the particle model satisfy
		\begin{equation}\label{cod2'}
			X_N(0)\in \Delta_N \ \text{and} \  \lim_{N\to \infty}\bar{W}_2(\bar{L}_N(0),\mu_0)=0,
		\end{equation}
		then $\{P^N\}$ is a Cauchy sequence in $\P(\C([0,\infty), \P(\T)))$ endowed with $\W_2$. Moreover, we can arrange the initial positions for the particles such that for $\forall \ N<M$,
		\begin{equation}\label{prop7a}
			\E[W_2^2(L_N(t), L_M(t))]\leq  \frac{C\beta t}{N}+\frac{4}{N},
		\end{equation}
		and 
		\begin{equation}\label{prop7b}
			{\mathbb{W}}^2_2(P^N,P^M)<C N^{-\frac{1}{2}},
		\end{equation}
		when $N$ is large enough.
	\end{theorem}

		 \quad In general, one may derive the hydrodynamic limits or the fluctuating hydrodynamic limits by tightness arguments so that a subsequence convergence can be obtained; and the uniqueness of the limits relies on the uniqueness of solutions to the limit equation. Instead, we derive the Cauchy convergence of the distribution $\{P^N\}$ of $L_N(t)$ in $\P(\C([0,\infty), \P(\T)))$ endowed with a 2-Wasserstein distance $\W_2$. 
		 
	\quad The proof of \eqref{prop7a} is intuitively based on a contraction estimate on the transport distance between two gradient flows $\mu_t^1, \,  \mu_t^2\in \P_{ac}(\T)$ of $\textbf{W}$, that is,
	\begin{equation*}
		\begin{split}
			d_tW_2^2(\mu_t^1, \mu_t^2)&= \int_{[0,1]^2}\<x_1-x_2, \H\mu_t^1(x_1)\>d\gamma_t^{1,2}(x_1,x_2)+  \int_{[0,1]^2}\<x_2-x_1, \H\mu_t^2(x_2)\>d\gamma_t^{2,1}(x_2,x_1)\leq 0,
		\end{split}
	\end{equation*}
	by the derivative formula of $W_2$ (see \cite{villani2009optimal} Theorem 23.9 or \cite{Ambrosio2005GradientFI} Theorem 8.4.7) and the geodesically convexity of $\textbf{W}$ (see \cite{V03, Ambrosio2005GradientFI, Carrillo2011AMA, li2020}), where $\gamma_t^{i,j}$ is the optimal transport plan from $\mu_t^i$ to $\mu_t^j$ and $\H\mu$ is the Hilbert transform of $\mu$ on the torus, given by
	$$\H\mu(x)=P.V.\int_{0}^1\cot(\pi(x-y))d\mu(y).$$
	However, this estimate does not apply to the two empirical measures in our case. Also, the optimal mass transport plan from projected particles $(\{X^i_N\})_{i=1}^N$ to $(\{X^i_M\})_{i=1}^M$ is not explicitly known, which causes technical complexity in representing $W_2^2(L_N(t),L_M(t))$. To overcome these difficulties, we consider the empirical measure of \eqref{model1} on the universal space $\R$, denoted as
	$$\bar{L}_N(t)=\frac{1}{N}\sum_{i=1}^{N}\delta_{X_N^i(t)}.$$
	One of the key estimates in Proposition \ref{ratecon} is the following control along the empirical measure processes $\bar{L}_N(t)$ and $\bar{L}_M(t)$,
	\begin{equation*}\label{S1f}
		\int_{[0,1]^2}\<x_1-x_2, \H\bar{L}_N(t)(x_1)\>d\gamma_t^{1,2}(x_1,x_2)+  \int_{[0,1]^2}\<x_2-x_1, \H\bar{L}_M(t)(x_2)\>d\gamma_t^{2,1}(x_2,x_1)\leq \frac{2}{N}+\frac{2}{M},
	\end{equation*}
where $\gamma_t^{1,2}$ is an optimal transport plan from $\bar{L}_N(t)$ to $\bar{L}_M(t)$ and $\gamma_t^{2,1}$ is an optimal transport plan from $\bar{L}_M(t)$ to $\bar{L}_N(t)$. See Lemma \ref{lem3} for a proof. Based on the literatures we know about the convergence rates of empirical measure processes in some related models, for instance, unitary Dyson Brownian motions, the estimate \eqref{prop7a} on $\E[W_2^2(L_N(t), L_M(t))]$ seems to be new. Note that the noises are not averaged out as $N$ goes to infinity in our model, which increases the instability in approximation. We point out that the convergence rate estimate \eqref{prop7a} strongly relies on our special noise in form of $e_k(\frac{i}{N})dB_t^k$ in the particle model \eqref{model1}, which produces smaller variation around the limit than other common noises, for instance $e_k(X^i_N(t))dB_t^k$. In such sense, our model might be a good approximator of Dean Kawasaki type systems. 
	
	\quad As a consequence, the distributions $\{P^N\}$ admits a unique limit distribution $P$. Given an initial measure $\mu_0$ and $\gamma >\frac{3}{2}$, we can apply Skorohod representation theorem to get a unique (in distribution) fluctuating hydrodynamic limit $p_t$ with distribution $P$, associated with a new probability space $(\tilde{\Omega},\tilde{\mathcal{F}}, \{\tilde{\mathcal{F}}_t\}_{t\geq 0},\tilde{\P})$, corresponding to a construction of the initial empirical measures $\{L_N(0)\}$. Actually, we prove in Proposition \ref{prop2} that the limit distribution $P$ is independent of the initial empirical measures $\{L_N(0)\}$ under certain admissibility conditions. This proposition allows us to identify a unique limit distribution $P_\beta$ of our particle model \eqref{model1} only depending on macroscopic parameters $\beta>\lambda_0^\gamma$ and $\gamma>\frac{3}{2}$. We also prove in Proposition \ref{prop5.1} that as $\beta'$ approaching $\beta$, $P_{\beta'}$ converges to $P_\beta$ under $2-$Wasserstein distance $\W_2$ on $\P(\C([0,\infty),\P(\T)))$.
	
	
	\quad Theorem \ref{introthm1} also implies the fluctuating hydrodynamic limit does not have mass concentrations. Indeed, when $\beta>8\pi\sqrt{K_1^\gamma K_2^\gamma}$, we rigorously prove in Proposition \ref{noatom} that $p^\beta_t$ contains no atoms for all $t\in [0,\infty)$ almost surely, provided that the initial condition $\mu_0$ is regular enough. Roughly saying, the proof is based on a simple fact that for a sequence of measure $\{\mu_n\}_{n\in \N}$ weakly converging to $\mu$ in $\P(\T)$ with an atom at some point $x\in \T$, they must gather in a small neighbourhood around $x$ with a very high density, which will strongly increase the logarithmic interaction energy of $\mu_n$ and hence violate the assertion that $\textbf{W}(L_N(t))-C_{\gamma,\beta} t$ is a supermartingale. When $\beta>\frac{1}{2}\pi K_1^\gamma$, we prove in Proposition \ref{abs} that $p^\beta_t$ is absolutely continuous w.r.t Leb for all $t\geq 0$ almost surely, based on the assertion that $S^N(L_N(t))-2\beta \pi t$ is a supermartingale. 
	
	\begin{remark}
		When $\gamma-\frac{3}{2}$ is small enough, we have $8\pi\sqrt{K_1^\gamma K_2^\gamma}<\frac{1}{2}\pi K_1^\gamma$, then $p_t^\beta$ is almost surely non-atomic for $\beta\in (\lambda_0^\gamma, \frac{1}{2}\pi K_1^\gamma]$ and absolutely continuous for $\beta\in (\frac{1}{2}\pi K_1^\gamma, \infty)$. While when $8\pi\sqrt{K_1^\gamma K_2^\gamma}\geq \frac{1}{2}\pi K_1^\gamma$, we have $\lambda_0^\gamma=\frac{1}{2}\pi K_1^\gamma$ and $p_t^\beta$ is almost surely absolutely continuous for $\beta\in (\lambda_0^\gamma, \infty)$. This result fits our intuition that as $\gamma$ becomes smaller, the noise becomes rougher in space, which tends to make the fluctuating hydrodynamic limit less regular. However, due to the limitations of our method, we could neither determine the exact pathwise regularity of  $p^\beta_t$ between these thresholds nor the sharp value of thresholds. 
	\end{remark}	
	
	\quad Based on Proposition \ref{noatom} and Proposition \ref{abs}, we prove that 
		\begin{theorem}\label{maintheorem}
		Suppose that $\beta>\lambda_0^\gamma, \, \gamma >\frac{3}{2}$ and $\mu_0$ satisfies condition 
		\begin{equation}\label{regucon}
			\textup{\textbf{(J):}}  \ \ 	\mu_0 \textup{\ has\ a\ density\ } \rho_0 \textup{\ satisfying \ } ||\rho_0||_{\L^{\infty}}<\infty,
		\end{equation}   
	then the fluctuating hydrodynamic limit $p^\beta_t$ of the particle model \eqref{model1} is a non-trivial solution to the martingale problem (Definition \ref{MPs}) of the following Dean-Kawasaki type equation
			\begin{equation}\label{withrepul}
			\partial_t\mu=K_2^\gamma \partial_{xx}\mu-\beta\partial_x(\mu \H\mu) dt-\partial_x(\mu \cdot \dot{\xi}_{\mu}^\gamma),
		\end{equation}
	\end{theorem}
The noise part is given by
\begin{equation*}\label{noise1}
	\dot{\xi}_{\mu}^\gamma(t,x)=\sum_{k=-\infty}^{+\infty}\frac{1}{|k|^{\gamma}}e_k(F_{\mu_t}(x))\dot{W}_t^k.
\end{equation*}
which is natural according to the construction of the stochastic moving frame. The definitions of $\dot{\xi}_{\mu}^\gamma$ and the martingale problem of \eqref{withrepul} will be introduced in Subsection \ref{sec2}. We also convince the reader that the noise in \eqref{withrepul} is indeed a spectral approximation to the original Dean-Kawasaki equation by the formal consistency with \eqref{onex1} as $\gamma \to 0$. 
	
	\vskip 3mm

	\quad As the inverse temperature rises, the regularization effect of repulsive interactions on the measure becomes stronger. In addition to the qualitative results for the regularity of $p_t^\beta$, we further give the decay estimates on the expectations of $\textup{\textbf{S}}(p^\beta_t)$ and $\textup{\textbf{H}}_p(p^\beta_t)$. 
	\begin{thm}[Proposition \ref{abs}, Theorem \ref{lpregular}]
			Suppose that $\gamma>\frac{3}{2}$, $\mu_0\in \P(\T)$ and the initial empirical measure of the particle model \eqref{model1} satisfies the condition \eqref{cod2'}. If $\beta>\frac{1}{2}\pi K_1^\gamma$ and $\mu_0$ satisfes condition $\textbf{\textup{(J)}}$ , then the fluctuating hydrodynamic limit $p^\beta_t$ is absolutely continuous w.r.t $Leb$ for all $t\geq 0$, $\tilde{\P}-$a.s., and satisfies
			\begin{equation}\label{entineq}
				\tilde{\E}[\textup{\textbf{S}}(p^\beta_t)]\leq e^{-2 C_1^\beta t} C_0+C, \ \ \ \forall \, t\geq 0,
			\end{equation}
			where $ C_0=\ln ||\rho_0||_{\L^\infty}$ and $C=\frac{\beta \pi}{C_1^\beta}$. Furthermore, if $\beta>\frac{p}{2}\pi K_1^\gamma$ for $p>1$, then given any $\mu_0\in \P(\T)$, $p^\beta_t$ is absolutely continuous w.r.t Leb almost surely for any $t>0$ and
		\begin{equation*}\label{introest}
			\tilde{\E}\big[\textup{\textbf{H}}_p(p^\beta_t)\big]\leq \frac{1}{C(1-e^{-2C t})},\ \ \  \forall \, \beta>\frac{p}{2}\pi K_1^\gamma,
		\end{equation*}
	for some positive constants $C$ depending on $p$.
	\end{thm}
 \quad The estimate \eqref{entineq} is not trivial. In fact, for a solution $\mu_t$ to \eqref{withrepul} with a density $\rho_t$, a formal application of It\^o's formula on $\textbf{S}(\mu_t)$ yields
	\begin{equation*}\label{itoformal}
		d\textbf{S}(\mu_t)=\beta \int_0^1\partial_x\rho_t\cdot \H \mu_t dx dt+\frac{K^\gamma_1}{2}\int_0^1 \rho^3_t dx dt + dM_t,
	\end{equation*}
	where $M_t$ is a local martingale. It is not hard to check that 
	\begin{equation*}
		\int_{0}^1\partial_x\rho(x)\H \mu (x) dx=-||\rho||^2_{\dot{\mathcal{H}}_{\frac{1}{2}}}.
	\end{equation*}
	where $\dot{\mathcal{H}}_{s}$ is the homogeneous fractional Sobolev
	space on $\T$ equipped with the norm $$||f||_{\dot{\mathcal{H}}_s}:=\Big(\sum_{k\in \Z , k\neq0}|2\pi k|^{2s}(\hat{f}(k))^2\Big)^\frac{1}{2} .$$ 
	On the other hand, by Sobolev inequality for fractional derivatives on $\T$ \cite{BO13}, we have
	\begin{equation}\label{eqsob}
		C_q||\rho||_{\mathcal{L}^q}\leq ||\rho||_{\mathcal{H}_{\frac{1}{2}}},  \  \ \  \frac{1}{2}<q<\infty,
	\end{equation}
	where $\mathcal{H}_{s}$ is the inhomogeneous fractional Sobolev space  on $\T$ equipped with the norm
	$$||f||_{\mathcal{H}_{s}}:=\Big(\sum_{k\in \Z} (1+4\pi^2 k^2)^s(\hat{f}(k))^2\Big)^\frac{1}{2}.$$
	Since $||f||^2_{\mathcal{H}_{\frac{1}{2}}}\leq ||f||^2_{\dot{\mathcal{H}}_{\frac{1}{2}}}+||f||^2_{\mathcal{L}^2}$, by plugging this into the inequality \eqref{eqsob} one could expect that
	$$\frac{d}{dt}\tilde{\E}[\textbf{S}(\mu_t)]\leq -\beta C_q\tilde{\E}||\rho_t||_{\mathcal{L}^q}^2+\beta\tilde{\E}||\rho_t||^2_{\mathcal{L}^2}+\frac{K^\gamma_1}{2}\tilde{\E}||\rho_t||_{\mathcal{L}^3}^3.$$
	Following this inequality, it seems hard to neutralize the It\^o's correction term $\frac{K^\gamma_1}{2}\E||\rho_t||_{\mathcal{L}^3}^3$ by the regularization effect of repulsive interactions without using informations from the approximation model. Our method is to derive a uniform decay estimate like \eqref{entineq} for a suitable substitution of $L_N(t)$, which shares the same limit with $L_N(t)$ at the same time. In fact, we define the auxiliary empirical measure process associated with $L_N(t)$ by 
	$$\sum_{i=0}^{N}\frac{1}{N|X_N^i-X_N^{i+1}|}\mathbbm{1}_{(X_N^i,X_N^{i+1}]}.$$
	Here, we just denote this measure as $\nu_t^N$ at the moment. $\nu_t^N$ is a local mollification of Dirac measure and shares the limit with $L_N(t)$. In such sense, $\textbf{S}(\nu_t^N)$ can be seen as a microscopic counterpart of $\textbf{S}(p^\beta_t)$. It is easy to see that $S^N(L_N(t))\geq \textbf{S}(\nu_t^N)$, then the estimate \eqref{ugr'} in Theorem \ref{introthm1} gives the decay estimate we want. Formally, one of the key ingredients in proving \eqref{ugr'} is the following inequality 
	\begin{equation}\label{aaa}
		-\int_\T  \partial_x(\H L_N(t))d\nu_t^N\leq C-C'\int_{\T}(\nu_t^{N})^3 dx,
	\end{equation}
  based on which one can expect that the It\^o's correction term can be neutralized if $\beta$ is large enough. See Theorem \ref{introthm1} and Proposition \ref{abs} for details. 
  
   
   \quad We also use the similar argument to derive an estimate on the expectation of R\'enyi entropy $\textbf{H}_p(p_t^\beta)$ in Theorem \ref{lpregular}. The repulsive interactions regularize any initial measure instantaneously to an absolutely continuous measure and keep $\E[\textbf{H}_p(p^\beta_t)]$ uniformly bounded on $t\in [\eps, \infty)$ for any positive number $\eps>0$ when $\beta>\frac{p}{2}\pi K_1^\gamma$. This interesting regularization effect can also be seen in the non-stochastic case, see e.g. \cite{Carrillo2011AMA}. Here, we give a decay estimate. 

\begin{remark}
	The $N$ limit of \eqref{aaa} can be seen as the following inequality 
	$$-\int_{0}^{1}\partial_x(\H \rho) \rho dx\leq C-C' \int_{0}^{1}\rho^3 dx,$$
	for a general measure with density $\rho$. However, we do not know whether this inequality holds for nonnegative function $\rho\in \L^3\cap \dot{\mathcal{H}}_{\frac{1}{2}}$.
\end{remark}

	
	\begin{remark}
		The well-posedness problem of \eqref{withrepul} seems difficult. As a stochastic Fokker-Planck equation, its drift and diffusion coefficients are not Lipschitz w.r.t both measure and space variables; also, unlike the viscous case considered in \cite{fehrman2021}, the Laplacian term in \eqref{withrepul} totally comes from It\^o's correction, and hence no regularization effect from viscosity helps. Nonetheless, we give a new viewpoint to deal with the term $\sqrt{\mu}$ in fluctuation part. More specifically, to renormalize the noise term $\partial_x(\sqrt{\mu} \dot{\xi})$ in the original equation \eqref{onex1}, \cite{fehrman2021} uses a spectral approximation to the white noise and the renormalized equation in $\T$ is given by
		\begin{equation}\label{com1}
			d_t \mu =K\partial_{xx}\mu dt-\partial_x(\mu\Gamma(\mu)) dt-\partial_x(\sqrt{\mu} \sum_{k\in \Z}\frac{1}{|k|^\gamma}e_k ) dW^k_t.
		\end{equation}
		Instead, we adopt a new measure dependent spectral approximation to the noise, which is living in the tangent space $\textbf{T}_{\mu}\P(\T)$. As a comparison, the renormalized equation in this way should formally satisfy
		\begin{equation}\label{com2}
			d_t \mu =K\partial_{xx}\mu dt-\partial_x(\mu\Gamma(\mu)) dt-\partial_x(\mu \sum_{k\in \Z}\frac{1}{|k|^\gamma}e_k(F_{\mu}) ) dW^k_t,
		\end{equation}
		based on the construction of stochastic moving frames. For both two ways of spectral approximations, formal computations yield the consistency of the quadratic variations of the martingale solutions to \eqref{com1} and \eqref{com2} with that of the original Dean-Kawasaki type equation as $\gamma\to 0$. For the case $\beta=0$, \eqref{withrepul} and its associated stochastic moving frames is closely related with a stochastic system of conservation law, which will be specified in subsection \ref{subsec2.1}.
	\end{remark}



	\subsection{Related literature}
	\quad We introduce some other examples of diffusions constructed by methods in various settings, along with some particle approximations. Stuem and von Renesse \cite{von2009} construct an Entropic measure on $\P(\T)$ by using Dirichlet form. The associated diffusion formally corresponds to a Langevin dynamic, where the drift term is the negative gradient of Boltzmann entropy. Their Entropic measure supports on singular measures. In \cite{andres2010}, a particle approximation to the Wasserstein diffusion in \cite{von2009} is constructed with electrostatic interactions from neighbour. Konarovskyi and Von Renesse \cite{konarovskyi2015, konarovskyi2017system, KON24} propose the Arrita flow, Modified Arrita flow (MAF) and their variants for a field of coalescing Brownian motions, which capture the main features of the Brownian motions on $\P(\T)$. In their approximation model, mass carrying Brownian particles interact through a collision mechanism aggregating the mass and scaling the diffusivity in inverse proportional way. \cite{marx2018} gives a regularized approximation to MAF, which keeps similar properties of Wasserstein diffusions while has better regularity. Besides, Cornalba, Shardlow and Zimmer \cite{cornalba2019, cornalba2020} regularize the Dean-Kawasaki model via the second order Langevin dynamic derivation and get well-posedness for a regularized undamped equivalent of \eqref{onex1}. Also, \cite{REN2024104339} construct a Wasserstein diffusion, whose stationary measure is a fully supported Gaussian measure, by taking advantage of an Ornstein-Uhlenbeck (O-U) type Dirichlet form.
	
	\quad On the other hand, we also briefly introduce some related literatures on particle systems with  logarithmic interactions. In random matrix theory \cite{mehta2004random,Anderson_Guionnet_Zeitouni_2009,Fo2010}, the eigenangle distributions of the $\beta-$circular ensembles on $U(N)$ are given by
	
	\begin{equation}\label{Gibbs}
		\frac{1}{Z_N}e^{-\frac{\beta}{N^2}\sum_{l,j=1}^N\ln|e^{i\theta_l}-e^{i\theta_j}|}d\theta.
	\end{equation}
	
	The associated Langevin dynamic are known as the unitary Dyson Brownian motion \cite{dyson1962brownian}, which has generated much interests from various areas in mathematics and physics, including quantum chaotic systems, two-dimensional Yang–Mills theory and free probability theory \cite{biane1997free,voiculescu1998lectures}. Our particle model shares the same drift part with the unitary Dyson Brownian motion, but with a specially designed noise. Moreover, the measure \eqref{Gibbs} coincides with the Gibbs distribution of Coulomb gas on $\S^1$, and \eqref{LIE} is just the sum of logarithmic interactions between particles confined to $\S^1$. There are intensive studies on the log and Coulomb gas model (see e.g. \cite{forrester2010log,bourgade2014universality, ES2012,ES2013, LS2017,guionnet2019rigidity}), we only mention a few aspects due to the limitation of the scope of this paper. For the classical Coulomb gas model on the plane $\R^2$, the macroscopic behaviour of this system is totally determined by the minimizer of the energy with an external potential:
	$$\beta\int_{\R^2}\int_{\R^2}\ln|x-y| d\mu(x)d\mu(y)+\beta\int_{\R^2}V d\mu,$$
	which is independent of $\beta$, while the microscopic behaviour seems more complicated. As shown by \cite{LS2017}, one can see a tendency to disorder states due to the effect of entropy, when the configurations
	are dilated by $N^{-\frac{1}{2}}$. In particular, as $\beta \to 0$, the entropy dominates and
	the	system converges to a Poisson process \cite{Lebl2015LogarithmicCA,LS2017}, which means the interaction energy no longer influences the microscopic configurations. In our stochastic model on $\T$, the competition between the repulsive interaction energy and the Dean-Kawasaki type noise follows a similar pattern as the temperature varies, which can be seen from the pathwise regularity of $p_t^\beta$. More precisely, as $\beta$ gets smaller, the pathwise regularity of $p_t^\beta$ gets worse.

	\subsection{Organization}
	
	\quad In Section \ref{sec1.2}, we firstly formally derive the form of stochastic moving frames by using stochastic parallel translations in $\P_2(\R)$. Then we introduce the definition of the noise $\xi^\beta_{\mu}$ and the solution to the martingale problem of \eqref{withrepul}. As a remark, we show the formal consistency of solutions to the martingale problem of \eqref{withrepul} with solutions to the martingale problem of \eqref{onex1} in the case $\Gamma=0$. 
	
	\quad In Section \ref{sec3}, we prove Theorem \ref{introthm1} and Theorem \ref{ratecon}. Also, we prove in Proposition \ref{prop2} that given admissible $\gamma$ and $\beta$, the limit distribution $P$ is independent of the choice of the initial empirical measure $\{L_N(0)\}$. Finally, we prove in Proposition \ref{prop5.1} that the limit distributions satisfy $\W^2_2(P_{\beta}, P_{\beta'})\leq C|\beta-\beta'|$ for any $\beta, \beta'>\lambda_0^\gamma$.  

	\quad In Section \ref{sec6}, We firstly prove in Proposition \ref{noatom} that the fluctuating hydrodynamic limit $p^\beta_t$ is non-atomic for all $t\in [0,\infty )$ almost surely when $\beta>\lambda_0^\gamma$, and is a solution to the martingale problem of \eqref{withrepul} (Theorem \ref{maintheorem}). Then we prove Proposition \ref{abs} and Theorem \ref{lpregular}. 
	
	
	\vskip4mm
	\section{Formulations}\label{sec1.2}
	\quad In this section, we will firstly introduce the formal construction of stochastic moving frames without consideration of well-posedness and regularity. The main purpose is to explain why we propose the noise part of the particle model \eqref{model1}. Then we will introduce the definition of the solution to the martingale problem of \eqref{withrepul}.
	
	\subsection{Stochastic moving frames}\label{subsec2.1}	
	\quad Starting from Lagrangian viewpoint, we heuristically introduce the intrinsic construction of Wasserstein diffusion by constructing the corresponding stochastic flows. We will construct a so called stochastic moving frame $\displaystyle \{v_k(t)\}_{k\in \N}$ such that the push forward measure process $\mu_t^N$ induced by the stochastic flow satisfying
	\begin{equation*}\label{flow}
		dX^N_t=\sum_{k=-N}^{N}v_k(t,X^N_t) \circ dW_t^k,
	\end{equation*}
can formally solve the martingale problem of a Dean-Kawasaki type equation. 
	
	\quad  Recall that on a Riemannian manifold $M^d$ equipped with Levi-Civita connection $\bm{\nabla}$, Brownian motions are projections of stochastic horizontal curves $U_t$ in the orthonormal frame bundle $\mathcal{O}(M^d)$ down to the manifold, satisfying
	\begin{equation}\label{insde}
		dU_t=\sum_{i=1}^d H_i(U_t)\circ dW_t^i,
	\end{equation}
	where $H_i$ are fundamental horizontal vector fields, with a initial frame $U_0=(x_0, \{\textbf{v}_k(0)\}_{k=1}^d)\in F(M^d)$. SDE \eqref{insde} is equivalent to the following stochastic equations
	\begin{align}
		&dX_t=\sum_{j=1}^{d}\textbf{v}_j(t,X_t)\circ dW_t^j \label{coor1}\\
		&\bm{\nabla}_{\circ dX_t}\textbf{v}_k(t)=0,  \ \ \ \forall \, k=1,\ldots,d  . \label{coor2}
	\end{align}
	\eqref{coor1} gives dynamics on $M^d$ while \eqref{coor2} gives dynamics of the vector basis. To realize this idea to $M=\P_2(\R)$ for an initial value $\displaystyle \big(\mu_0, \{e_k(F_{\mu_0})\}_{k\in \N}\big)$, we adopt Lagranian viewpoint to interpret dynamics on $\P_2(\R)$ as a push forward measure process $\mu_t^N$ induced by a stochastic flow $X_t^N$, which is driven by some stochastic moving vector fields $\displaystyle \{v_k(t)\}_{k=-N}^N$ in $\L^2(\mu_t^N)$, i.e. 
	\begin{equation*}\label{flow'}
		dX^N_t=\sum_{k=-N}^{N}v_k(t,X^N_t) \circ dW_t^k,
	\end{equation*} 
	with $v_k(0)=e_k(F_{\mu_0})$. Note that $\{e_k(F_{\mu})\}_{k\in \N}$ is an O.N.B. in $\L^2(\mu_0)$ if $(F_{\mu})_{\#}\mu=Leb$, since 
	$$\int_{\R} e_i(F_{\mu})e_j(F_{\mu}) d\mu=\delta_{i,j}.$$
	On the other hand, $\displaystyle \{v_k(t)\}_{k=-N}^N$ are determined by stochastic parallel translations of $\displaystyle \{e_k(F_{\mu_0})\}_{k=-N}^N$ along $\mu^N_t$ on $\P_2(\R)$, equipped with Lott's Levi-Civita connection \cite{lot2006som}.
	In \cite{fan2021}, it is shown that when $X_t^N$ is a stochastic diffeomorphic flow on a compact manifold, stochastic parallel transport vector field $u_t$ of $u_0$ along $\mu^N_t$ should satisfy
	\begin{equation}\label{stopara}
		d_t u_t+\sum_{k=-N}^{N}\Pi_{\mu^N_t}\big(\nabla_{v_k(t)} u_t\big)\circ dW^k_t=0,
	\end{equation}
	where $\Pi_{\mu}$ is defined as the projection operator from $\L^2(d\mu)$ onto the gradient type of vector field. In $\R$, $\Pi_\mu$ is trivial, i.e. $\Pi_{\mu}\big(w\big)=w$.  Then \eqref{stopara} actually becomes
	\begin{equation}\label{transport}
		d_t u_t+\partial_xu_t\circ dX^N_t=0.
	\end{equation} 
	Therefore the stochastic horizontal curve $\big(\mu^N_t, \{v_k(t)\}_{k\in \N}\big)$ in $\mathcal{O}(\P(\R))$ should formally satisfy
	\begin{equation}\label{spt}
		\left\{
		\begin{aligned}
			&dX^N_t=\sum_{k=-N}^{N}v_k(X^N_t)\circ dW_t^i,\\
			&\mu^N_t=(X^N_t)_{\#}\mu_0,\\
			&d_t v_k(t)+\partial_xv_k(t)\circ dX^N_t=0, \ \ \ \ \ \forall \, k\in \N,
		\end{aligned}
		\right.
	\end{equation}
	with the initial values $v_k(0)=e_k(F_{\mu_0})$ for $k\in \N$. Note that in \eqref{transport}, the stochastic parallel translations $u_t$ of $u_0$ along $\mu^N_t$ is the unique strong solution to \eqref{transport} satisfying
	\begin{equation}\label{ins}
		u(t,x)=u(0,(X^N_t)^{-1}(x)),
	\end{equation}
	by the characteristic method if $X^N_t$ is a diffeomorphic flow. According to this representation formula \eqref{ins} of stochastic parallel translations, we construct the stochastic moving frame $\{v_k(t), \, k\in \N\}$ as $\displaystyle v_k(t,x)=e_k\big(0, ( X_t^N)^{-1}(x)\big)$. Under the assumption that $X^N_t$ keeps order, we further have 
	\begin{equation*}\label{eqXt}
		(X^N_t)^{-1}=F^{-1}_{\mu_0}\circ F_{\mu^N_t}.
	\end{equation*} 
	Thus, we construct $\{v_k(t)\}_{k\in\N}$ as 
	\begin{equation}\label{ek}
		v_k(t,x):=e_k( F_{\mu^N_t}(x)), \ \ \ \forall \, x\in \R.
	\end{equation}
	Assuming that $\mu_t^N$ has a density $\rho_t^N $, we immediately have
	\begin{equation*}\label{partiale}
		\partial_x v_k(t,x)=\partial_x e_k(F_{\mu^N_t}(x))\cdot \rho^N_t(x).
	\end{equation*}

	\quad We will formally check that $\{\mu^N_t, \{e_k(F_{\mu_t^N})\}_{k\in \N}\}$ satisfies \eqref{spt}, i.e. $\{\mu^N_t, \{e_k(F_{\mu_t^N})\}_{k\in \N}\}$  is indeed a stochastic horizontal curve. From now on, we treat the stochastic flow $X_t^N$ as a diffeomorphic flow so that we can proceed the formal derivations, although this assumption seems impossible at the moment. Note that 
	\begin{equation}\label{iden2}
		\sum_{k=-N}^{N}\partial_xe_k( F_{\mu^N_t}(X_t^N)) \rho_t^N(X_t^N)e_k( F_{\mu^N_t}(X_t^N)) =0,
	\end{equation}
	the It\^o's equation of $X_t^N$ reads
	\begin{equation*}\label{Xt}
		dX^N_t=\sum_{k=-N}^{N}v_k(t,X^N_t) dW_t^i.
	\end{equation*}
	Next, we check that for each $k\in \N$, $v_k(t)$ satisfies
	\begin{equation}\label{spt2}
		d v_k(t)+\partial_x v_k(t)\circ\Big(\sum_{j=-N}^{N} v_j(t) dW^j_t\Big)=0.
	\end{equation}
	Indeed, a formal application of It\^o's formula on $X_t^N$ shows that, for $\forall \, f\in \C_c^\infty(\R)$,
	\begin{equation*}
		\begin{split}
			d_tf(X_t^N)&=\sum_{k=-N}^{N} f'(X_t^N)e_k( F_{\mu^N_t}(X_t^N)) dW_t^k+\frac{1}{2}f''(X_t^N)\sum_{k=-N}^{N}e^2_k( F_{\mu^N_t}(X_t^N)) dt.
		\end{split}
	\end{equation*}
	Since
	\begin{equation}\label{iden1}
		\frac{1}{2}\sum_{k=-N}^{N}e^2_k( F_{\mu^N_t}(X_t^N)) =N+\frac{1}{2},
	\end{equation}
	we have
	$$d_tf(X_t^N)=\sum_{k=-N}^{N} f'(X_t^N)e_k( F_{\mu^N_t}(X_t^N)) dW_t^k+(N+\frac{1}{2})f''(X_t^N) dt.$$
	Note that in sense of distributional derivative, we have
	$$\big(\mathbbm{1}_{(-\infty,x]}\big)'(y)=-\delta_{x}(y).$$
	Hence we obtain, by formally using stochastic Fubini theorem and integration by parts, that
	\begin{equation*}
		dF_{\mu^N_t}=(N+\frac{1}{2})\partial_x \rho^N_t dt -\sum_{k=-N}^{N} \rho^N_te_k(F_{\mu_t^N})dW_t^k.
	\end{equation*}
	It yields again by It\^o's formula and \eqref{iden1} that
	\begin{equation*}
		\begin{split}
			&d_te_k(F_{\mu_t^N})\\
			&=-\partial_xe_k(F_{\mu^N_t})\cdot \rho^N_t\sum_{j=-N}^N e_j(F_{\mu_t^N})dW_t^j+(N+\frac{1}{2})\partial_xe_k(F_{\mu^N_t})\partial_x\rho^N_t dt\\
			&+(N+\frac{1}{2})\partial_x^2 e_k(F_{\mu^N_t})\big(\rho^N_t\big)^2dt.
		\end{split}
	\end{equation*}
	On the other hand, we have for $\forall \, k\in \N$,
	\begin{equation*}
		\begin{split}
			&\partial_x v_k(t)\circ\Big(\sum_{j=-N}^{N} v_j(t) dW^j_t\Big)\\
			&=\partial_x e_k(F_{\mu^N_t})\rho^N_t\cdot\sum_{j=-N}^{N} e_j(F_{\mu_t^N}) dW^j_t-\frac{1}{2}\sum_{j=-N}^{N}\partial^2_xe_k(F_{\mu^N_t})\big(\rho^N_t\big)^2e^2_j(F_{\mu_t^N}) dt\\
			&-\frac{1}{2}\sum_{j=-N}^{N}\partial_x e_k(F_{\mu^N_t})\partial_x\big(\rho^N_te_j(F_{\mu_t^N})\big)e_j(F_{\mu_t^N}) dt.
		\end{split}
	\end{equation*}
	Now, using \eqref{iden1} and \eqref{iden2}, we have formally checked \eqref{spt2}. Therefore, we obtain that $U_t=\{\mu^N_t, \{e_k(F_{\mu_t^N})\}_{k\in \N}\}$ is a horizontal curve satisfying \eqref{spt}. As $\mu_t^N$ is the projection of $U_t$ down to $\P_2(\R)$,  we may say $\mu^N_t$ is a degenerated Brownian motion on $\P_2(\R)$. 
	
	\quad We also give a formal verification that $\mu^N_t$ indeed shares the same fluctuation with the martingale solution to the Dean-Kawasaki type equation. We formally obtain, by applying \eqref{ek} into \eqref{spt}, that
	\begin{equation}\label{DKNeq}
		\partial_t \mu^N_t= (N+\frac{1}{2})\partial_x^2\mu_t^N-\partial_x\big(\mu^N_t \sum_{j=-N}^N e_j(F_{\mu^N_t}) \dot{W}^j_t\big)
	\end{equation}
	whose martingale solution satisfies, for $\forall \, f\in \C^\infty_c(\R)$, 
	$$d_t\langle f,\mu^N_t\rangle=\sum_{k=-N}^N \langle f', e_k( F_{\mu^N_t}) \rangle_{\mu^N_t} dW_t^k +(N+\frac{1}{2})\langle f'', \mu_t^N\rangle dt, $$
	with quadratic variation process
	$$< M_t(f)>=\int_{0}^{t}\sum_{k=-N}^{N}\big(\int_0^1f'e_k(F_{\mu^N_s})d\mu^N_s\big)^2ds.$$
	If $(F_{\mu^N_t})_{\#}\mu^N_t=Leb$, we have
	$$<M^N_t(f)>=\int_{0}^{t}\sum_{k=-N}^N \Big(\widehat{f'(F_{\mu_s^N}^{-1})}(k)\Big)^2ds.$$
	Let $N\to \infty$, by Parseval's identity, we formally obtain
	\begin{equation}\label{parse}
		\lim\limits_{N\to \infty}< M^N_t(f)> = \int_{0}^{t}||f'(F_{\mu_s}^{-1})||_{\L^2(Leb)}^2 ds.
	\end{equation}
	Assume that $(F_{\mu})_{\#}\mu =Leb$, we immediately have
	\begin{equation}\label{eqtran}
		\int_0^1 \Big| f'(F_{\mu_s}^{-1}(x))\Big|^2 dx =\int_{0}^{1}|f'|^2 d\mu_s.
	\end{equation}
	Applying \eqref{eqtran} into \eqref{parse} gives
	$$	\lim\limits_{N\to \infty}< M^N_t(f)> =\int_{0}^{t}||f'||^2_{\L^2(\mu_s)} ds.$$

	\quad Heuristically, we yield both by the theory of intrinsic construction of Brownian motion and the consistency of quadratic variation that the push forward measure process $\mu_t$ induced by stochastic flow $X_t$ satisfying 
	\begin{equation*}
		dX_t=\sum_{k=-\infty}^{\infty} e_k(F_{\mu_t}(X_t)) dW_t^k, 
	\end{equation*}
	is a Wasserstein diffusion and has a noise part whose quadratic variation process coincides with \eqref{charac}. However, according to \eqref{DKNeq}, we need to renormalize so that the drift part of $\mu_t$ is well posed as $N\to \infty$. One way to renormalize is to replace by a spectral approximation. In fact, we consider the stochastic flow $X^\gamma_t$ with a coefficient $\gamma>\frac{1}{2}$ governed by
	\begin{equation}\label{noise2}
		dX^\gamma_t=\sum_{k=-\infty}^{\infty} \frac{1}{|k|^\gamma}e_k(F_{\mu^\gamma_t}(X^\gamma_t)) dW_t^k,
	\end{equation}
	where $\mu^\gamma_t=(X^\gamma_t)_{\#} \mu_0$.

	\begin{remark}\label{rem2.1}
		Adding weights $a_k=\frac{1}{|k|^\gamma}$ on $k-$th mode to \eqref{spt}, we obtain an equation of stochastic $"(N+1)\times (N+1)"$ consevation law for $(\mu^N, v_1, \ldots, v_N)$:
		\begin{equation}\label{scl1}
			\left\{
			\begin{aligned}
				d_t \mu^N&=-\sum_{j=-N}^{N}\frac{1}{|j|^{\gamma}}\partial_x(\mu^N v_j)\circ dW_t^j,\\
				d_t (\mu^N v_i) &=- \sum_{j=-N}^{N}\frac{1}{|j|^{\gamma}}\partial_x(\mu^N v_i v_j)\circ dW_t^j, \ \ \forall \, i=1,\ldots, N
			\end{aligned}
			\right.
		\end{equation} 
		with a initial value $(\mu_0, v^0_1\circ F_{\mu_0}, \ldots, v^0_N\circ F_{\mu_0})$. According to our analysis above, once we find a weak solution $\mu_t^N$ to \eqref{DKNeq} with weights $\frac{1}{|k|^\gamma}$, then $\{\mu_t^N, v^0_1\circ F_{\mu_t^N}, \ldots, v^0_N\circ F_{\mu_t^N}\}$ is also a weak solution to the system \eqref{scl1}. This construction of solutions can also be generalized to the case $N\to \infty$.
	\end{remark}
	
	

	\begin{remark}
		According to Otto's Riemannian structure on the Wasserstein space $\P(\T)$, a tangent vector in $\textbf{T}_{\mu}\P(\T)$ should be a gradient type vector field in $\L^2(\mu)$. However, the stochastic moving vector $e_k(F_\mu)$ is not a gradient type function on $\T$. In such sense, $\sum_{k\in \Z}\frac{1}{|k|^{\gamma}}e_k(F_{\mu_t})dB_t^k$ is not exactly $d\textbf{B}_t^Q$. As mentioned in Subsection \ref{subsec1.1}, it is hard to directly construct a gradient type moving frame on $\P(\T)$.
	\end{remark}

	\vskip2mm

	\subsection{Definitions of the noise $\dot{\xi}_{\mu}^\gamma$ and the martingale problem }\label{sec2}
	
	\quad In order to describe the fluctuating hydrodynamic limits of the particle models \eqref{model1}, we give the definition of the noise $\dot{\xi}_{\mu}^\gamma$ and the martingale problem of \eqref{withrepul} for initial measure $ \mu_0$. Based on \eqref{noise2}, we define $\xi_{\mu}^{\gamma}$ as 
	\begin{equation*}
		\xi_{\mu}^\gamma(t,x)=\sum_{k=-\infty}^{+\infty}\frac{1}{|k|^{\gamma}}e_k(F_{\mu})W_t^k.
	\end{equation*}
	
	Given $\mu\in \P(\T)$, the law of $\xi_{\mu}^{\gamma}$ can be determined by the correlation kernel
	\begin{equation*}
		\bar{Q}_{\mu}^{\gamma}(x,y):=\sum_{i=-\infty}^{\infty}\frac{1}{|i|^{2\gamma}}e_i(F_\mu(x))e_i(F_\mu(y))=\sum_{k=0}^{\infty}\frac{2}{k^{2\gamma}}\cos(2\pi k(F_{\mu}(x)-F_{\mu}(y))),
	\end{equation*}
	because, for $\forall \, f \in \C^2(\T)$, the quadratic variation process of the martingale $$M_t(f)=\int_0^1f(x)\xi_{\mu}^\gamma(t,x)d\mu(x) $$
	is determined by
	\begin{equation*}
		\frac{d}{dt}<M_t(f)>=\int_{0}^1\int_{0}^1f(x)f(y)\bar{Q}_{\mu}^{\gamma}(x,y)d\mu(x) d\mu(y).
	\end{equation*}

	$\xi_{\mu}^\gamma$ is a $(\bar{Q}_{\mu}^{\gamma})^{\frac{1}{2}}-$Wiener process in $\L^2(\mu)$ and satisfies 
	\begin{equation*}
		\E[\xi^\gamma_\mu(t, x)]=0; \quad \E[\xi^\gamma_\mu(t,x)\xi^\gamma_\mu(s,y)]=(t\land s) \cdot\bar{Q}_{\mu}^{\gamma}(x,y).
	\end{equation*}
	
	We denote its time derivative in distribution as $\dot{\xi}^\gamma_\mu$ . It can be proved by Doob's inequality that $\dot{\xi}^\gamma_\mu (t,\cdot)\in \L^2(\mu) $ a.s. .

	\quad Note that if $\mu$ does not have atoms, then 
	\begin{equation*}
		Leb=(F_{\mu})_{\#}\mu.
	\end{equation*} 
	By the definition of the pushforward measure, 
	\begin{equation*}
		\int_{0}^{1}e_i(F_{\mu})e_j(F_{\mu})d\mu=\int_0^1e_ie_jdx=\delta_{ij},
	\end{equation*}
	thus $\{e_k(F_{\mu})\}_{k\in \Z}$ constitute an O.N.B. of $\L^2(d\mu)$. In such case, we have for $\forall \, k\in \Z$,
	\begin{equation*}
		\begin{split}
			\int_{0}^1\bar{Q}_{\mu}^{\gamma}(x,y)e_k(F_\mu(y)) d\mu(y)&=\int_{0}^1\Big(\sum_{i=-\infty}^{\infty}\frac{1}{|i|^{2\gamma}}e_i(F_\mu(x))e_i(F_\mu(y))\Big)e_k(F_\mu(y)) d\mu(y)\\
			&=\int_{0}^1\Big(\sum_{i=-\infty}^{\infty}\frac{1}{|i|^{2\gamma}}e_i(F_\mu(x))e_i(F_\mu(y))\Big)e_k(F_\mu(y)) dy\\
			&=\frac{1}{|k|^{2\gamma}}e_k(F_\mu(x)) .
		\end{split}
	\end{equation*}
	
	Hence, $\dot{\xi}_\mu^\gamma$ is indeed a spectral approximation to the cylindrical noise in $\L^2(\mu)$. According to the correlation kernel of the noise, we give the definition of solutions to the martingale problem of \eqref{withrepul}.

	
	\begin{definition}\label{MPs}
		We say a continuous $\P(\T)$-valued process $\mu_t$ is a solution to the martingale problem $(MP)^{K,\beta,\gamma}_{\mu_0}$ of \eqref{withrepul}, if there exists a filtered probability space $(\Omega, \mathcal{F}, \{\mathcal{F}_t\}_{t\geq 0}, \P)$ such that 
	for $\forall \phi \in \C^{2}(\T)$,
		\begin{equation*}
			M_t(\phi):=\langle\mu_t,\phi \rangle-\langle \mu_0,\phi\rangle -\frac{\beta}{2}\int_0^t \int_{0}^{1}\int_{0}^{1}\frac{\phi'(x)-\phi'(y)}{\tan(\pi(x-y))}d\mu_s(x)d\mu_s(y)ds-K\int_0^t \langle \mu_s,\phi''\rangle d s
		\end{equation*}
		is a $\F_t-$adapted martingale, whose quadratic variation process is given by 
		\begin{equation*}
			<M_t(\phi)>=\int_0^t Q_{\mu_s}^{\gamma} (\phi,\phi)d s ,
		\end{equation*}
	where $\displaystyle Q_{\mu_s}^{\gamma} (\phi,\phi)=\int_{0}^{1}\int_{0}^{1}\phi'(x)\phi'(y)\bar{Q}_{\mu_s}^{\gamma}(x,y) d\mu_s(x)d\mu_s(y)$.
	\end{definition}

	\vskip 2mm
	\quad Note that when $\mu_t$ does not have atoms, we have $Leb=(F_\mu)_{\#}\mu$, which yields
	\begin{equation*}\label{repre2.2}
		Q_{\mu_s}^{\gamma} (\phi, \phi)=\int_{0}^1\int_{0}^1\phi'(G_{\mu_s}(x))\phi'(G_{\mu_s}(y))\Big(\sum_{k=0}^{\infty}\frac{2}{k^{2\gamma}}\cos(2\pi k(x-y))\Big)dxdy,
	\end{equation*}
	where $G_{\mu_s}$ is the quantile function of $\mu_s$ satisfying $G_{\mu_s}\circ F_{\mu_s}(x)= x$ . Due to
	$$\displaystyle |\sum_{k=0}^{\infty}\frac{2}{k^{2\gamma}}\cos(2\pi k(x-y))|<2K_2^\gamma,$$  
	we have
	\begin{equation}\label{qbeta}
		\begin{split}
			Q^{\gamma}_{\mu_s} (\phi, \phi)&=\int_{0}^1\int_{0}^1\phi'(G_{\mu_s}(x))\phi'(G_{\mu_s}(y))\Big(\sum_{k=-\infty}^{\infty}\frac{1}{|k|^{2\gamma}}e_k(x)e_k(y)\Big)dxdy \\
			&=\sum_{k=-\infty}^{+\infty}\frac{1}{|k|^{2\gamma}}\int_{0}^1\phi'(G_{\mu_s}(x))e_k(x)dx\int_{0}^1\phi'(G_{\mu_s}(y))e_k(y)dy\\
			&=\sum_{k=-\infty}^{+\infty}\frac{1}{k^{2\gamma}}|\widehat{\phi'(G_{\mu_s})}(k)|^2,
		\end{split}
	\end{equation}
	where the fourier coefficients are defined as 
	\begin{equation*}
		\begin{split}
			\widehat{f}(k)&=\sqrt{2}\int_{0}^{1}f(x)\sin(2\pi k x) dx , \quad k=1,2,\ldots ; \\
			\widehat{f}(0)&=\int_{0}^{1}f(x)dx ;\\
			\widehat{f}(k)&=\sqrt{2}\int_{0}^{1}f(x)\cos(2\pi k x) dx , \quad k=-1,-2,\ldots .
		\end{split}
	\end{equation*}
	\vskip1.5mm
	\begin{remark}[Formal consistency with \eqref{onex1}]
		According to \eqref{qbeta}, when $\gamma=0$, the quadratic variation of $M_t(\phi)$ formally becomes 
		\begin{equation*}
			<M_t(\phi)>=\int_0^t||\phi'(G_{\mu_s}(x))||^2_{\L^2} d s =\int_{0}^{t}||\phi'||^2_{\L^2(\mu_s)} ds.
		\end{equation*}
	\end{remark}

	\section{Particle model and rate of convergence}\label{sec3}
	\quad In this section, we will mainly prove the well-posedness of particle model on $\T$ and the convergence rate of the associated empirical measure process. For convenience of reading, we write the particle model again. For $N>0$ , the particle model satisfies
	\begin{equation}\label{Napprox}
		d X_N^i(t)=\frac{\beta}{N}\sum^N_{j=1, j \neq i}\cot\big(\pi(X_N^i(t)-X_N^j(t))\big) d t +  \sum_{k=-N}^{N}\frac{1}{|k|^{\gamma}}e_k(\frac{i}{N})dW_t^k  , 
	\end{equation}
	for $i=1,\ldots,N$, with the initial value $X_N(0)\in \Delta_N$. The It\^o's correction term equals 0 because the diffusion coefficients are constants. $\gamma, \beta$ are the model's macroscopic parameters. For simplicity of notations, we denote $X_N^i(t)$ as $X_t^i$ in the proof of Theorem \ref{introthm1}. Let
	\begin{equation*}
		K_1^{N,\gamma}=\sum_{j=1}^{N}\frac{4\pi^2}{j^{2\gamma-2}} ;  \ \ \ 
		K_2^{N,\gamma}=\sum_{j=1}^{N}\frac{1}{j^{2\gamma}} ;
	\end{equation*} 
and
	\begin{equation*}
		\begin{split}
			a_m&=\sum_{i=1}^{N-m}\left|\frac{1}{\sin(\pi(X_t^i-X_t^{i+m}))}\right|^2 ;\\
			b_m&=\sum_{i=1}^{N-m}\left|\frac{1}{\tan(\pi(X_t^i-X_t^{i+m}))}\right|^2 .
		\end{split}
	\end{equation*}
	\subsection{Proof of Theorem \ref{introthm1}}
	\quad Before going to the proof, we give an estimate on $b_m+b_{N-m}$, which quantifies how little $b_m+b_{N-m}$ is compared to $b_1+b_{N-1}$.
	\begin{lemma}\label{lem3.1}
		We have for $\forall \, m\in \{1,2,\ldots, [\frac{N}{2}]\}$,  
		\begin{equation}\label{est2}
			b_m+b_{N-m}\leq \left(\frac{1}{m^2}+\frac{1}{(N-m)^2}\right) (b_1+b_{N-1}) +\frac{N}{m^2}.
		\end{equation}
	\end{lemma}
	\begin{proof}
		 If $0<X_t^{i+M}-X_t^i<\frac{1}{2}$, then due to convexity of the function $\frac{1}{x^2}$ , we have
		\begin{equation}\label{estconv}
			\begin{split}
				\frac{1}{|X_t^i-X_t^{i+M}|^2}&=\frac{1}{|\sum_{l=1}^{M}X_t^{i+l-1}-X_t^{i+l}|^2}\\
				&\leq \frac{1}{M^3}	\left(\sum_{l=1}^{M}\frac{1}{|X_t^{i+l-1}-X_t^{i+l}|^2}\right).
			\end{split}
		\end{equation}
	Note that for $x\in (0,\frac{1}{2})$, we have
	\begin{equation}\label{compa}
		\frac{1}{\tan^2(x)}<\frac{1}{x^2} \leq \frac{1}{\sin^2(x)} =\frac{1}{\tan^2(x)}+1.
	\end{equation}
It then follows from \eqref{estconv} that 
\begin{equation}\label{tan1}
	\frac{1}{\tan^2(\pi(X_t^{i+M}-X_t^i))}<\frac{1}{M^3}\left(\sum_{l=1}^{M}\frac{1}{\tan^2(\pi(X_t^{i+l}-X_t^{i+l-1}))}\right)+\frac{1}{M^2}.
\end{equation}
	If $\frac{1}{2}<X_t^{i+M}-X_t^i<1$, we let $x_t^{i}=X_t^i+1$ and $x_t^{i+M}=X_t^{i+M}$. We have $-\frac{1}{2}<x^{i+M}_t-x_t^i<0$. Using the same argument as above, we obtain
	\begin{equation}\label{tan2}
		\frac{1}{\tan^2(\pi(X_t^{i+M}-X_t^i))}<\frac{1}{(N-M)^3}\sum_{l=1}^{N-M}\frac{1}{\tan^2(\pi(X_t^{i+M+l}-X_t^{i+M+l-1}))}+\frac{1}{(N-M)^2},
	\end{equation}
We denote $X_t^{k+N}=X_t^k$ in the sequel for simplicity of notation. For $i\in {1,2, \ldots, N}$, let 
	\begin{equation*}
	\begin{split}
		B^i_{+}&=\left\{j: 0<|X_N^i-X_N^j|\leq\frac{1}{2},\  j\neq i\right\} \\
		B^i_{-}&=\left\{j: \frac{1}{2}<|X_N^j-X_N^i|\leq1,\  j\neq i\right\} .
	\end{split}
\end{equation*}
Then combined with \eqref{tan1} and \eqref{tan2}, we obtain
\begin{equation*}
		\begin{split}
		b_m+b_{N-m}=&\sum_{i=1}^{N}\frac{1}{\tan^2(\pi|X_t^{i+m}-X_t^i|)}\\
		=&\left(\sum_{ i+m\in B^i_{+}}+\sum_{ i+m\in B^i_{-}}\right)\frac{1}{\tan^2(\pi|X_t^{i+m}-X_t^i|)}\\
		=&\sum_{ i+m\in B^i_{+}} \left(\frac{1}{m^3}\left(\sum_{l=1}^{m}\frac{1}{\tan^2(\pi(X_t^{i+l}-X_t^{i+l-1}))}\right)+\frac{1}{m^2}\right)\\
		+&\sum_{ i+m\in B^i_{-}}\left(\frac{1}{(N-m)^3}\left(\sum_{l=1}^{N-m}\frac{1}{\tan^2(\pi(X_t^{i+m+l}-X_t^{i+m+l-1}))}\right)+\frac{1}{(N-m)^2}\right)\\
		\leq &\frac{N}{m^2}+S_+ +S_-,
	\end{split}
\end{equation*}
where 
$$S_+=\sum_{ i+m\in B^i_{+}} \frac{1}{m^3}\left(\sum_{l=1}^{m}\frac{1}{\tan^2(\pi(X_t^{i+l}-X_t^{i+l-1}))}\right),$$
$$S_-=\sum_{ i+m\in B^i_{-}}\frac{1}{(N-m)^3}\left(\sum_{l=1}^{N-m}\frac{1}{\tan^2(\pi(X_t^{i+m+l}-X_t^{i+m+l-1}))}\right).$$
Since each $\tan^2(\pi(X_t^{i+1}-X_t^{i}))$ can be counted at most $m$ times in $S_+$ and $N-m$ times in $S_-$, it yields
		\begin{equation*}
			b_m+b_{N-m}\leq \frac{N}{m^2}+\sum_{i=0}^{N}\tan^2(\pi(X_t^{i+1}-X_t^{i})) \left(\frac{1}{m^2}+\frac{1}{(N-m)^2}\right).
		\end{equation*}
	
	We finish the proof.
	\end{proof}
	\quad It is worth mentioning that the main difficulty in proving that $\textup{\textbf{W}}(L_N(t))-Ct$ is a supermartingale for some constant $C$ is to obtain the following estimate
	\begin{equation*}
		\begin{split}
			&\sum_{l,j=1, \, l\neq j}^{N}\sum_{k=-N}^{N}\frac{1}{|k|^{2\gamma-2}}\left|\frac{e_k(\frac{l}{N})-e_k(\frac{j}{N})}{\tan(\pi(X_t^l-X_t^j))}\right|^2\\
			&\leq \frac{C}{2N}\sum_{l,j=1, \, l\neq j}^{N}\frac{1}{|\tan(\pi(X_t^l-X_t^{j}))|^2}+o( N^2).
		\end{split}
	\end{equation*}
	Our strategy of proof starts from decomposing the sum on the left side according to different scales of spacing between particles.

	\begin{proof}[\textbf{Proof of Theorem \ref{introthm1}}]
		We follow the classical argument in \cite{rog1993}, \cite{Guionnet} and \cite{li2020}. We construct the truncated process. For $R>100$, let $\phi_R(x)$ be a $\C^2(\R)$ function which satisfies $\phi_R(x)=\cot(\pi x)$ for $x\in( -1+\frac{1}{R},-\frac{1}{R})\cup(\frac{1}{R},1-\frac{1}{R})$. Then the following SDE
		\begin{equation*}
			d X_{R,N}^{i}(t)=\frac{\beta}{N}\sum^N_{j=1, j \neq i}\phi_R(X_{R,N}^{i}(t)-X_{R,N}^{j}(t)) d t + \sum_{k=-[N^c]}^{[N^c]}\frac{1}{|k|^{\gamma}}e_k(\frac{i}{N})dW_t^k,
		\end{equation*}
		with initial value $X_{R,N}^{i}(0)=x^i$ for $1\leq i\leq N$, has a unique strong solution $X_{R,N}(t)$. Let
		\begin{equation*}
			\begin{split}
				\tau_{R}:=\inf\left\{t:\min_{l\neq j}|e^{2\pi iX_{R,N}^{l}(t)}-e^{2\pi iX_{R,N}^{j}(t)}|\leq R^{-1}\right\} .
			\end{split}
		\end{equation*}
		Then $\tau_R$ is monotone increasing in $R$ and $X_{R,N}(t)=X_{R',N}(t)$  for all $t\leq \tau_R$ and $R<R'$. 
		We construct $X_N(t)=X_{R,N}(t)$ on $t\in [0,\tau_R)$. According to the similar argument in \cite{Guionnet}, to prove $X_N(t)$ is the unique solution to \eqref{Napprox}, we only need to show that $(X_t)_{t\in[0,\infty)}$ never collide and $|X_N^1(t)-X_N^N(t)|<1$. Consider the Lyapunov function, as in \cite{rog1993,li2020}, $F(x_1,...,x_N)=-\frac{1}{N^2}\sum_{l\neq j}\ln|e^{2\pi ix_l}-e^{2\pi ix_j}|$. By It\^o's formula, we have
		\begin{equation*}
			\begin{split}
				&	d F(X_t^1,...X_t^N)\\
				&=-\frac{\beta \pi}{N^{2}}\sum_{l=1}^{N}\sum^N_{j=1, j \neq l}\cot(\pi(X_t^l-X_t^j))dX_t^l +\frac{1}{2N^2}\left(\sum_{l=1}^{N}\sum^N_{j=1, j \neq l}\frac{\pi^2}{\sin^2(\pi(X_t^l-X_t^j))}d_t<X_t^l>\right)\\
				&-\frac{1}{2N^2}\left(\sum_{l=1}^{N}\sum^N_{j=1, j \neq l}\frac{\pi^2}{\sin^2(\pi(X_t^l-X_t^j))}d_t< X_t^l, X_t^j> \right).
			\end{split}
		\end{equation*}
		Note that for the above three terms (denoted as $\mathbf{R}_\beta$,$\mathbf{C}^1_\gamma$ and $\mathbf{C}^2_\gamma$) , we have 
		\begin{equation*}
			\begin{split}
				\mathbf{R}_\beta&=d_t M_N(t)-\frac{\beta \pi}{2N^{3}}\sum_{l=1}^{N}\sum^N_{j,k=1,\, j,k\neq l}\cot(\pi(X_t^l-X_t^j))\cot(\pi(X_t^l-X_t^k))dt\\
				\mathbf{C}^1_\gamma&=\frac{1}{2N^2}\sum_{l=1}^{N}\sum^N_{j=1, j \neq l}\left(\frac{\pi^2}{\sin^2(\pi(X_t^l-X_t^j))}\sum_{k=-N}^{N}\frac{1}{|k|^{2\gamma}}|e_k(\frac{l}{N})|^2\right)dt\\
				\mathbf{C}^2_\gamma&=-\frac{1}{2N^2}\sum_{l=1}^{N}\sum^N_{j=1, j \neq l}\left(\frac{\pi^2}{\sin^2(\pi(X_t^l-X_t^j))}\sum_{k=-N}^{N}\frac{1}{|k|^{2\gamma}}e_k(\frac{l}{N})e_k(\frac{j}{N})\right)dt ,
			\end{split}
		\end{equation*}
		where $M_N(t)$ is a local martingale. It is known that, if $C=A+B$, then 
		$$\cot C=\frac{1-\tan A\tan B}{\tan A+\tan B}.$$
		Thus, we have 
		\begin{equation}\label{eqcot}
			-\cot A\cdot \cot B+ \cot A\cdot \cot C+ \cot B\cdot \cot C=-1.
		\end{equation}
		It follows that 
		$$\sum_{l=1}^{N}\sum^N_{j,k\neq l. j\neq k}\cot(\pi(X_t^l-X_t^j))\cot(\pi(X_t^l-X_t^k))=-\frac{N(N-1)(N-2)}{6}> -\frac{1}{6} N^3.$$
	Since
		\begin{equation*}
			\begin{split}
				&\sum_{l=1}^{N}\sum^N_{j,k=1,\, j,k\neq l }\cot(\pi(X_t^l-X_t^j))\cot(\pi(X_t^l-X_t^k))\\
				&=\sum_{l=1}^{N}\sum^N_{j \neq l }\cot^2(\pi(X_t^l-X_t^j))+\sum_{l=1}^{N}\sum^N_{j,k\neq l. j\neq k}\cot(\pi(X_t^l-X_t^j))\cot(\pi(X_t^l-X_t^k)),
			\end{split}
		\end{equation*}
		
we obtain
		\begin{equation}\label{est5}
			\begin{split}
				&d_t F(X_t^1,...X_t^N)\\
				&=\frac{\pi}{2N^2}\left(-\frac{\beta}{N}\sum_{1\leq l< j\leq N} \cot^2(\pi(X_t^l-X_t^j))+\sum_{1\leq l< j\leq N}\sum_{k=-N}^{N}\frac{\pi}{|k|^{2\gamma}}\left|\frac{e_k(\frac{l}{N})-e_k(\frac{j}{N})}{\sin(\pi(X_t^l-X_t^j))}\right|^2\right) dt\\
				&+d_t M_N(t)+C_\beta dt,
			\end{split}
		\end{equation}
	where $C_\beta <\frac{\beta \pi}{12}$. 
	
	\quad Next, we want to prove that $F(X_{t\wedge \tau_{R}}^1,...X_{t\wedge \tau_{R}}^N)-C t\wedge \tau_{R}$ is a super-martingale for some constant $C$. We divide the sum $\displaystyle \sum_{1\leq l< j\leq N}\sum_{k=-N}^{N}\frac{1}{|k|^{2\gamma}}\left|\frac{e_k(\frac{l}{N})-e_k(\frac{j}{N})}{\sin(\pi(X_t^l-X_t^j))}\right|^2$ into three parts:
		
		\begin{equation*}
			\begin{split}
				(A)_\gamma&=\left(\sum_{M=1}^{M_1-1}+\sum_{M=N+1-M_1}^{N}\right)\sum_{i=1}^{N}\sum_{k=-N}^{N}\frac{1}{|k|^{2\gamma}}\left|\frac{e_k(\frac{i}{N})-e_k(\frac{i+M}{N})}{\sin(\pi(X_t^i-X_t^{i+M}))}\right|^2\\
				(B)_\gamma&=\left(\sum_{M=M_1}^{M_2-1}+\sum_{M=N+1-M_2}^{N-M_1}\right)\sum_{i=1}^{N}\sum_{k=-N}^{N}\frac{1}{|k|^{2\gamma}}\left|\frac{e_k(\frac{i}{N})-e_k(\frac{i+M}{N})}{\sin(\pi(X_t^i-X_t^{i+M}))}\right|^2\\
			(C)_\gamma&=\sum_{M=M_2}^{N-M_2}\sum_{i=1}^{N-M}\sum_{k=-N}^{N}\frac{1}{|k|^{2\gamma}}\left|\frac{e_k(\frac{i}{N})-e_k(\frac{i+M}{N})}{\sin(\pi(X_t^i-X_t^{i+M}))}\right|^2.
			\end{split}
		\end{equation*}
		In the sequel, we denote $X_t^{k+N}=X_t^k$ for $k\geq 0$. 
		
		\textbf{Estimate $(A)_\gamma$:}	Note that 
		\begin{equation}\label{est1}
			|e_k(\frac{i}{N})-e_k(\frac{i+M}{N})|\leq |k|\frac{2\sqrt{2}\pi M}{ N} .
		\end{equation}
		It follows
		\begin{equation*}
			\begin{split}
				(A)_\gamma\leq& \left(\sum_{M=1}^{M_1-1}+\sum_{M=N+1-M_1}^{N}\right) \sum_{i=1}^{N-M}\sum_{k=-N,k\neq 0}^{N}\frac{1}{|k|^{2\gamma-2}}\frac{8\pi^2 M^2}{ N^2} \cdot \frac{1}{\sin^2(\pi(X_t^i-X_t^{i+M}))} \\
				=& \frac{2K_1^{N,\gamma}}{N^2}\sum_{M=1}^{M_1-1}M^2 (a_M +a_{N-M}) \\
				<&\frac{2K_1^{N,\gamma}M_1^2}{N^{2}}\sum_{M=1}^{M_1-1}(b_M+b_{N-M}+N).
			\end{split}
		\end{equation*}

 We can choose $M_1$ such that 

\begin{equation*}\label{assump1}
	M_1^2< N^{1-\eps}.
\end{equation*}  
Then we obtain
		\begin{equation}\label{esta1}
			(A)_\gamma\leq \frac{C'}{N^{1+\epsilon}}\left(\sum_{M=1}^{M_1-1}+\sum_{M=N+1-M_1}^{N}\right)b_M<\frac{1}{N^{\frac{\eps}{2}}}Q_N+C'N^{\frac{1}{2}},
		\end{equation}
		where we denote 
		\begin{equation*}
			\begin{split}
				Q_N&=\frac{1}{N}\sum_{1\leq l< j\leq N} \frac{1}{|\tan(\pi(X_t^i-X_t^{i+m}))|^2}.
			\end{split}
		\end{equation*}

		\textbf	{Estimate $(C)_\gamma$: } 	When $M$ is large,  \eqref{est1} is not enough to estimate $(B)_\gamma$ and $(C)_\gamma$. Since 
		$$|e_k(x)-e_k(y)|^2+|e_{-k}(x)-e_{-k}(y)|^2=8\sin^2(k\pi(x-y))<8,$$
		we obtain that
		\begin{equation*}
			\begin{split}
				(C)_\gamma&<\sum_{M=M_2}^{N-M_2}\sum_{i=1}^{N-M}\frac{8K^{N,\gamma}_2}{|\sin(\pi(X_t^i-X_t^{i+M}))|^2}.
 			\end{split}
		\end{equation*}
\eqref{compa} yields
	\begin{equation*}
		\begin{split}
			(C)_\gamma&<8K^{N,\gamma}_2\sum_{M=M_2}^{N-M_2}(b_M+N-M)\\
			&=8K^{N,\gamma}_2\sum_{M=M_2}^{N-M_2}b_M+8K_2^{N,\gamma}\sum_{M=M_2}^{N-M_2}(N-M). 
		\end{split}
	\end{equation*}
		Taking $M_2=\alpha N$ for $\alpha\in (0,1)$, we obtain by using \eqref{est2} and \eqref{compa} that
		\begin{equation}\label{estc1}
			\begin{split}
				(C)_\gamma&\leq 8K_2^{N,\gamma} \sum_{M=\alpha N}^{\frac{1}{2}N}\left( \Big(\frac{1}{M^2}+\frac{1}{(N-M)^2}\Big) (b_1+b_{N-1}) +\frac{N}{M^2} \right)+8K_2^{N,\gamma}N^2\\
				&<8K_2^\gamma Q_N  N\sum_{M=\alpha N}^{(1-\alpha)N} \frac{1}{M^2}+8K_2^{\gamma}N^2+\frac{8K_2^\gamma}{\alpha}\\
				&<\frac{8K_2^\gamma }{\alpha} Q_N-8K_2^\gamma (\sum_{M=(1-\alpha) N}^{N}\frac{1}{M^2})NQ_N+9K_2^\gamma N^2.
			\end{split}
		\end{equation}
		
		\textbf{Estimate $(B)_\gamma$:} Based on \eqref{est1} and \eqref{compa}, we have 
		\begin{equation*}
			\begin{split}
				(B)_\gamma&\leq \left(\sum_{M=M_1}^{\alpha N-1}+\sum_{M=N+1-\alpha N}^{N-M_1}\right) \sum_{i=1}^{N-M}\frac{2M^2}{N^2\sin^2(\pi(X_t^i-X_t^{i+M}))}\sum_{k=1}^{N}\frac{4\pi^2}{|k|^{2\gamma-2}}\\
				&\leq 2 K^{N,\gamma}_1 \sum_{M=M_1}^{\alpha N-1}\frac{M^2}{N^2}(b_M+b_{N-M}+N) .
			\end{split}
		\end{equation*}
		It then follows from \eqref{est2} that
		\begin{equation*}
			\begin{split}
				(B)_\gamma&< \frac{2 K^{\gamma}_1}{N^2}\sum_{M=M_1}^{\alpha N}M^2\left( \Big(\frac{1}{M^2}+\frac{1}{(N-M)^2}\Big) (b_1+b_{N-1}) +\frac{N}{M^2} \right)+\frac{2 K^{\gamma}_1}{N}\sum_{\alpha N}^{N}M^2\\
				&<\frac{2 K^{\gamma}_1 (b_1+b_{N-1})}{N^2} \sum_{M=M_1}^{\alpha N} \Big(1+\frac{M^2}{(N-M)^2}\Big) +K^{\gamma}_1N^2 +2K_1^\gamma \alpha.
			\end{split}
		\end{equation*}
		and hence
		\begin{equation}\label{estb1}
			(B)_\gamma<2 K_1^\gamma \alpha Q_N+2 K_1^\gamma \alpha^2(\sum_{M=(1-\alpha) N}^{N}\frac{1}{M^2})NQ_N+2 K^{\gamma}_1N^2.
		\end{equation}
	 Combined with \eqref{esta1}, \eqref{estc1} and \eqref{estb1}, we conclude that 
		\begin{equation*}
			(A)_\gamma+(B)_\gamma+(C)_\gamma<(\frac{1}{N^{\frac{\eps}{2}}}+2 K_1^\gamma \alpha +\frac{8K_2^\gamma }{\alpha})Q_N+C_\gamma N^2+ (-8K_2^\gamma +2 K_1^\gamma \alpha^2)(\sum_{M=(1-\alpha) N}^{N}\frac{1}{M^2})NQ_N,
		\end{equation*}
	for some constant $C_\gamma$ depending on $\gamma$. Taking $\alpha=2\sqrt{\frac{K_2^\gamma}{K_1^\gamma}}$, we obtain
		\begin{equation}\label{ABC}
		(A)_\gamma+(B)_\gamma+(C)_\gamma<8\sqrt{K_1^\gamma K_2^\gamma} Q_N +C_{\gamma}N^2,
	\end{equation}
	 for $N$ large enough. Therefore, when $\beta>8\pi \sqrt{K_1^\gamma K_2^\gamma}$,  $F(X_{t\land\tau_{R}}^1,...X_{t\land\tau_{R}}^N)-\frac{\pi}{2}(C_\gamma+ \frac{\beta }{6})t\land\tau_{R}$ is a super-martingale. The diffusion process $\{e^{2\pi iX_N^j(t\land \tau_{R})}\}_{j=1,...,N}$ on the torus is well defined almost surely, and for $T>0$, we let
		\begin{equation*}
			S=\{\tau_{R}\leq T\}.
		\end{equation*}
		Then
		\begin{equation*}
			\begin{split}
				F(X_0)+\frac{\pi}{2}(C_\gamma+ \frac{\beta }{6}) \tau_{R}\land T
				&\geq \E[F(X_{\tau_{R}\land T})]\\
				&=\E[F(X_{\tau_{R}})\mathbbm{1}_{S}]+\E[F(X_{T})\mathbbm{1}_{S^{c}}]\\
				&\geq -\frac{1}{N^2}\ln(\frac{1}{R})\P(S)-\frac{1}{2N^2}(N^2-N-2)\log2 \cdot \P(S)\\
				&-\frac{1}{2N^2}(N^2-N)\log2 \cdot \P(S^{c})\\
				&=\frac{1}{N^2}(\ln R+\log2)\P(S)-\frac{N-1}{2N}\log2 .
			\end{split}
		\end{equation*}
		Therefore,
		\begin{equation}\label{pest}
			\P(\tau_R\leq  T)\leq \frac{N^2(F(X_N(0))+\frac{\pi}{2}(C_\gamma+ \frac{\beta }{6})T+\log2)}{\ln R+\log2}.
		\end{equation}

		Next, following the same argument above, we show that when $\beta>\frac{1}{2}\pi K_1^\gamma$, $S^N(L_N(t\land \tau_R))-Ct\land \tau_R$ is a supermartingale for some positive number $C$. In fact, by It\^o's formula,
		\begin{equation*}
			\begin{split}
				&d_t S^N(L_N(t))\\&=-\sum_{i=1}^{N}\frac{\pi}{N\tan(\pi(X^i_N(t)-X_N^{i-1}(t)))}\circ(dX^i_N(t)-dX_N^{i-1}(t))\\
				&=-\sum_{i=1}^{N}\frac{\pi}{N\tan(\pi(X^i_N(t)-X_N^{i-1}(t)))}\frac{\beta}{N}\\
				&\cdot \left(\sum_{j=1,j\neq i}^{N}\frac{1}{\tan(\pi(X^i_N(t)-X_N^{j}(t)))}-\sum_{j=1,j\neq i-1}^{N}\frac{1}{\tan(\pi(X^{i-1}_N(t)-X_N^{j}(t)))}\right)dt\\
				&+\frac{1}{2N}\sum_{i=1}^{N}\frac{\pi^2}{\sin^2(\pi(X^i_N(t)-X_N^{i-1}(t)))}\sum_{k=-N}^{N} \frac{1}{k^{2\gamma}}\left|e_k(\frac{i+1}{N})-e_k(\frac{i}{N})\right|^2 dt+dM_N(t)\\
				&=\beta a(t)dt+b(t)dt+dM_N(t).
			\end{split}
		\end{equation*}
		Here, we just denote $X_N^0$ as $X_N^N-1$ for simplicity of notations. Since \eqref{est1}, we have
		$$b(t)\leq \frac{1}{N}\sum_{i=1}^{N}\frac{\pi^2}{\sin^2(\pi(X^i_N(t)-X_N^{i-1}(t)))}\sum_{k=1}^{N} \frac{4\pi^2}{N^2}\frac{1}{k^{2\gamma-2}}.$$
		Following the notations in the introduction, $a(t)$ formally equals $\displaystyle -\int_\T  \partial_x(\H L_N(t))\nu_t^N(x) dx$. Note that 
		\begin{equation*}
			\begin{split}
				a(t)&=-\frac{2\pi }{N^2}\sum_{i=1}^{N}\frac{1}{\tan^2\big(\pi(X^i_N(t)-X_N^{i-1}(t))\big)}\\
				&-\frac{\pi }{N^2}\sum_{i=1}^{N}\frac{1}{\tan\big(\pi(X^i_N(t)-X_N^{i-1}(t))\big)}\\
				&\cdot \left(\sum_{j=1,j\neq i,i-1}^{N}\frac{1}{\tan(\pi(X^i_N(t)-X_N^{j}(t)))}-\sum_{j=1,j\neq i-1,i}^{N}\frac{1}{\tan(\pi(X^{i-1}_N(t)-X_N^{j}(t)))}\right)\\
				&=a_1(t)+a_2(t).
			\end{split}
		\end{equation*}
		
	In the sequel, we denote 
		$$Z_N^i(t):=\cot(\pi(X^{i}_N(t)-X^{i-1}_N(t))).$$ 
		By \eqref{eqcot}, we have
		\begin{equation*}
			\begin{split}
				a_2(t)&=\frac{\pi }{N^2}\sum_{i=1}^{N}\sum_{j=1,j\neq i,i-1}^{N}\left(1+\cot(\pi(X^i_N(t)-X_N^{j}(t)))\cot(\pi(X^{i-1}_N(t)-X_N^{j}(t))\right)\\
				&\leq \pi +  \frac{\pi }{N^2}\sum_{i=1}^{N}\sum_{j=1,j\neq i,i-1}^{N}\frac{1}{\tan(\pi(X^i_N(t)-X_N^{j}(t))\tan(\pi(X^{i-1}_N(t)-X_N^{j}(t))}\\
				&=\pi +a_3(t).
			\end{split}
		\end{equation*}
		It follows from Lemma \ref{cotineq} for $p=0$  in Appendix, that
		$$a_1(t)+a_3(t)<-\frac{1}{N}\frac{2\pi }{N^2}\sum_{i=1}^{N}\big(Z_N^i(t)\big)^2.$$
		Then it yields
		\begin{equation}\label{est6}
			\begin{split}
				\beta a(t)+b(t)&\leq K_1^\gamma  \frac{1}{N}\sum_{i=1}^{N}\frac{\pi^2}{N^2\sin^2(\pi(X^i_N(t)-X_N^{i-1}(t)))}+\beta \pi- \frac{2\beta}{N\pi}\frac{\pi^2}{N^2}\sum_{i=1}^{N}\big(Z_N^i(t)\big)^2\\
				&\leq  -2  \frac{\beta-\frac{1}{2}\pi K_1^\gamma}{\pi}\sum_{i=1}^N\frac{1}{N}\frac{\pi^2}{N^2\sin^2(\pi(X^i_N(t)-X_N^{i-1}(t)))}+\beta\pi(1+\frac{2}{N}).
			\end{split}
		\end{equation}
		So, when $\beta>\frac{1}{2}\pi K_1^\gamma$, $S^N(L_N(t))-2\beta \pi t\wedge \tau_R$ is a supermartingale. Then by the same argument, we have
		\begin{equation*}
			\begin{split}
				S^N(L_N(0))+2\beta \pi T\wedge \tau_R 
				&\geq \E[S^N(L_N(\tau_{R}\land T))]\\
				&=\E[S^N(L_N(\tau_{R}))\mathbbm{1}_{S}]+\E[S^N(L_N(T))\mathbbm{1}_{S^{c}}]\\
				&\geq (\frac{1}{N}\ln\frac{2\pi R}{N}+\frac{N-1}{N}\ln \frac{2\pi}{N})\P(S)+\ln \frac{2\pi}{N} \P(S^c)\\
				&=\P(S)\frac{1}{N}\Big(\ln\frac{2\pi R}{N}-\ln\frac{2\pi}{N}\Big)+(2-\frac{1}{N})\ln \frac{2\pi}{N}
			\end{split}
		\end{equation*}
	Therefore, 
	\begin{equation}\label{pest'}
		\P(\tau_R\leq T)\leq \frac{N(S^N(L_N(0))+2\beta \pi T-\ln \frac{2\pi}{N})}{\ln\frac{2\pi R}{N}-\ln\frac{2\pi}{N}}.
	\end{equation}
For fixed $ T$, letting $ R\to \infty$, both \eqref{pest} and \eqref{pest'} yield $\P(\tau_{\infty}\leq  T)=0$. It follows that when $\beta>\lambda_0^\gamma$, $\{(e^{2\pi iX^j_{t\land T}})\}_{j=1,\ldots, N}$ never collide. Then letting $T\to \infty$ , there is no collision among the particles $\{e^{2\pi iX_N^j(t)}\}_{j=1,...,N}$ for all $t\in[0,+\infty)$. Coming back to the original process , this means $\{X^j_N(t)\}_{j=1,\ldots, N}$ never collide and $|X_N^1(t)-X^N_N(t)|<1$. Also, by continuity of  the trajectories of  $X_N(t)$ , we have $X_N(t)\in \Delta_N$ for all $t\geq 0$. 

\quad Note that $F(X_N(t))=\textbf{W}(L_N(t))$. Now, based on \eqref{est5}, \eqref{ABC} and the well-posedness of $X_N(t)$, $\textbf{W}(L_N(t))-\frac{\pi}{2}(C_\gamma +\frac{\pi}{6}) t$ is a supermartingale. Also, by \eqref{est6} and $\frac{1}{x^2}>-\ln x$, we have
\begin{equation*}\label{sn}
	\frac{d}{dt}\E[S^N(L_N(t))]<-2 \frac{\beta-\frac{1}{2}\pi K_1^\gamma}{\pi} \E[S^N(L_N(t))]+2\beta \pi.
\end{equation*}
By Gronwall inequality, 
\begin{equation*}
	\E[S^N(L_N(t))]\leq e^{-2 C_1^\beta t}S^N(L_N(0)) +(1-e^{-2C_1^\beta t})\frac{\beta \pi}{C_1^\beta},
\end{equation*}
where $C_1^\beta=\frac{\beta-\frac{1}{2}\pi K_1^\gamma}{\pi}$. \eqref{ugr'} is proved. We finish the proof .

	\end{proof}
	\vskip 3mm
	\begin{remark}
		For $\gamma<1$, one can prove that the particle model \eqref{model1} has collisions in finite time almost surely. Indeed, let $Y^i_t\in [0,1]$ be the solution to
		
$$dY^i_t=\frac{\beta}{N}\cot\big(\pi Y_t^i(t)\big) d t +\sum_{k=-N}^{N}\frac{1}{|k|^{\gamma}}\Big(e_k(\frac{i+1}{N})-e_k(\frac{i}{N})\Big) dW_t^k, $$
with $Y_0^i=X_N^{i+1}(0)-X_N^{i}(0)$. By the argument in \cite{cepa2001brownian} Theorem 3.1 , we can deduce that, when $\gamma<1$, $Y_t^i$ touches $0$ and $1$ a.s. and 
$$0\leq X_N^{i+1}-X_N^i\leq Y_t^i;  \  \ Y_t^N\leq X_N^{N}-X_N^1\leq 2\pi.$$
However, we do not know whether there are collisions in finite time for the case $\gamma\in[1,\frac{3}{2}]$.
	\end{remark}

	\subsection{Rate of convergence}\label{subsec4}

	\quad In this subsection, we will give the convergence rate of the particle approximation, and hence prove that there is a unique fluctuating hydrodynamic limit of the particle model \eqref{model1}. Firstly, we give an lemma for preparation. By \cite{V03} Proposition 2.24, the optimal transport plan from $\bar{L}_M(t)$  to $\bar{L}_N(t)$, as two discrete measures in $\P_2(\R)$ for $N<M$, should have a support which is cyclically monotone; and by Rockafellar's theorem on cyclical monotonicity, there exist two nondecreasing surjective maps:
	$$\tau^k:\, \{1,2,\ldots,M\}\rightarrow \{1,2,\ldots, N\} , \\\ \ \  k=1,2,$$ 
	and "mass allocation" functions $a^k:\, \{1,2,\ldots,M\}\rightarrow [0,1]$  for $k=1,2$, such that the amount $a^k(m)$ of the weight of the $m-$th atom at $X_M^m(t)$ is transported to $X_N^{\tau^k(m)}(t)$ along a straight line. It should be satisfied that
	$$\tau^2-\tau^1=1, \ \ \ \ \ \  a^1+a^2=1,$$
	and 
	\begin{equation}\label{cod5}
		\frac{\tau^1(m)-1}{N}\leq \frac{m}{M}\leq \frac{\tau^2(m)}{N}.
	\end{equation}
	The optimal transport plan $\gamma_t^{1,2}$ from $\bar{L}_M(t)$  to $\bar{L}_N(t)$ is given by
	\begin{equation}\label{optmap}
		\sum_{i=1}^{M}\sum_{k=1}^{2}\frac{a^k(i)}{M}\delta_{(X_M^i(t),X_N^{\tau^k(i)}(t))};
	\end{equation}
and the optimal transport plan $\gamma_t^{2,1}$ from $\bar{L}_N(t)$  to $\bar{L}_M(t)$ is given by
\begin{equation}\label{optmap2}
	\sum_{i=1}^{M}\sum_{k=1}^{2}\frac{a^k(i)}{M}\delta_{(X_N^{\tau^k(i)}(t),X_M^i(t))}.
\end{equation}
	Notice that $(X_M)$ and $(X_N)$ both keep orders, thus $\tau^k$ and $a^k$ actually do not change over time and \eqref{optmap}, \eqref{optmap2} are well defined. Let the interaction energy on $\P(\R)$ be 
	$$\bar{\textbf{W}}(\mu)=-\int_\R\int_\R \ln\big|e^{2\pi ix}-e^{2\pi iy}\big|d\mu(x)d\mu(y).$$
Obviously, $\bar{\textbf{W}}$ is not geodesically convex on $\P_2(\R)$.

	
	\begin{lemma}\label{lem3} 
		For $\forall \, M>N>1$, we have almost surely that  
		\begin{equation*}\label{S1}
			\int_{[0,1]^2}\<x_1-x_2, \H\bar{L}_N(t)(x_1)\>d\gamma_t^{2,1}(x_1,x_2)+  \int_{[0,1]^2}\<x_2-x_1, \H\bar{L}_M(t)(x_2)\>d\gamma_t^{1,2}(x_2,x_1)\leq \frac{2}{N}+\frac{2}{M}.
		\end{equation*}
	\end{lemma}
\begin{proof}
	By the representation \eqref{optmap} and \eqref{optmap2}, we have
	\begin{equation*}
		\begin{split}
			&	\int_{[0,1]^2}\<x_2-x_1, \H\bar{L}_N(t)(x_1)\>d\gamma_t^{2,1}(x_1,x_2)+  \int_{[0,1]^2}\<x_1-x_2, \H\bar{L}_M(t)(x_2)\>d\gamma_t^{1,2}(x_2,x_1)\\
			&=\sum_{i=1}^{M}\frac{a^1(i)}{M}(X_M^i(t)-X_N^{\tau^1(i)}(t))\Big(\frac{1}{M}\sum_{j=1,j\neq i}^{M}\cot\big(\pi(X_M^i(t)-X_M^j(t))\big)\Big) \\
			&-\sum_{i=1}^{M}\frac{a^1(i)}{M}(X_M^i(t)-X_N^{\tau^1(i)}(t))\Big(\frac{1}{N}\sum_{j=1,j\neq \tau^1(i)}^{N}\cot\big(\pi(X_N^{\tau^1(i)}(t)-X_N^j(t))\big)\Big)\\
			&+\sum_{i=1}^{M}\frac{a^2(i)}{M}(X_M^i(t)-X_N^{\tau^2(i)}(t))\Big(\frac{1}{M}\sum_{j=1,j\neq i}^{M}\cot\big(\pi(X_M^i(t)-X_M^j(t))\big)\Big)\\
			&-\sum_{i=1}^{M}\frac{a^2(i)}{M}(X_M^i(t)-X_N^{\tau^2(i)}(t))\Big(\frac{1}{N}\sum_{j=1,j\neq \tau^2(i)}^{N}\cot\big(\pi(X_N^{\tau^2(i)}(t)-X_N^j(t))\big)\Big).
		\end{split}
	\end{equation*}
Denote the right hand side of the equality as $S$. We only need to prove $S\leq \frac{2\beta}{N}+\frac{2\beta}{M}$. By the well-posedness of the particle model, we have for any $(t,\omega)\in (0,\infty)\times \Omega$ that
	$$\bar{\textbf{W}}(\bar{L}_M(t))=-\frac{1}{M^2}\sum_{i<j}^{M}2\ln\Big|\sin\big(\pi (X_M^i(t)- X_M^j(t))\big)\Big|<\infty;$$
	$$\bar{\textbf{W}}(\bar{L}_N(t))=-\frac{1}{N^2}\sum_{ i< j}^{N}2\ln\Big|\sin\big(\pi (X_N^i(t)- X_N^j(t))\big)\Big|<\infty.$$
	
	Let 
	$$Z_{t}^{i,j}(\lambda)=\lambda X_N^j(t)+(1-\lambda)X_M^i(t),  \ \ \  \ \lambda\in[0,1],$$
	then 
	$$\bar{L}_{Z_t}(\lambda)=\sum_{i=1}^{M}\frac{a^1(i)}{M}\delta_{Z_t^{i,\tau^1(i)}(\lambda)}+\sum_{i=1}^{M}\frac{a^2(i)}{M}\delta_{Z_t^{i,\tau^2(i)}(\lambda)}$$
	is the geodesic from $\bar{L}_M(t)$ to $\bar{L}_N(t)$ according to the optimal transport plan \eqref{optmap}.
	Note that 
	\begin{equation*}
		\begin{split}
			&\bar{\textbf{W}}(\bar{L}_{Z_t}(\lambda))\\
			&=-\sum_{ i\leq j}^M\sum_{k,l=1}^2\frac{2a^k(i)a^l(j)}{M^2}\ln\Big|\sin \pi\big( Z_t^{i,\tau^k(i)}(\lambda)- Z_t^{j, \tau^l(j)}(\lambda)\big)\Big| \qquad \qquad \qquad \qquad \qquad \qquad 
		\end{split}
	\end{equation*}
			\begin{equation*}
				\begin{split}
			&=-\sum_{(i,j,k,l)\in I_1}\frac{2a^k(i)a^l(j)}{M^2}\ln\Big|\sin \pi\Big((1-\lambda)(X_M^i(t)-X_M^j(t))\Big)\Big|\\
			&-\sum_{(i,j,k,l)\in I_2}\frac{2a^k(i)a^l(j)}{M^2}\ln\Big|\sin \pi\Big(\lambda(X_N^{\tau^k(i)}(t)-X_N^{\tau^l(j)}(t))\Big)\Big|\\
			&-\sum_{(i,j,k,l)\in I_c}\frac{2a^k(i)a^l(j)}{M^2}\ln\Big|\sin \pi\Big((1-\lambda)(X_M^i(t)-X_M^j(t))+\lambda(X_N^{\tau^k(i)}(t)-X_N^{\tau^l(j)}(t))\Big)\Big|.
		\end{split}
	\end{equation*}
	where the sum for each line is $D_1(1-\lambda)$, $D_2(\lambda)$, $D_c(\lambda)$, and 
	$$I_1=\{(i,j,k,l): \ i,j\in \{1,2,\ldots,M\},\, k,l\in \{1,2\},\, i<j, \, \tau^k(i)=\tau^l(j), \, a^k(i)>0, \, a^l(j)>0\},$$
	$$I_2=\{(i,j,k,l): \ i,j\in \{1,2,\ldots,M\},\, k,l\in \{1,2\},\, i=j, \, k<l, \, a^k(i)>0, \, a^l(i)>0\},$$
	$$I_c=\{(i,j,k,l): \ i,j\in \{1,2,\ldots,M\}, \, k,l\in \{1,2\}, \, i< j, \, a^k(i)>0, \, a^l(j)>0\}/I_1.$$
	It follows from Theorem \ref{introthm1} that $(X_M(t)) \in \Delta_M, \  (X_N(t))\in \Delta_N$, and hence for $(i,j,k,l)\in  I_c$, we have 
	$$-1<X_M^i(t)-X_M^j(t)<0 \, ; \ \ \  -1<X_N^{\tau^k(i)}(t)-X_N^{\tau^l(j)}(t)<0.$$
	Since $-\ln|\sin (\pi x)|$ is convex on $(-1,0)$, we have by Jensen inequality that 
	\begin{equation*}\label{conequ}
		\begin{split}
			&D_c(\lambda)\\
			&\leq -\sum_{(i,j,k,l)\in I_c}\frac{2a^k(i)a^l(j)}{M^2}\left[(1-\lambda)\ln\big|\sin\pi(X_M^i(t)-X_M^j(t))\big|+\lambda\ln\big|\sin\pi(X_N^{\tau^k(i)}(t)-X_N^{\tau^l(j)}(t))\big| \right]\\
			&=(1-\lambda)\bar{\textbf{W}}(\bar{L}_M(t))+\lambda\bar{\textbf{W}}(\bar{L}_N(t))+(1-\lambda)\sum_{(i,j,k,l)\in I_1}\frac{2a^k(i)a^l(j)}{M^2}\ln\big|\sin\pi(X_M^i(t)-X_M^j(t))\big|\\
			&+\lambda \sum_{(i,j,k,l)\in I_2}\frac{2a^k(i)a^l(j)}{M^2}\ln\Big|\sin \pi\Big(X_N^{\tau^k(i)}(t)-X_N^{\tau^l(j)}(t)\Big)\Big|\\
			&=(1-\lambda)\bar{\textbf{W}}(\bar{L}_M(t))+\lambda\bar{\textbf{W}}(\bar{L}_N(t)+\sum_{(i,j,k,l)\in I_1}\frac{2a^k(i)a^l(j)}{M^2}\ln\big|\sin\pi(X_M^i(t)-X_M^j(t))\big|+\lambda A_t,
		\end{split}
	\end{equation*}
	where we denote
	\begin{equation*}
		\begin{split}
			A_t&=-\sum_{(i,j,k,l)\in I_1}\frac{2a^k(i)a^l(j)}{M^2}\ln\big|\sin\pi(X_M^i(t)-X_M^j(t))\big|\\
			&+\sum_{(i,j,k,l)\in I_2}\frac{2a^k(i)a^l(j)}{M^2}\ln\Big|\sin \pi\Big(X_N^{\tau^k(i)}(t)-X_N^{\tau^l(j)}(t)\Big)\Big|,
		\end{split}
	\end{equation*}
	and the equality in the second line is due to 
	$$\bar{\textbf{W}}(\bar{L}_M(t))=-\sum_{(i,j,k,l)\in I_1\cup I_c}\frac{2a^k(i)a^l(j)}{M^2}\ln\big|\sin\pi(X_M^i(t)-X_M^j(t))\big|$$
	and 
	$$\bar{\textbf{W}}(\bar{L}_N(t))=-\sum_{(i,j,k,l)\in I_2\cup I_c}\frac{2a^k(i)a^l(j)}{M^2}\ln\big|\sin\pi(X_N^{\tau^k(i)}(t)-X_M^{\tau^l(j)}(t))\big|.$$
	Therefore,
	\begin{equation*}
		\begin{split}
			&\bar{\textbf{W}}(\bar{L}_{Z_t}(\lambda))-D_1(1-\lambda)-D_2(\lambda)-\sum_{(i,j,k,l)\in I_1}\frac{2a^k(i)a^l(j)}{M^2}\ln\big|\sin\pi(X_M^i(t)-X_M^j(t))\big|\\
			\leq & (1-\lambda)\bar{\textbf{W}}(\bar{L}_M(t))+\lambda\bar{\textbf{W}}(\bar{L}_N(t)) +\lambda A_t.
		\end{split}
	\end{equation*}
	
	This is equivalent to 
	\begin{equation}\label{equmn}
		\begin{split}
			&\bar{\textbf{W}}(\bar{L}_{Z_t}(\lambda))-D_2(\lambda)-\bar{\textbf{W}}(L_M(t))-D_1(1-\lambda)\\
			&-\sum_{(i,j,k,l)\in I_1}\frac{2a^k(i)a^l(j)}{M^2}\ln\big|\sin\pi(X_M^i(t)-X_M^j(t))\big|\\
			&\leq \lambda\big(\bar{\textbf{W}}(\bar{L}_N(t))-\bar{\textbf{W}}(\bar{L}_M(t))\big)+\lambda A_t.
		\end{split}
	\end{equation}
	
	Analogously, we can obtain by replacing $\lambda$ by $1-\lambda$ in \eqref{equmn} that 
	\begin{equation}\label{equnm}
		\begin{split}
			&\bar{\textbf{W}}(\bar{L}_{Z_t}(1-\lambda))-D_1(\lambda)-\bar{\textbf{W}}(\bar{L}_N(t))-D_2(1-\lambda)\\
			&-\sum_{(i,j,k,l)\in I_2}\frac{2a^k(i)a^l(j)}{M^2}\ln\big|\sin\pi(X_N^{\tau^k(i)}(t)-X_N^{\tau^l(j)}(t))\big|\\
			&\leq \lambda\big(\bar{\textbf{W}}(\bar{L}_M(t))-\bar{\textbf{W}}(\bar{L}_N(t))\big)-\lambda A_t.
		\end{split}
	\end{equation}

	Let $\lambda \to 0$, we have
	\begin{equation}\label{equsu1}
		\begin{split}
			& \lim_{\lambda\to 0^+}\frac{1}{\lambda}\left[\bar{\textbf{W}}(\bar{L}_{Z_t}(\lambda))-D_2(\lambda)-\bar{\textbf{W}}(\bar{L}_M(t))\right] \\
			&=\sum_{i=1}^{M}\frac{a^1(i)}{M}(X_M^i(t)-X_N^{\tau^1(i)}(t))\Big(\frac{\pi}{M}\sum_{j=1,j\neq i}^{M}\cot\big(\pi(X_M^i(t)-X_M^j(t))\big)\Big)\\
			&+\sum_{i=1}^{M}\frac{a^2(i)}{M}(X_M^i(t)-X_N^{\tau^2(i)}(t))\Big(\frac{\pi}{M}\sum_{j=1,j\neq i}^{M}\cot\big(\pi(X_M^i(t)-X_M^j(t))\big)\Big);
		\end{split}
	\end{equation}
	and 
	\begin{equation}\label{equsu2}
		\begin{split}
			& \lim_{\lambda\to 0^+}\frac{1}{\lambda}\left[\bar{\textbf{W}}(\bar{L}_{Z_t}(1-\lambda))-D_1(\lambda)-\bar{\textbf{W}}(\bar{L}_N(t))\right] \\
			&=\sum_{i=1}^{M}\frac{a^1(i)}{M}(X_M^i(t)-X_N^{\tau^1(i)}(t))\Big(\frac{\pi}{N}\sum_{j=1,j\neq \tau^1(i)}^{N}\cot\big(\pi(X_N^{\tau^1(i)}(t)-X_N^j(t))\big)\Big)\\
			&+\sum_{i=1}^{M}\frac{a^2(i)}{M}(X_M^i(t)-X_N^{\tau^2(i)}(t))\Big(\frac{\pi}{N}\sum_{j=1,j\neq \tau^2(i)}^{N}\cot\big(\pi(X_N^{\tau^2(i)}(t)-X_N^j(t))\big)\Big).
		\end{split}
	\end{equation}
	
	Also, a simple calculation gives
	\begin{equation}\label{equsu3}
		\begin{split}
			& \lim_{\lambda\to 0^+}\frac{1}{\lambda}\left[-D_1(1-\lambda)-\sum_{(i,j,k,l)\in I_1}\frac{2a^k(i)a^l(j)}{M^2}\ln\big|\sin\pi(X_M^i(t)-X_M^j(t))\big|\right] \\
			&=-\sum_{(i,j,k,l)\in I_1}\frac{2\pi a^k(i)a^l(j)}{M^2}\cot\pi(X_M^i(t)-X_M^j(t)) \cdot (X_M^i(t)-X_M^j(t)) ;
		\end{split}
	\end{equation}
	and
	\begin{equation}\label{equsu4}
		\begin{split}
			& \lim_{\lambda\to 0^+}\frac{1}{\lambda}\left[-D_2(1-\lambda)-\sum_{(i,j,k,l)\in I_2}\frac{2a^k(i)a^l(j)}{M^2}\ln\big|\sin\pi(X_N^{\tau^k(i)}(t)-X_N^{\tau^l(j)}(t))\big|\right] \\
			&=-\sum_{(i,j,k,l)\in I_2}\frac{2\pi a^k(i)a^l(j)}{M^2}\cot\pi(X_N^{\tau^k(i)}(t)-X_N^{\tau^l(j)}(t))\cdot (X_N^{\tau^k(i)}(t)-X_N^{\tau^l(j)}(t)).
		\end{split}
	\end{equation}
	
	Now, it follows from \eqref{equmn}, \eqref{equnm}, \eqref{equsu1}, \eqref{equsu2}, \eqref{equsu3} and \eqref{equsu4} that 
	\begin{equation}\label{sumequ1}
		\begin{split}
			&\sum_{i=1}^{M}\frac{a^1(i)}{M}(X_M^i(t)-X_N^{\tau^1(i)}(t))\Big(\frac{\pi}{M}\sum_{j=1,j\neq i}^{M}\cot\big(\pi(X_M^i(t)-X_M^j(t))\big)\Big)\\
			&+\sum_{i=1}^{M}\frac{a^2(i)}{M}(X_M^i(t)-X_N^{\tau^2(i)}(t))\Big(\frac{\pi}{M}\sum_{j=1,j\neq i}^{M}\cot\big(\pi(X_M^i(t)-X_M^j(t))\big)\Big)\\
			&\leq \bar{\textbf{W}}(\bar{L}_N(t))-\bar{\textbf{W}}(\bar{L}_M(t))+ A_t\\
			&+\sum_{(i,j,k,l)\in I_1}\frac{2\pi a^k(i)a^l(j)}{M^2}\cot\pi(X_M^i(t)-X_M^j(t)) \cdot (X_M^i(t)-X_M^j(t)), 
		\end{split}
	\end{equation}
	and
	\begin{equation}\label{sumequ2}
		\begin{split}
			&\sum_{i=1}^{M}\frac{a^1(i)}{M}(X_M^i(t)-X_N^{\tau^1(i)}(t))\Big(\frac{\pi}{N}\sum_{j=1,j\neq \tau^1(i)}^{N}\cot\big(\pi(X_N^{\tau^1(i)}(t)-X_N^j(t))\big)\Big)\\
			&+\sum_{i=1}^{M}\frac{a^2(i)}{M}(X_M^i(t)-X_N^{\tau^2(i)}(t))\Big(\frac{\pi}{N}\sum_{j=1,j\neq \tau^2(i)}^{N}\cot\big(\pi(X_N^{\tau^2(i)}(t)-X_N^j(t))\big)\Big)\\
			&\leq \bar{\textbf{W}}(\bar{L}_M(t))-\bar{\textbf{W}}(\bar{L}_N(t))- A_t\\
			&+\sum_{(i,j,k,l)\in I_2}\frac{2\pi a^k(i)a^l(j)}{M^2}\cot\pi(X_N^{\tau^k(i)}(t)-X_N^{\tau^l(j)}(t))\cdot (X_N^{\tau^k(i)}(t)-X_N^{\tau^l(j)}(t)).
		\end{split}
	\end{equation}
	We obtain by summing up \eqref{sumequ1} and \eqref{sumequ2} that 
	\begin{equation*}
		\begin{split}
	S&\leq \sum_{(i,j,k,l)\in I_1}\frac{2\pi a^k(i)a^l(j)}{M^2}\cot\pi(X_M^i(t)-X_M^j(t)) \cdot (X_M^i(t)-X_M^j(t))\\
			&+\sum_{(i,j,k,l)\in I_2}\frac{2\pi a^k(i)a^l(j)}{M^2}\cot\pi(X_N^{\tau^k(i)}(t)-X_N^{\tau^l(j)}(t))\cdot (X_N^{\tau^k(i)}(t)-X_N^{\tau^l(j)}(t)).
		\end{split}
	\end{equation*}
	Note that for $x\in (-1,0)$, it holds
	$$\cot\pi x \cdot x<\frac{1}{\pi}.$$
	Thus, 
	$$S\leq \sum_{(i,j,k,l)\in I_1}\frac{2 a^k(i)a^l(j)}{M^2}+\sum_{(i,j,k,l)\in I_2}\frac{2 a^k(i)a^l(j)}{M^2}.$$
	Next, we estimate the amounts of $I_1$ and $I_2$. In fact, there are at most $\frac{M}{N}+1$ particles in $\{X^i_M(t)\}_{i=1}^M$ transported to one particle in $\{X^i_N(t)\}_{i=1}^N$; one particle in $\{X^i_M(t)\}_{i=1}^M$ could be transported to at most two particles in $\{X^i_N(t)\}_{i=1}^N$. Hence, we have
	$$|I_1|\leq \frac{N}{2}(\frac{M}{N}+1)\frac{M}{N},$$
	and 
	$$|I_2|\leq M.$$
	We then obtain
	\begin{equation*}
		S(t,\omega)< \frac{2}{N}+\frac{2}{M},  \ \ \  \forall \, (t,\omega)\in (0,\infty)\times \Omega.
	\end{equation*}
\end{proof}

\vskip 3mm

Now, we prove Theorem \ref{ratecon}.

	\begin{proof}[\textbf{Proof of Theorem \ref{ratecon}}]
	By the representation \eqref{optmap}, we have,
		$$\bar{W}_2^2(\bar{L}_N(t), \bar{L}_M(t))=\frac{1}{M}\sum_{i=1}^{M}a^1(i)|X_M^i(t)-X_N^{\tau^1(i)}(t)|^2+\frac{1}{M}\sum_{i=1}^{M}a^2(i)|X_M^i(t)-X_N^{\tau^2(i)}(t)|^2.$$
	Using the formula above, we now compute the evolution of $\E[\bar{W}_2(\bar{L}_N(t), \bar{L}_M(t))]$. By It\^o's formula, we have
		\begin{equation*}
			\begin{split}
				&d \bar{W}_2^2(\bar{L}_N(t),\bar{L}_M(t))\\
				&=\sum_{i=1}^{M}\frac{a^1(i)}{M}(X_M^i(t)-X_N^{\tau^1(i)}(t))\circ d_t(X_M^i(t)-X_N^{\tau^1(i)}(t)) \\
				&+\frac{a^2(i)}{M}\sum_{i=1}^{M}(X_M^i(t)-X_N^{\tau^2(i)}(t))\circ d_t(X_M^i(t)-X_N^{\tau^2(i)}(t)) \\
				&=\sum_{i=1}^{M}\frac{a^1(i)}{M}(X_M^i(t)-X_N^{\tau^1(i)}(t))\Big(\frac{\beta}{M}\sum_{j=1,j\neq i}^{M}\cot\big(\pi(X_M^i(t)-X_M^j(t))\big)\Big) dt\\
				&-\sum_{i=1}^{M}\frac{a^1(i)}{M}(X_M^i(t)-X_N^{\tau^1(i)}(t))\Big(\frac{\beta}{N}\sum_{j=1,j\neq \tau^1(i)}^{N}\cot\big(\pi(X_N^{\tau^1(i)}(t)-X_N^j(t))\big)\Big)dt \qquad  \\
				&+\sum_{i=1}^{M}\frac{a^2(i)}{M}(X_M^i(t)-X_N^{\tau^2(i)}(t))\Big(\frac{\beta}{M}\sum_{j=1,j\neq i}^{M}\cot\big(\pi(X_M^i(t)-X_M^j(t))\big)\Big)dt\\
				&-\sum_{i=1}^{M}\frac{a^2(i)}{M}(X_M^i(t)-X_N^{\tau^2(i)}(t))\Big(\frac{\beta}{N}\sum_{j=1,j\neq \tau^2(i)}^{N}\cot\big(\pi(X_N^{\tau^2(i)}(t)-X_N^j(t))\big)\Big)dt\\
				&+\sum_{i=1}^{M}\sum_{k=-N}^{N}\frac{a^1(i)}{2M|k|^{2\gamma}}\Big(e_k(\frac{i}{M})-e_k(\frac{\tau^1(i)}{N})\Big)^2+\frac{a^2(i)}{2M|k|^{2\gamma}}\Big(e_k(\frac{i}{M})-e_k(\frac{\tau^2(i)}{N})\Big)^2 dt \ \  \\
				&+\sum_{i=1}^{M}\sum_{|k|=N+1}^{M}\frac{a^1(i)+a^2(i)}{2M|k|^{2\gamma}}e_k^2(\frac{i}{M}) dt \\
				&+\sum_{i=1}^{M}\frac{a^1(i)}{M}(X_M^i(t)-X_N^{\tau^1(i)}(t))\Big(\sum_{k=-M}^{M}\frac{1}{|k|^{\gamma}}e_k(\frac{i}{M})dW_t^k-\sum_{k=-N}^{N}\frac{1}{|k|^\gamma}e_k(\frac{\tau^1(i)}{N})dW_t^k\Big)\\
				&+\sum_{i=1}^{M}\frac{a^2(i)}{M}(X_M^i(t)-X_N^{\tau^2(i)}(t))\Big(\sum_{k=-M}^{M}\frac{1}{|k|^{\gamma}}e_k(\frac{i}{M})dW_t^k-\sum_{k=-N}^{N}\frac{1}{|k|^\gamma}e_k(\frac{\tau^1(i)}{N})dW_t^k\Big)\\
				&=S_1 dt+S_2 dt+dM_t,
			\end{split}
		\end{equation*}
		where we denote the last two line as a differential of a local martingale $M_t$ and the sum of the two lines above the martingale terms as $S_2$ and the sum of the other terms as $S_1$. Also denote the quadratic variation process as $<M>_t$. Next, we will estimate these three terms.
		
		\textbf{Estimate on $S_2$:} Note that for $l=1,2$, we have
		\begin{equation}\label{iden3}
			\Big(e_k(\frac{i}{M})-e_k(\frac{\tau^l(i)}{N})\Big)^2+\Big(e_{-k}(\frac{i}{M})-e_{-k}(\frac{\tau^l(i)}{N})\Big)^2=8\sin^2\big(k\pi(\frac{i}{M}-\frac{\tau^l(i)}{N})\big).
		\end{equation}
		Therefore, combing with \eqref{cod5}, we get
		\begin{equation}\label{equs2}
			\begin{split}
				S_2&\leq \sum_{i=1}^{M}\frac{4\pi^2}{MN^2}\sum_{k=1}^{N}\frac{1}{|k|^{2\gamma-2}} +\frac{1}{2M}\sum_{i=1}^{M}\sum_{k=N+1}^{M}\frac{1}{|k|^{2\gamma}}\\
				&\leq CN^{1-2\gamma}.
			\end{split}
		\end{equation}

		\textbf{Estimate on $<M>_t$:} Note that by \eqref{cod5}, \eqref{iden3} and Cauchy inequality, we have
		\begin{equation*}
			\begin{split}
				&\frac{d}{dt}<M>_t\\
				&\leq\sum_{i=1}^{M}\frac{(a^1(i))^2}{M}(X_M^i(t)-X_N^{\tau^1(i)}(t))^2\Big(\sum_{k=-N}^{N}\frac{1}{|k|^{2\gamma}}\big(e_k(\frac{i}{M})-e_k(\frac{\tau^1(i)}{N})\big)^2 +\sum_{|k|=N+1}^{M}\frac{1}{|k|^{2\gamma}}e_k^2(\frac{i}{M})\Big)\\
				&+\sum_{i=1}^{M}\frac{(a^2(i))^2}{M}(X_M^i(t)-X_N^{\tau^2(i)}(t))^2\Big(\sum_{k=-N}^{N}\frac{1}{|k|^{2\gamma}}\big(e_k(\frac{i}{M})-e_k(\frac{\tau^2(i)}{N})\big)^2 +\sum_{|k|=N+1}^{M}\frac{1}{|k|^{2\gamma}}e_k^2(\frac{i}{M})\Big) \\
				&\leq \sum_{i=1}^{M}\frac{a^1(i)}{M}(X_M^i(t)-X_N^{\tau^1(i)}(t))^2\Big(\sum_{k=1}^{N}\frac{8\pi^2}{k^{2\gamma-2}N^2}+\sum_{k=N+1}^{\infty}\frac{1}{k^{2\gamma}}\Big)  \ \ \ \ \ \ \ \ \ \ \ \ \ \ \ \ \ \ \ \ \ \ \ \ \ \ \ \ \ \ \ \  \\
				&+\sum_{i=1}^{M}\frac{a^2(i)}{M}(X_M^i(t)-X_N^{\tau^2(i)}(t))^2\Big(\sum_{k=1}^{N}\frac{8\pi^2}{k^{2\gamma-2}N^2}+\sum_{k=N}^{\infty}\frac{1}{k^{2\gamma}}\Big).
			\end{split}
		\end{equation*}
It yields 
	\begin{equation}\label{equs4}
	\frac{d}{dt}<M>_t\leq C\bar{W}_2^2(\bar{L}_N(t), \bar{L}_M(t))\cdot(N^{1-2\gamma}+N^{1-2\gamma})	\leq \frac{C}{N^{2\gamma-1}}\bar{W}_2^2(\bar{L}_N(t), \bar{L}_M(t)).
	\end{equation}

	\textbf{Estimate on $S_1$:} 
	By Lemma \ref{lem3}, we have
	\begin{equation}\label{equs3}
		S_1(t,\omega)< \frac{2\beta}{N}+\frac{2\beta}{M},  \ \ \  \forall \, (t,\omega)\in (0,\infty)\times \Omega.
	\end{equation}
		
		\quad We denote $S_1+S_2$ as $S_t$; we have, for $\forall \, t>0$, that
		$$\max_{s\in[0,t]}\bar{W}_2^2(\bar{L}_N(s), \bar{L}_M(s)) \leq \bar{W}_2^2(\bar{L}_N(0), \bar{L}_M(0)) +\int_{0}^{t}(S_r)_{+} dr +\max_{s\in[0,t]}|M_s|.$$ 
		
		We then estimate $\E[\max_{s\in[0,t]}\bar{W}_2^2(\bar{L}_N(s), \bar{L}_M(s)) ]$. Combing with the estimates \eqref{equs2}, \eqref{equs3}, \eqref{equs4} and Burkholder-Davis-Gundy inequality, we obtain
		\begin{equation}\label{est33.5}
			\begin{split}
				&\E[\max_{s\in [0,t]} \bar{W}_2^2(\bar{L}_N(s), \bar{L}_M(s)) ] \\
				&\leq C_0^N+N^{1-2\gamma}t+4\beta N^{-1} t+CN^{\frac{1}{2}(1-2\gamma)}\E[\Big(\int_0^t\bar{W}_2^2(\bar{L}_N(s), \bar{L}_M(s))ds\Big)^\frac{1}{2}]\\
				&\leq C_0^N+C'N^{-1}t+CN^{\frac{1}{2}(1-2\gamma)}t^{\frac{1}{2}} \E[\max_{s\in [0,t]} \bar{W}_2(\bar{L}_N(s), \bar{L}_M(s)) ],
			\end{split}
		\end{equation}
		
		where 
		$$C^N_0=4\bar{W}_2^2(\bar{L}_N(0),\mu_0)\geq \bar{W}_2^2(\bar{L}_N(0),\bar{L}_M(0)),$$

		when $N$ is large enough. Note that the condition $\lim_{N\to \infty}\bar{W}_2(\bar{L}_N(0),\mu_0)=0$ implies that $\lim\limits_{N\to \infty}C_0^N\to 0$, therefore for small $\eps>0$, we can choose $N$ large enough such that 
		\begin{equation}\label{con3}
			C_0^N<\frac{1}{2}\eps, \ \ \ \ 2C^2N^{1-2\gamma}\eps^{-1}+4C'N^{-1}\eps^{-1}<\frac{1}{2}\eps.
		\end{equation}
		
		According to \eqref{est33.5}, we then obtain by Cauchy inequality that
		\begin{equation}\label{est34}
			\begin{split}
				\E[\max_{s\in [0,t]} \bar{W}_2^2(\bar{L}_N(s), \bar{L}_M(s)) ]\leq \eps, \ \ \ \ \ \ \forall   \, t<\frac{1}{\eps}.
			\end{split}
		\end{equation}
		
		Since $L_N(t)$ and $L_M(t)$ are the projections of $\bar{L}_N(t)$ and $\bar{L}_M(t)$ on the torus respectively, it holds
		\begin{equation}\label{obv2}
			W_2^2(L_N(t), L_M(t))\leq \bar{W}_2^2(\bar{L}_N(t),\bar{L}_M(t)).
		\end{equation}
		It follows from \eqref{est34}, \eqref{obv2} and H\"older inequality that
		\begin{equation}\label{obv1}
			\begin{split}
				\mathbb{W}^2_2(P^N,P^M)&\leq \E\left[ \Big(\sum_{n=1}^{\infty}\frac{1}{2^n}\max_{t\in[0,n]}W_2(L_N(t), L_M(t))\wedge 1\Big)^2\right]\\
				&\leq \sum_{n=1}^{\infty}\frac{1}{2^n}\E[\big|\max_{t\in[0,n]}W_2(L_N(t), L_M(t))\wedge 1\big|^2]\\
				&\leq \sum_{n=\frac{1}{\eps}+1}^\infty \frac{1}{2^n}+\E[\max_{t\in[0,\frac{1}{\eps}]}W_2^2(L_N(t), L_M(t))]\\
				&\leq C'\eps^2+C\eps< C'' \eps.
			\end{split}
		\end{equation}
		Therefore, $\{P^N\}$ is a Cauchy sequence in $\big(\P(\C([0,\infty), \P(\T))), \W_2\big)$.

		\quad Furthermore, we can arrange the initial positions for particles such that $C^N_0\leq \frac{4}{N}$. Indeed, if $\mu_0$ does not having atoms, we construct
		\begin{equation}\label{inivalue}
			X_N^i(0)=F_{\mu_0}^{-1}(\frac{i}{N}). 
		\end{equation}
		Then the initial empirical measure $\bar{L}_N(0)=\frac{1}{N}\sum_{i=1}^{N}\delta_{\{X_N^i(0)\}}$ weakly converges to $\mu_0$ and
		\begin{equation}\label{est13}
			\begin{split}
				\bar{W}_2^2(L_N(0),\mu_0)&\leq\sum_{i=1}^{N}\int_{F_{\mu_0}^{-1}(\frac{i-1}{N})}^{F_{\mu_0}^{-1}(\frac{i}{N})}(F_{\mu_0}^{-1}(\frac{i}{N})-x )^2 d\mu_0(x) \\
				&\leq \sum_{i=1}^{N}\int_{\frac{i-1}{N}}^{\frac{i}{N}}(F_{\mu_0}^{-1}(\frac{i}{N})-F_{\mu_0}^{-1}(x))^2 dx\\
				&\leq \sum_{i=1}^{N}\frac{1}{N}\Big[F_{\mu_0}^{-1}(\frac{i}{N})-F_{\mu_0}^{-1}(\frac{i-1}{N})\Big]^2\\
				&\leq \frac{1}{N}\sum_{i=1}^{N}\big(F_{\mu_0}^{-1}(\frac{i}{N})-F_{\mu_0}^{-1}(\frac{i-1}{N})\big)= \frac{1}{N}.
			\end{split}
		\end{equation}
		If $\mu_0$ is a discrete measure, for instance $\mu_0=\delta_{0}$, then we can construct 
		\begin{equation*}\label{ln0}
			\bar{L}_{N}(0)=\frac{1}{N} \sum_{j=1}^{N}\delta_{\frac{j}{N^2}}.
		\end{equation*}
		Thus, a similar estimate as \eqref{est13} gives
		\begin{equation}\label{est13.5}
			\bar{W}_2^2(L_N(0), \mu_0)\leq \frac{C}{N^2}.
		\end{equation}
		Since any measure $\mu_0$ can be uniquely decomposed as a sum of a discrete measure and a measure without atoms, we obtain by \eqref{est13} and \eqref{est13.5} that 
		\begin{equation}\label{ini}
		C^N_0\leq \frac{4}{N}.
		\end{equation}
	
		Combining with \eqref{equs2}, \eqref{equs3}, \eqref{obv2} and \eqref{ini}, we have
		\begin{equation*}\label{est12}
			\E[W_2^2(L_N(t), L_M(t))]\leq \frac{4\beta t}{N}+\frac{4}{N},
		\end{equation*}
	when $N$ is large enough. \eqref{prop7a} is proved.
		Note that taking $C_0^N<\frac{4}{N}$ and choosing $N$ such that $N^{-\frac{1}{2}}<\eps$, then \eqref{con3} and \eqref{est34} show that for $t<N^{\frac{1}{2}}$, 
		\begin{equation*}
			\E[\max_{s\in[0,t]}W_2^2(L_N(t), L_M(t))]\leq CN^{-\frac{1}{2}}.
		\end{equation*}
		Plugging this estimate into \eqref{obv1} gives
		$$\W_2^2(P^N, P^M)<CN^{-\frac{1}{2}}.$$

		\eqref{prop7b} is proved. 
		The proof is finished.
	\end{proof}
	\vskip3mm

	\quad Note that $\C([0,\infty), \P(\T))$ equipped with $d$ is a complete separable metric space, therefore $\P(\C([0,\infty), \P(\T))$ equipped with $\W_2$ is a Polish space \cite{V03}. There exists a unique limit $P$ such that $\mathbb{W}_2(P^N, P)\to 0$. Furthermore, By Skorohod representation theorem and uniqueness of strong solution to SDEs \eqref{Napprox}, we can find a new probability space $(\tilde{\Omega},\tilde{\mathcal{F}}, \{\tilde{\mathcal{F}}_t\}_{t\geq 0},\tilde{\P})$ with processes $\displaystyle\{\tilde{X}_N(t), \, t\geq 0\}_{N\in \N^+}$ and independent Brownian motions $\{\tilde{W}^k_t\}_{k\in \Z}$ defined on it, such that for $\forall \, N>0$, 
	$$(\{X_N(t)\}, W_t^1,W_t^{-1},\ldots,W_t^{N},W_t^{-N}) \stackrel{dist}{=}(\{\tilde{X}_N(t)\}, \tilde{W}_t^1,\tilde{W}_t^{-1},\ldots,\tilde{W}_t^{N},\tilde{W}_t^{-N}),$$
	where $\tilde{X}_N(t)$ is the unique solution to 
	\begin{equation}\label{model2}
		d \tilde{X}_N^i(t)=\frac{\beta}{N}\sum^N_{j=1, j \neq i}\cot\big(\pi(\tilde{X}_N^i(t)-\tilde{X}_N^j(t))\big) d t +  \sum_{k=-N}^{N}\frac{1}{|k|^{\gamma}}e_k(\frac{i}{N})d\tilde{W}_t^k .
	\end{equation}
	Meanwhile, there is a limit process $p_t : \, (\tilde{\Omega},\tilde{\mathcal{F}}, \{\tilde{\mathcal{F}}_t\}_{t\geq 0},\tilde{\P})\to \C([0,\infty),\P(\T))$ with the distribution $P$ such that $\displaystyle \tilde{L}_N(t)=\frac{1}{N}\sum_{i=1}^{N}\delta_{\{\tilde{X}^i_N(t)\}}$ converges to $p_t$ in $\C([0,\infty),\P(\T))$ under metric $d$ almost surely. For simplicity, we also denote $\tilde{L}_N(t)$ as $p^N_t$. 
	
	\quad We point out in the next proposition that the choice of the initial empirical measure $L_N(0)$ actually has no influence on the limit distribution of the model, as long as they satisfy the admissible conditions \eqref{cod2'}. 
	\begin{proposition}\label{prop2}
		Given model parameters $\gamma>\frac{3}{2}$ and $\beta>\lambda_0^\gamma$, then for any two distinct initial empirical measures $\{L^1_N(0)\}, \{L^2_N(0)\}$ satisfying the admissibility conditions \eqref{cod2'}, the corresponding limit distributions $P_1$ and $P_2$ exists and 
		$$\W_2(P_1,P_2)=0.$$ 
	\end{proposition}
	
	\quad  In fact, let $\bar{L}_N^{L_N^1(0)}(t)$ and $\bar{L}_N^{L_N^2(0)}(t)$ be the two empirical measure processes induced by different initial measures $\{L^1_N(0)\},\{L^2_N(0)\}$. The situation here is more simpler than that in Theorem \ref{ratecon} since the optimal transport plan from $\bar{L}_N^{L_N^1(0)}(t)$ to $\bar{L}_N^{L_N^2(0)}(t)$ is a one-to-one map. One can derive by mimicking the proof above that the nonnegative semimartingale $Y_t^N=\bar{W}_2^2(\bar{L}_N^{L_N^1(0)}(t), \bar{L}_N^{L_N^2(0)}(t))$ satisfies
	$$dY_t^N=S_t^N dt+\sigma_t^N dW_t,$$
	where 
	$$0<S_t^N<N^{\kappa'},$$
	and 
	$$<Y^N>_t\leq N^{\kappa'} \int_{0}^{t}Y_s^N ds,$$
	for some negative constant $\kappa'$ depending on $\gamma$. Also, $C_0^N$ goes to $0$ as $N\to \infty$ due to the admissibility condition \eqref{cod2'}. We omit the complete proof since it is almost the same to the proof of Theorem \ref{ratecon} except a little modifications. This proposition implies that once the macroscopic parameters $\beta$, $\gamma$ and an appropriate initial measure $\mu_0$ are given, one always get a unique limit distribution $P$ no matter how you design the admissible initial empirical measure of the approximation model. 
	
	\vskip 3mm
	
		\quad Due to Theorem \ref{ratecon}, there always exists a unique fluctuating hydrodynamic limit of the particle model \eqref{model1} for $\gamma>\frac{3}{2}$ and $\beta>\lambda_0^\gamma$. As mentioned in the introduction, we denote the corresponding distribution of $p^\beta_t$ in $\C([0,\infty), \P(\T))$ as $P_{\beta}$. There is no confusion for these notations due to Proposition \ref{prop2}. We will prove that the curve $\beta \in (\lambda_0^\gamma, \infty) \, \to  \, P_\beta$ in $\P(\C([0,\infty), \P(\T)))$ is at least a continuous curve.
			\begin{proposition}\label{prop5.1}
			Suppose that $\gamma >\frac{3}{2}$, then the limit distributions $\{P_{\beta}\}_{\beta>\lambda_0^\gamma}$ of the model \eqref{model1} satisfies
			\begin{equation*}\label{ident1}
				\W_2(P_{\beta'}, P_\beta)\leq C\sqrt{|\beta'-\beta|}, \ \ \ \ \ \forall \, \beta, \beta'>\lambda_0^\gamma,
			\end{equation*}
			for positive number $C$ depending on $\beta$ and $\beta'$.
		\end{proposition}
		\begin{proof}
			We use the argument in the proof of Theorem \ref{ratecon} and a coupling method. Due to Proposition \ref{prop2}, we can take $L_N(0)$ satisfying \eqref{inivalue} without loss of generality. Under the identical choice of probability space $(\Omega, \F, \{\F_t\}_{t\geq 0},\P)$ and Brownian motions $\{B_t^k\}_{k\in \Z}$, we denote the solutions to \eqref{model1} and the associated empirical measure process respectively as $X^{\beta}_N(t)$ and $\bar{L}^\beta_N(t)$ for $\beta>\lambda_0^\gamma$ ($X^{\beta'}_N(t)$ and $\bar{L}^{\beta'}_N(t)$ for $\beta'>\lambda_0^\gamma$). Without loss of generality, we assume $\beta'>\beta$. According to Theorem \ref{introthm1}, both $\{X^{\beta,i}_N(t)\}_{i=1}^{N}$ and $\{X^{\beta',i}_N(t)\}_{i=1}^{N}$ keep orders, and thus we have
			\begin{equation*}
				\bar{W}_2^2(\bar{L}^\beta_N(t), \bar{L}^{\beta'}_N(t))=\frac{1}{N}\sum_{i=1}^{N}(X_N^{\beta,i}(t)-X_N^{\beta',i}(t))^2
			\end{equation*}
			Subsequently, by It\^o's formula,
			\begin{equation*}
				\begin{split}
					&d\bar{W}_2^2(\bar{L}^\beta_N(t), \bar{L}^{\beta'}_N(t))\\
					&=\sum_{i=1}^{N}\frac{1}{N}(X_N^{\beta,i}(t)-X_N^{\beta',i}(t))\circ \big(d X_N^{\beta,i}(t)-dX_N^{\beta',i}(t)\big)  \\
					&=\sum_{i=1}^{N}\frac{1}{N}(X_N^{\beta,i}(t)-X_N^{\beta',i}(t))\cdot \Big(\frac{\beta}{N}\sum_{j=1, j\neq i}^{N}\cot(\pi(X_N^{\beta,i}(t)-X_N^{\beta,j}(t))) \Big) dt\\
					&-\sum_{i=1}^{N}\frac{1}{N}(X_N^{\beta,i}(t)-X_N^{\beta',i}(t))\cdot \Big(\frac{\beta'}{N}\sum_{j=1, j\neq i}^{N}\cot(\pi(X_N^{\beta',i}(t)-X_N^{\beta',j}(t))) \Big) dt\\
					&+\sum_{i=1}^{N}\frac{1}{2N}\sum_{k=-N}^{N}\frac{1}{|k|^{2\gamma}}\Big(e_k(\frac{i}{N})-e_k(\frac{i}{N})\Big)^2 dt \\
					&+\sum_{i=1}^{N}\frac{1}{N}(X_N^{\beta,i}(t)-X_N^{\beta',i}(t)) \cdot\sum_{k=-N}^{N}\frac{1}{|k|^{\gamma}}\Big(e_k(\frac{i}{N})-e_k(\frac{i}{N})\Big)dW_t^k.\\
				\end{split}
			\end{equation*}
			Note that the sum of the last two lines equals $0$, and hence $\bar{W}_2^2(\bar{L}^\beta_N(t), \bar{L}^{\beta'}_N(t))$ is  a continuous process with bounded variation. Analogous to \eqref{equmn} and \eqref{equnm} (letting $\lambda\to 0$), we obtain 
			\begin{equation*}
				\beta(\bar{\textbf{W}}(\bar{L}^{\beta'}_N(t))-\bar{\textbf{W}}(\bar{L}^\beta_N(t)))\geq \sum_{i=1}^{N}\frac{1}{N}(X_N^{\beta,i}(t)-X_N^{\beta',i}(t))\cdot \Big(\frac{\beta}{N}\sum_{j=1, j\neq i}^{N}\cot(\pi(X_N^{\beta,i}(t)-X_N^{\beta,j}(t))) \Big);
			\end{equation*}
			and 
			\begin{equation*}
				\beta'(\bar{\textbf{W}}(\bar{L}^\beta_N(t))-\bar{\textbf{W}}(\bar{L}^{\beta'}_N(t)))\geq \sum_{i=1}^{N}\frac{1}{N}(X_N^{\beta',i}(t)-X_N^{\beta,i}(t))\cdot \Big(\frac{\beta'}{N}\sum_{j=1, j\neq i}^{N}\cot(\pi(X_N^{\beta',i}(t)-X_N^{\beta',j}(t))) \Big).
			\end{equation*}
			Therefore, 
			\begin{equation*}
				\frac{d}{dt}\bar{W}_2^2(\bar{L}^\beta_N(t), \bar{L}^{\beta'}_N(t))\leq 	\beta\big(\bar{\textbf{W}}(\bar{L}^{\beta'}_N(t))-\bar{\textbf{W}}(\bar{L}^\beta_N(t))\big)+	\beta'\big(\bar{\textbf{W}}(\bar{L}^\beta_N(t))-\bar{\textbf{W}}(\bar{L}^{\beta'}_N(t))\big).
			\end{equation*}
			Since $\bar{L}_N^\beta(0)=\bar{L}_N^{\beta'}(0)=L_N(0)$ and $\textup{\textbf{W}}(\bar{L}_N^\beta(t))-C_{\gamma,\beta}t$, $\textup{\textbf{W}}(\bar{L}_N^{\beta'}(t))-C_{\gamma,\beta'}t$ are supermartingales according to Theorem \ref{introthm1}, we have 
			$$\left|\E\left[\bar{\textbf{W}}(\bar{L}^{\beta'}_N(t))-\bar{\textbf{W}}(\bar{L}^\beta_N(t))\right]\right|\leq Ct. $$
			It follows that 
			$$\E[\bar{W}_2^2(\bar{L}^\beta_N(t), \bar{L}^{\beta'}_N(t))]\leq \frac{C}{2}|\beta'-\beta|t^2.$$
			Hence, by \eqref{obv2} and \eqref{obv1}, we obtain 
			\begin{equation}\label{est30}
				\W_2^2(P^N_\beta, P^N_{\beta'})\leq \sum_{n=1}^{\infty}\frac{1}{2^n}\E[\max_{t\in[0,n]}\bar{W}_2^2(\bar{L}^{\beta'}_N(t),\bar{L}^\beta_N(t))]\leq C|\beta'-\beta|,
			\end{equation}
			where $P^N_{\beta'}$ and $P^N_\beta$ are the distributions of $\bar{L}_N^{\beta'}(t)$ and $\bar{L}_N^\beta(t)$ respectively. By \eqref{prop7b}, $\W_2^2(P_\beta,P_\beta^N)\to 0$ and $\W_2^2(P_{\beta'},P_{\beta'}^N)\to 0$ as $N\to \infty$. It then follows from \eqref{est30} that 
			$$\W_2^2(P_{\beta'}, P_\beta)\leq C'|\beta'-\beta|,$$
			for some constant $C'$. The proof is finished.
		\end{proof}
	\vskip5mm

	\section{Limit equation and regularity estimates}\label{sec6}
In this section, we will show that \eqref{model1} is indeed an approximation particle model of the Dean-Kawasaki type equation \eqref{withrepul} and study the pathwise regularity of the fluctuating hydrodynamic limit.

	\subsection{Limit equation of $p^\beta_t$ for $\beta>\lambda_0^\gamma$}\label{sec5}
	\quad In this subsection, we prove that when $\beta>\lambda_0^\gamma$, the fluctuating hydrodynamic limit $p^\beta_t$ of the particle model is a non-trivial solution to $(MP)^{K_2^\gamma, \beta,\gamma}_{\mu_0}$ associated with $(\tilde{\Omega},\tilde{\mathcal{F}}, \tilde{\P})$ and $\{\tilde{W}^k_t\}_{k\in \Z}$, if the initial measure $\mu_0$ satisfies the condition $\textbf{(J)}$. As a necessary step for this convergence, we firstly prove there is no atom in $p^\beta_t$ almost surely. Since the admissible parameter $\beta$ has no influence on the result below, we denote the limit measure valued process as $p_t$ with its distribution $P$ for simplicity of notations at the moment. 
	
	\begin{lemma}\label{noatom}
		Suppose that $ \gamma>\frac{3}{2}$, $\beta>\lambda_0^\gamma$ and $\mu_0$ satisfies condition \textbf{\textup{(J)}}. If the particle model satisfies the conditions \eqref{cod2'},
		then for $\tilde{\P}-$a.s., the limit measure valued process $p^\beta_t$ is non-atomic for all $t\geq 0$. 
	\end{lemma}
	
	\begin{proof}
	Since the limit distribution $P$ does not depend on the initial empirical measures $\{L_N(0)\}$ by Proposition \ref{prop2}, we pick $\{L_N(0)\}$ satisfying \eqref{inivalue} without loss of generality. Let $U^T=\{\omega:\, \exists \, t\in [0,T], \, \mathrm{such \ that} \ p_t(\omega) \mathrm{ \ is \ atomic}\}$. 	To prove the lemma, it suffices to prove for $\forall \, T>0$, $\tilde{\P}(U^T)=0$. If  the measurable set $U^T$ has positive measure, i.e. $\tilde{\P}(U^T)=C>0$ . We define
		\begin{equation*}
			U^T_i=\{\omega\in U^T \, : \exists \, t, x,\  \mathrm{such \ that}\  p_t(\omega,dx)=\eta \delta_x \ \mathrm{with} \ \eta>\frac{1}{i}\} ,
		\end{equation*}
		then it is obvious that $U^T_i\subset U^T_{i+1}$ and $\displaystyle \bigcup_{i=1}^{\infty}U^T_i=U^T$. Thus we can find some $k\in \N$  such that $\tilde{\P}(U^T_k)>\frac{C}{2}$. 
		Let $E_x^N=(x-\frac{1}{2N},x+\frac{1}{2N})$. Note that for $\tilde{\P}-$ a.s. $\omega$, the curve $p^n_t(\omega,\cdot)$ weakly converges to $p_t(\omega,\cdot)$ uniformly in $t\in[0,T]$. If $p_t(\omega,dx)=\eta \delta_x dx$ for $t\in [0,T]$, then for each $N$, there exists $\mathcal{N}(N,\omega,t,x)$ such that $\forall \, n\geq \mathcal{N}(N,\omega,t,x)$,  
		\begin{equation*}
			\int_{E_x^N}p_t^n(\omega, dy)>\frac{\eta}{2} .
		\end{equation*}
		Based on this observation, we define 
		\begin{equation*}
			U_k^{n,N,T}=\{\omega:\exists \, t\in [0,T], \,  x\in \T , \ \mathrm{such \ that \ for } \, \forall \,j\geq n , \int_{E_x^N}p^j_t(\omega, dy) >\frac{1}{2k}\} ,
		\end{equation*}
		then we must have $U_k^{n,N,T}\subset U_k^{n+1,N,T}$ and $\displaystyle U^T_k\subset\bigcup_{n=1}^{\infty}U_k^{n,N,T}$. Thus for each $N$, there exists $m_N$ such that $\tilde{\P}(U_k^{m_N,N, T})>\frac{C}{3}=C'$.  Now let 
		\begin{equation*}
			\bar{U}_k^{m_N,N,T}=\{\omega:\exists \, t\in [0,T], \, x\in \T, \ \mathrm{such \ that \ for } \, \forall \,j\geq m_N , \int_{E_x^N}L_j(\omega,t ,dy) >\frac{1}{2k}\}.
		\end{equation*}
		Remember that $L_N(t)$ has the same distribution with $p^N_t$. We must have 
		\begin{equation*}
			\P(\bar{U}_k^{m_N,N,T})=\tilde{\P}(U_k^{m_N,N,T})=C'>0.
		\end{equation*}
		On the other hand, we define  a stopping time 
		\begin{equation*}
			\tau_{m_N}^k:=\inf\{t: \ \min_{j}|e^{2\pi iX_{m_N}^{j}(t)}- e^{2\pi iX_{m_N}^{j+\frac{m_N}{2k}} (t)}|\leq \frac{1}{2N}\}.
		\end{equation*}
		Denote 
		\begin{equation*}
			A_T=\{\tau_{m_N}^k\leq T\} ,
		\end{equation*}
		then, according to Theorem \ref{introthm1}, we have
		\begin{equation*}
			\begin{split}
				&\textbf{W}(L_{m_N}(0))+\frac{\pi}{2}(C_\gamma+ \frac{\beta }{6}) \tau_{m_N}^k\land T\\
				&\geq \E[\textbf{W}(L_{m_N}(\tau_{m_N}^k\land T)]\\
				&=\E[\textbf{W}(L_{m_N}(\tau_{m_N}^k))\mathbbm{1}_{A_T}]+\E[\textbf{W}(L_{m_N}(T))\mathbbm{1}_{A_T^{c}}]\\
				&\geq -\frac{1}{2m_N^2}\frac{m_N}{2k}(\frac{m_N}{2k}-1)\ln(\frac{1}{2N})\P(A_T)\\
				&-\frac{1}{2m_N^2}(m_N^2-m_N-\frac{m_N}{2k}(\frac{m_N}{2k}-1))\log2\P(A_T)\\
				&-\frac{1}{2m_N^2}(m^2_N-m_N)\log2\P(A_T^{c})\\
				&=\frac{1}{2m_N^2}\frac{m_N}{2k}(\frac{m_N}{2k}-1)(\ln (2N)+\log2)\P(A_T)-\frac{m_N-1}{2m_N}\log2 .
			\end{split}
		\end{equation*}
		It follows that 
		\begin{equation}\label{estpt}
			\P(\tau_{m_N}^k\leq  T)\leq \frac{\textbf{W}(L_{m_N}(0))+\frac{\pi}{2}(C_\gamma+ \frac{\beta }{6}) T+\log2}{\frac{\frac{m_N}{2k}-1}{4k m_N}(\ln(2N)+\log2)}\leq \frac{10k^2(C_T+2\textbf{W}(L_{m_N}(0)))}{\ln(2N)+\log2}.
		\end{equation}
		
		Now, we need to estimate $\textbf{W}(L_{m_N}(0))$. Under the condition $\textbf{\textup{(J)}}$, we have
		\begin{equation*}
			\begin{split}
				&\bar{\textbf{W}}(L_N(0))\\
				&=-\frac{1}{N^2}\sum_{i=1}^{N}\sum_{j=1, j\neq i}^{N}\ln\big|e^{2\pi i F_{\mu_0}^{-1}(\frac{i}{N})}-e^{2\pi i F_{\mu_0}^{-1}(\frac{j}{N})}\big| \\
				&=-\frac{1}{N^2}\sum_{i=1}^{N}\sum_{j=1, j\neq i}^{N}\ln\Big(2\sin\left(\pi \Big[\big| F_{\mu_0}^{-1}(\frac{i}{N})-F_{\mu_0}^{-1}(\frac{j}{N})\big|\wedge\big(1-\big| F_{\mu_0}^{-1}(\frac{i}{N})-F_{\mu_0}^{-1}(\frac{j}{N})\big|\big)\Big]\right)\Big). 
			\end{split}
		\end{equation*}
		Since $-\ln|\sin(x\pi)|< -\ln|x|$ for $|x|\leq\frac{1}{2}$ and 
		$$\big|F_{\mu_0}^{-1}(x)-F_{\mu_0}^{-1}(y)\big|\geq \frac{1}{||\rho_0||_{\L^\infty}}|x-y|,$$
	we have
		$$\bar{\textbf{W}}(L_N(0))\leq \ln(2||\rho_0||_{\L^\infty})-\frac{2}{N^2}\sum_{i=1}^{N}\sum_{j=1, j\neq i}^{N} \ln \Big( |\frac{i-j}{N}|\wedge(1-|\frac{i-j}{N}|)\Big).$$
		Note that, for $i\in \{1,2,\ldots, N\}$,
		$$-\frac{1}{N}\sum_{j=1, j\neq i}^{N} \ln \Big( |\frac{i-j}{N}|\wedge(1-|\frac{i-j}{N}|)\Big)<2\int_{0}^{1}-\ln x dx=2.$$
		Therefore,
		\begin{equation}\label{estw}
			\bar{\textbf{W}}(L_{N}(0))<\ln(2||\rho_0||_{\L^\infty})+2 .
		\end{equation}
	Plugging \eqref{estw} into \eqref{estpt}, we can choose $N_T$ such that 
		\begin{equation*}
			\P(\tau_{m_{N_T}}^k\leq T)<\frac{1}{2}C'.
		\end{equation*} 
		However, for each $\omega \in \bar{U}_k^{m_{N_T},N_T}$, there must be at least $\frac{m_{N_T}}{2k}$ particles included in some interval $(x-\frac{1}{2N_T}, x+\frac{1}{2N_T})$, which means $\tau_{m_{N_T}}^k(\omega )\leq T$. Therefore , $$\P(\tau_{m_{N_T}}^k\leq T)\geq \P(\bar{U}_k^{m_{N_T},N_T})=C'.$$ Contradiction! We finished the proof.
	\end{proof}
	\begin{remark}
		We remark that condition \textbf{\textup{(J)}} is a technical condition to make sure the boundedness of  $\textbf{W}(L_{N}(0))$.
	\end{remark}
	\vskip4mm
	As mentioned in the introduction, we define an auxiliary $\P_{ac}(\T)-$valued process on $(\tilde\Omega, \tilde\F, \tilde\P)$ associated with the convergent sequence
	$\displaystyle p^N_t=\frac{1}{N}\sum_{i=1}^{N}\delta_{\{\tilde{X}_N^i(t)\}}$. Let
	$$q_N(t)=\sum_{i=1}^{N}\frac{\frac{1}{N}}{|\tilde{X}_N^i(t)- \tilde{X}_N^{i-1}(t)|}\mathbbm{1}_{[\tilde{X}_N^{i-1}(t),\tilde{X}_N^i(t))}.$$
	Here, we represent the measure by its density $q_N$ without confusion. We will show that when $\beta>\frac{1}{2}\pi K_1^\gamma$, the fluctuating hydrodynamic limit $p^\beta_t$ has a density with finite Boltzmann entropy almost surely for any time $t\geq 0$ when $\beta>\frac{1}{2}\pi K_1^\gamma$. 
	
	\begin{proposition}\label{abs}
		Suppose that $\gamma>\frac{3}{2}, \, \beta>\frac{1}{2}\pi K_1^\gamma$. If $\mu_0$ satisfies the condition $\textbf{\textup{(J)}}$ and the initial empirical measure of the particle model \eqref{model1} satisfies the condition \eqref{cod2'}, then the fluctuating hydrodynamic limit $p^\beta_t$ is absolutely continuous w.r.t $Leb$ for all $t\geq 0$, $\tilde{\P}-$a.s., and there exists constant $C_0$ such that \eqref{entineq} holds.
		
	\end{proposition}

\begin{proof}
	We pick the initial value $L_N(0)$ satisfying \eqref{inivalue} without loss of generality. Under the condition \textbf{\textup{(J)}}, we have 
	\begin{equation*}\label{cod5.3}
		S^N(L_N(0))=\frac{1}{N}\sum_{j=1}^{N}\ln \frac{\pi}{N\sin (\pi(X_N^j(t)-X_N^{j-1}(t)))}\leq \frac{1}{N}\sum_{j=1}^{N}\ln ||F'_{\mu_0}||_{\L^\infty}=\ln||\rho_0||_{\L^\infty},
	\end{equation*} 
	Thus, combing with \eqref{ugr'}, we obtain
	\begin{equation}\label{estinv}
		\E[S^N(L_N(t))]\leq e^{-2 C_1^\beta t}C_0+C,
	\end{equation}
	for constants $ C_0=\ln ||\rho_0||_{\L^\infty}$, $C_1^\beta=\frac{\beta-\frac{1}{2}\pi K_1^\gamma}{\pi}$ and $C=\frac{\beta \pi}{C_1^\beta}$. 
	
	\quad On the other hand, by the inequality 
	$$|e^{2\pi i \tilde{X}_N^i(t)}-e^{2\pi i \tilde{X}_N^{i-1}(t)}|\leq 2\pi|\tilde{X}_N^i(t)- \tilde{X}_N^{i-1}(t)|,$$
	we have
	\begin{equation}\label{est7}
		\textbf{S}(q_N(t))\leq S^N(p^N_t).
	\end{equation}
	Since we have 
	\begin{equation*}
		W_2^2(p^N(t), q_N(t))\leq \frac{1}{N}\sum_{i=1}^{N}|\tilde{X}_N^i(t)- \tilde{X}_N^{i-1}(t)|^2 \leq \frac{1}{N},
	\end{equation*}
	for $\tilde{\P}-a.s.$, $q_N(t)$ and $p^N(t)$ weakly converges to the same limit curve $p^\beta_t$ in $\C([0,\infty), \P(\T))$ under metric $d$. Thus, by the lower semicontinuity of $\textbf{S}$ on $\P(\T)$, we have for $\forall \,  t\geq 0$,
	\begin{equation}\label{est5.1}
		\textbf{S}(p_t^\beta)\leq \liminf_{N}\textbf{S}(q_N(t)),
	\end{equation}
	Since $p_t^N$ and $L_N(t)$ share the same distribution, we have
	\begin{equation}\label{est5.12}
		\tilde{\E}[S^N(p^N_t)]=\E[S^N(L_N(t))].
	\end{equation}
Then \eqref{entineq} follows from \eqref{estinv}, \eqref{est7}, \eqref{est5.1}, \eqref{est5.12} and Fatou's lemma. Also, we have, for any $k$ large enough, 
	\begin{equation}\label{BC}
		\P(S^N(L_N(t))>k^2)\leq \frac{e^{-2 C_1^\beta t}C_0+C}{k^2},
	\end{equation}
	Let 
	$$I_N^k=\left\{\omega\in \tilde{\Omega}: \textbf{S}(q_N(t,\omega))>k^2\right\},$$
	then by the lower semicontinuity of entropy, we know
	$$\left\{\omega\in \tilde{\Omega}: \textbf{S}(p^\beta_t(\omega))=\infty\right\}\subset \bigcap_{k=1}^{\infty}\bigcup_{m=1}^\infty\bigcap_{N=m}^{\infty}I_N^k.$$
	Also, notice again that $p_t^N$ and $L_N(t)$ share the same distribution, and thus by  \eqref{BC}, we have
	$$\tilde{\P}\Big(S^N(p^N_t)>k^2\Big)=\P(S^N(L_N(t))>k^2)\leq \frac{e^{-2 C_1^\beta t}C_0+C}{k^2}.$$
	Combining with \eqref{est7}, we have
	$$\tilde{\P}\Big(\textbf{S}(q_N(t))>k^2\Big)\leq \tilde{\P}\Big(S^N(p^N_t)>k^2\Big),$$
	so $\tilde{\P}\Big(\bigcap_{N=m}^{\infty}I_N^k\Big)\leq \frac{e^{-2 C_1^\beta t}C_0+C}{k^2}$. Because this estimate is independent of $m$, we get 
	$$\tilde{\P}\Big(\bigcup_{m=1}^\infty\bigcap_{N=m}^{\infty}I_N^k\Big)\leq \frac{e^{-2 C_1^\beta t}C_0+C}{k^2},$$
	which means 
	$$\tilde{\P}\left\{\omega\in \tilde{\Omega}: \textbf{S}(p^\beta_t(\omega))=\infty\right\}=0.$$
	$p^\beta_t$ must be absolutely continuous w.r.t $Leb$ in $\T$.	

\end{proof}

	\vskip 3mm
	
		\quad Recall that $p^N_t$ is the empirical measure of the projection of  $(\tilde{X}_N(t))\in \R^N$ satisfying \eqref{model2}, i.e. $\displaystyle p_t^N=\frac{1}{N}\sum_{i=1}^{N}\delta_{\{\tilde{X}_N^i(t)\}}$, and we have
		$$\int_0^1 \H p_t^N dp_t^N =\frac{1}{N^2}\sum_{i=1}^N\sum_{j=1, \, j\neq i}^N \cot\big(\pi(\tilde{X}_N^i(t)-\tilde{X}_N^j(t))\big).$$
		Next, we prove $p_t^N$ converges (in distribution) to the solution to $(MP)^{K_2^\gamma, \beta, \gamma}_{\mu_0}$.

	\quad For $\phi\in \C^2(\T)$, we define
	\begin{equation*}
		Q^{\gamma,N}_{\mu}(\phi)=\int_{0}^1\int_{0}^1\phi'(x)\phi'(y)\bigg(1+\sum_{k=0}^{N}\frac{2}{k^{2\gamma}}\cos\big(2\pi k (\int_{0}^1\mathbbm{1}_{(x\land y,x\lor y]} (z) \mu(dz)) \big ) \bigg)d\mu(x)  d\mu(y)  .
	\end{equation*}
and a continuous function $F_{\phi}:\, \T^2\to \R$ by
\begin{equation*}
	F_\phi(x,y)=\left\{
	\begin{aligned}
		&\frac{\phi'(x)-\phi'(y)}{\tan(\pi(x-y))}, \ \ \ \ \ x\neq y\\
		&\frac{\phi''(x)}{\pi}, \ \ \ \ \ \ x = y.
	\end{aligned}
	\right.
\end{equation*}
We now prove the theorem of this subsection.

	\begin{proof}[\textbf{Proof of Theorem \ref{maintheorem}}]
		For a solution $\mu_t$ to $(MP)^{K_2^\gamma, \beta, \gamma}_{\mu_0}$, the generator $\mathbb{L}$ associated with $\langle\mu_t, \phi\rangle$ is 
		
		\begin{equation*}
			\mathbb{L}_\beta f= \frac{\beta}{2} \int_{0}^1\int_{0}^{1}F_{\phi}d\mu\otimes \mu \cdot  f'+K_2^\gamma\langle\mu_t, \phi''\rangle f'+\frac{1}{2}Q^\gamma_{\mu_t}(\phi) f'',  \quad \quad \forall \, f\in \C^2(\R) .
		\end{equation*} 
		
		Thus, according to the equivalent description of $\P(\T)-$valued process (Lemma 7.2.1 in \cite{daw1993}), we only need to prove for $\forall \, \G\in \mathcal{D}:=\{\G: \, \G(\mu)=g(\langle \mu, \phi\rangle), \ \phi\in \C^2(\T), \ g\in \C^{2}(\R)\} $, 
		
		\begin{equation*}
			M^\G_t(p^\beta):=\G(p^\beta_t)-\G(\mu_0)-\int_{0}^{t}\mathrm{\mathbf{L}}_\beta \G(p^\beta_s) ds
		\end{equation*}
		
		is a $\tilde{\P}-$local martingale,	where 
		$$\mathrm{\mathbf{L}}_\beta \G(\mu)=\frac{\beta}{2} \int_{0}^1\int_{0}^{1}F_{\phi}d\mu\otimes \mu \cdot  g'(\langle \phi,\mu \rangle )+K_2^{\gamma}g'(\langle \phi,\mu \rangle )\langle \mu, \phi''\rangle+\frac{1}{2}g''(\langle \mu, \phi\rangle)Q^\gamma_{\mu}(\phi).$$ 
		This suffices to prove that, for $\forall \, \displaystyle 0\leq s<t$ and any continuous function $H: \C([0,\infty), \P(\T))\ra \R$,
		\begin{equation}\label{estgoal}
			\tilde{\E}\big[\big(\G(p^\beta_t)-\G(p^\beta_s)-\int_{s}^{t}\mathrm{\mathbf{L}}_\beta \G(p^\beta_r) dr\big)\cdot H(p^\beta\big|_{[0,s]}) \big]=0 .
		\end{equation} 
			Note that by \eqref{qbeta} we have
		\begin{equation*}
			\begin{split}
				Q^\gamma_{\mu} (\phi)
				&=\sum_{k=-\infty}^{+\infty}\frac{1}{|k|^{2\gamma}}\int_{0}^1\phi'(G_{\mu}(x))e_k(x)dx\int_{0}^1\phi'(G_{\mu}(y))e_k(y)dy
			\end{split}
		\end{equation*}
		On the other hand, since we have proved that, $\tilde{\P}-$a.s. $p^\beta_t(\omega)$ has no discrete part for all $t\geq 0$ (Proposition \ref{noatom} and Proposition \ref{abs}), for $ \tilde{\P}-$a.s., 
		$$G_{p_t^N}(x)\to G_{p^\beta_t}(x),  \ \ \ \forall \, x\in [0,1], \, t\in (0,\infty).  $$
	Therefore, by dominated convergence theorem, for $\forall \, k$,
	$$\int_{0}^1\phi'(G_{p_t^N}(x))e_k(x)dx \to \int_{0}^1\phi'(G_{p^\beta_t}(x))e_k(x)dx,$$
	as $N\to \infty$. Also, 
	$$\left|\sum_{|k|>N}\frac{1}{|k|^{2\gamma}}\int_{0}^1\phi'(G_{p_t^N}(x))e_k(x)dx\int_{0}^1\phi'(G_{p_t^N}(y))e_k(y)dy\right|<4N^{1-2\gamma}||\phi'||^2_{\L^\infty}.$$
	Hence, we have $Q^\gamma_{p^N_t} (\phi)\to Q^\gamma_{p^\beta_t} (\phi)$ for $\P-a.s.$. Now, combing with the weak convergence of $p^N_t$ to $p^\beta_t$ and dominated convergence theorem, we obtain
		\begin{equation}\label{EST5}
			\int_{0}^{t}\frac{1}{2}g''(\langle p_r^N, \phi\rangle)Q^\gamma_{p_r^N}(\phi) dr\rightarrow \int_{0}^{t}\frac{1}{2}g''(\langle p^\beta_r, \phi\rangle)Q^\gamma_{p^\beta_r}(\phi) dr.
		\end{equation}
		
		Since 
		$$\langle \phi'\H p_r^N, p_r^N\rangle =\frac{1}{2}\int_{0}^1\int_{0}^{1}F_{\phi}dp_r^N\otimes p_r^N- \frac{1}{2N^2}\sum_{i=1}^{N}\frac{\phi''(\tilde{X}_N^i(r))}{\pi}\to\frac{1}{2}\int_{0}^1\int_{0}^{1}F_{\phi}dp_r^\beta\otimes p_r^\beta,$$
		 we obtain
		$$\langle \phi'\H p_r^N, p_r^N\rangle g'(\langle \phi,p_r^N \rangle )\to \langle \phi'\H p^\beta_r, p^\beta_r\rangle g'(\langle \phi,p^\beta_r\rangle), \ \ \ \ a.s.$$
	Since
		$$\int_{0}^{t}\left|\langle \phi'\H p_r^N, p_r^N\rangle g'(\langle \phi,p_r^N \rangle )\right|dr\leq \int_{0}^{t}\frac{||\phi''||_{\infty}}{\pi}(1+\frac{1}{N^2})||g'||_{\infty}dr,$$
		By dominated convergence theorem, we get
		\begin{equation}\label{est5.16}
			\beta \int_{0}^{t}\langle \phi'\H p_r^N, p_r^N\rangle g'(\langle \phi,p_r^N\rangle )dr\to \frac{\beta}{2}\int_{0}^{t}\Big(\int_{0}^1\int_{0}^{1}F_{\phi}dp_r^\beta\otimes p_r^\beta \Big)\cdot g'(\langle \phi,p^\beta_r\rangle )dr.
		\end{equation}
		It follows from \eqref{EST5}, \eqref{est5.16} and the weak convergence of $p^N_t$ to $p_t$ that
		\begin{equation}\label{est5.17}
			M_t^\G(p^N)\cdot H(p^N\big|_{[0,s]}) \to  M_t^\G(p^\beta)\cdot H(p^\beta\big|_{[0,s]}) , \ \ \ a.s. .
		\end{equation}
	Note that
		\begin{equation*}
			\tilde{\E}\big[\big|M_t^\G(p^N)-M_s^\G(p^N)\big|\cdot \big|H(p^\beta\big|_{[0,s]}) \big|\big]<\infty, 
		\end{equation*}
		
		then combing with \eqref{est5.16} and \eqref{est5.17}, we obtain again by dominated convergence theorem that
		\begin{equation}\label{est5.19}
			\begin{split}
				&\tilde{\E}\big[\big(\G(p^\beta_t)-\G(p^\beta_s)-\int_{s}^{t}\mathrm{\mathbf{L}}_\beta \G(p^\beta_r) dr\big)\cdot H(p^\beta\big|_{[0,s]}) \big]\\
				&=\lim_{N\to \infty}\tilde{\E}\big[\big(\G(p^N_t)-\G(p^N_s)-\int_{s}^{t}\mathrm{\mathbf{L}}_\beta \G(p^N_r) dr\big)\cdot H(p^N\big|_{[0,s]}) \big] .
			\end{split}
		\end{equation}

		Note that when $k>0$, $e_k( x)=-e_k(-x)$, which implies
		
		\begin{equation*}
			Q^\gamma_{\mu}(\phi)=\int_0^1\int_0^1\phi'(x)\phi'(y)\bigg(\sum_{k=1}^{\infty}\frac{2}{k^{2\gamma}}\cos\big(2\pi k (\int_0^1\mathbbm{1}_{(x\land y,x\lor y]} (z) d\mu(z)) \big ) \bigg)d\mu(x)  d\mu(y)  .
		\end{equation*}
Let
		
		
		
		$$\mathrm{\mathbf{L}}_\beta^N \G(\mu)=\beta\langle \phi'\H\mu, \mu\rangle g'(\langle \phi,\mu \rangle )+K_2^{N,\gamma}g'(\langle \phi,\mu \rangle)\langle \mu, \phi''\rangle+g''(\langle \mu, \phi\rangle)Q_{\mu}^{\gamma,N}(\phi)  .$$

		Since 
	
	\begin{equation*}\label{esteq1}
		\left|\sum_{k=N}^{\infty}\frac{2}{k^{2\gamma}}\cos\big(2\pi k (\int_{0}^1\mathbbm{1}_{(x\land y,x\lor y]} (z) \mu(dz)) \big)\right|<\sum_{k=N}^{+\infty}\frac{2}{k^{2\gamma}}\ra 0, \ \mathrm{as}\ N\ra \infty , 
	\end{equation*}
it then follows from \eqref{est5.19} that
		\begin{equation*}
			\begin{split}
				&\tilde{\E}\big[\big(\G(p^\beta_t)-\G(p^\beta_s)-\int_{s}^{t}\mathrm{\mathbf{L}}_\beta \G(p^\beta_r) dr\big)\cdot H(p^\beta\big|_{[0,s]}) \big]\\
				&=\lim_{N\to \infty}\tilde{\E}\big[\big(\G(p^N_t)-\G(p^N_s)-\int_{s}^{t}\mathrm{\mathbf{L}}_\beta^N \G(p^N_r)dr\big)\cdot H(p^N\big|_{[0,s]}) \big] .
			\end{split}
		\end{equation*}
		
		So, to prove \eqref{estgoal}, we only need to prove
		\begin{equation}\label{est16}
			\lim_{N\to \infty}\tilde{\E}\big[\big(\G(p^N_t)-\G(p^N_s)-\int_{s}^{t}\mathrm{\mathbf{L}}^N_\beta \G(p^N_r) dr\big)\cdot H(p^N\big|_{[0,s]}) \big] =0.
		\end{equation}
		
		In fact, due to $Law(p^N)=P^N$, we have
		\begin{equation*}
			\begin{split}
				&\tilde{\E}\big[\big(\G(p^N_t)-\G(p^N_s)-\int_{s}^{t}\mathrm{\mathbf{L}}_\beta^N \G(p^N_r)dr\big)\cdot H(p^N\big|_{[0,s]}) \big] \\
				&=\E\big[\big(\G(L_N(t))-\G(L_N(s))-\int_{s}^{t}\mathrm{\mathbf{L}}_\beta^N \G(L_N(r))dr\big)\cdot H(L_N\big|_{[0,s]}) \big] .
			\end{split}
		\end{equation*}
		On the other hand , we have by It\^o's formula that
		
		\begin{equation*}
			\begin{split}
				&d \G(L_N(t))\\
				&=g'(\langle \phi, L_N(t)\rangle) \frac{\beta}{N^{2}}\sum_{i=1}^{N}\phi'(X_N^i(t))\cdot \left(\sum_{j=1,j\neq i}^{N}\cot\big(\pi(X_N^i(t)-X_N^j(t))\big)\right) d t\\
				&+g'(\langle \phi, L_N(t)\rangle) \frac{1}{2N}\sum_{i=1}^{N}\phi''(X_N^i(t))\sum_{k=-N}^{N}\frac{1}{|k|^{2\gamma}}e_k^2(\frac{i}{N}) d t\\
				&+g''(\langle \phi, L_N(t)\rangle)\left(\frac{1}{N^2}\sum_{i,j=1}^N\phi'(X_N^i(t))\phi'(X_N^j(t))\sum_{k=0}^{N}\frac{2}{k^{2\gamma}}\cos(2\pi k\frac{i-j}{N}) \right) dt +dM_N^{g,\phi}(t)\\
				&=(I)dt+(J)dt+(K)dt+dM_N^{g,\phi}(t),
			\end{split}
		\end{equation*}
		where $M_N^{g,\phi}(t)$ is a $\P-$local martingale. one can check that
		
		\begin{equation*}
			(J)=K_2^{N,\gamma}g'(\langle\phi,L_N(t)\rangle)\langle L_N(t), \phi''\rangle ;
		\end{equation*}
		
		\begin{equation*}
			(K)=g''(\langle L_N(t), \phi\rangle)Q_{L_N(t)}^{\gamma,N}(\phi) ;
		\end{equation*}
		
		\begin{equation*}
			(I)=\beta g'(\langle \phi,L_N(t) \rangle)\langle \phi'\H L_N(t), L_N(t)\rangle.
		\end{equation*}
		Therefore, \eqref{est16} is proved. We finish the proof.

	\end{proof}

\begin{remark}
	According to our method, the well-posedness of the particle model for $t\in [0,\infty)$ is based on the fact that there is no collision between particles in finite time. One could reasonably separate the colliding particles so that the dynamics of particles can be continued in a similar way after collisions and the well-posedness (at least in sense of weak solution) for $t\in [0,\infty)$ might not be a problem even when $\beta\in [0, \lambda_0^\gamma)$. However, it might be difficult to prove the convergence of the empirical measure processes to weak solutions to \eqref{withrepul} with the effect of the separating mechanism. Also, it is still a problem whether there exists collisions in paths with positive probability when $\beta\in [0, \lambda_0^\gamma)$. 
\end{remark}


	
	\subsection{R\'enyi entropy of $p^\beta_t$ for $\beta>\frac{p}{2}\pi K_1^\gamma$}
	\quad 	In this subsection, we will give a $\L^p$ regularity estimate of the limit measure. For $N-$particle model \eqref{model1}, we construct 
	the mollified R\'enyi entropy $H_{p}^N(L_N(t)) $
	$$H_p^N(L_N(t))=\sum_{j=1}^{N}\frac{1}{N}\left(\frac{2\pi}{N|e^{2\pi i X_N^{j}(t)}-e^{2\pi i X_N^{j-1}(t)}|}\right)^{p-1}, \ \ \ \  \ p>1.$$

	\vskip 3mm

	\begin{theorem}\label{lpregular}
		Suppose that $\gamma>\frac{3}{2}, \beta>\frac{p}{2}\pi K_1^\gamma$ for $p>1$ and $\mu_0\in \P(\T)$. If the particle model \eqref{model1} satisfies the condition \eqref{cod2'},
		then the fluctuating hydrodynamic limit $p^\beta_t$ is a solution to the martingale problem $(MP)_{\mu_0}^{K_2^\gamma, \beta, \gamma}$ of \eqref{withrepul} and is absolutely continuous w.r.t Leb for all $t> 0$, $\tilde{\P}-$a.s.. Moreover, 
		\begin{equation}\label{est19}
			\left(\tilde\E [\textbf{H}_p(p^\beta_t)]\right)^{\frac{1}{p}}\leq  \frac{a_p}{b_p}\frac{1}{1-e^{-2\frac{b_p}{a_p}t}}, \ \ \ \ \forall \, t>0,
		\end{equation}
		where $b_p=\frac{p-1}{4}(\frac{2\beta}{\pi}-pK_1^\gamma)$ and $a_p=(p-1)\pi \beta+1$.
	\end{theorem}
	\begin{proof}
		We estimate $\E[ H_p^N(L_N(t))]$ by a similar method as we used in Proposition \ref{abs}. In fact, by It\^o's formula, 
		\begin{equation*}
			\begin{split}
				&d_t H_p^N(L_N(t))\\&=-\sum_{j=1}^{N}\frac{(p-1)2^{p-1}\pi^{p}Z_N^j(t)}{N^p|e^{2\pi i X_N^{j}(t)}-e^{2\pi i X_N^{j-1}(t)}|^{p-1}}\circ(dX^j_N(t)-dX_N^{j-1}(t))\\
				&=-\sum_{j=1}^{N}\frac{(p-1)2^{p-1}\pi^p Z_N^j(t)}{N^p|e^{2\pi i X_N^{j}(t)}-e^{2\pi i X_N^{j-1}(t)}|^{p-1}}\frac{\beta}{N}\cdot\\
				& \left(\sum_{l=1,l\neq j}^{N}\frac{1}{\tan(\pi(X^j_N(t)-X_N^{l}(t)))}-\sum_{l=1,l\neq j-1}^{N}\frac{1}{\tan(\pi(X^{j-1}_N(t)-X_N^{l}(t)))}\right)dt\\
				&+\sum_{j=1}^{N}\frac{(p-1)^2 2^{p-1}\pi^{p+1}(Z_N^j(t))^2}{2N^p|e^{2\pi i X_N^{j}(t)}-e^{2\pi i X_N^{j-1}(t)}|^{p-1}}\sum_{k=-N}^{N} \frac{1}{k^{2\gamma}}\left|e_k(\frac{j}{N})-e_k(\frac{j-1}{N})\right|^2 dt\\
				&+\sum_{j=1}^{N}\frac{(p-1)2^{p-1}\pi^{p+1}\big((Z_N^j(t))^2+1\big)}{2N^p|e^{2\pi i X_N^{j}(t)}-e^{2\pi i X_N^{j-1}(t)}|^{p-1}}\sum_{k=-N}^{N} \frac{1}{k^{2\gamma}}\left|e_k(\frac{j}{N})-e_k(\frac{j-1}{N})\right|^2 dt+dM_N(t)\\
				&=(2\pi)^{p-1}\beta a(t)dt+b(t)dt+c(t)dt+dM_N(t).
			\end{split}
		\end{equation*}
		Without repeating the details, we obtain by applying Lemma \ref{cotineq} for $\displaystyle f(x)=\frac{1}{|\sin( x\pi)|^{p-1}}$ that
		\begin{equation*}
			\begin{split}
				a(t)\leq&\sum_{j=1}^{N}\frac{(p-1)\pi}{N^{p}|e^{2\pi i X_N^{j}(t)}-e^{2\pi i X_N^{j-1}(t)}|^{p-1}}\\
				&-\sum_{j=1}^{N}\frac{2(p-1)\pi (Z_N^j(t))^2}{N^{p+2}|e^{2\pi i X_N^{j}(t)}-e^{2\pi i X_N^{j-1}(t)}|^{p-1}}
			\end{split}
		\end{equation*} 
		and 
		\begin{equation*}
			\begin{split}
				b(t)\leq&\sum_{j=1}^{N}\frac{(p-1)^22^{p-1}\pi^{p+1}K_1^{N,\gamma}(Z_N^j(t))^2}{N^{p+2}|e^{2\pi i X_N^{j}(t)}-e^{2\pi i X_N^{j-1}(t)}|^{p-1}},\\
			\end{split}
		\end{equation*} 
		\begin{equation*}
			\begin{split}
				c(t)\leq&\sum_{j=1}^{N}\frac{(p-1)2^{p-1}\pi^{p+1}K_1^{N,\gamma}}{N^{p+2}|e^{2\pi i X_N^{j}(t)}-e^{2\pi i X_N^{j-1}(t)}|^{p-1}\sin^2(\pi(X^j_N(t)-X_N^{j-1}(t)))}.
			\end{split}
		\end{equation*} 
		Notice that 
		\begin{equation*}
			\begin{split}
				b(t)+c(t)\leq&\sum_{j=1}^{N}\frac{(p-1)p2^{p-1}\pi^{p+1}K_1^{N,\gamma}(Z_N^j(t))^2}{N^{p+2}|e^{2\pi i X_N^{j}(t)}-e^{2\pi i X_N^{j-1}(t)}|^{p-1}}\\
				&+\sum_{j=1}^{N}\frac{(p-1)2^{p-1}\pi^{p+1}K_1^{N,\gamma}}{N^{p+2}|e^{2\pi i X_N^{j}(t)}-e^{2\pi i X_N^{j-1}(t)}|^{p-1}}.
			\end{split}
		\end{equation*} 
		 We have
		\begin{equation}\label{equ5.5}
			\begin{split}
				&\frac{d}{dt} \E\big[H_p^N(L_N(t))\big]\\
				\leq&(p-1)\pi \beta\E\big[H_p^N(L_N(t))\big]+\frac{(p-1)\pi^2K_1^\gamma}{N^2}\E\big[H_p^N(L_N(t))\big]\\
				&+\sum_{j=1}^{N}\frac{(p-1)2^{p-1}\pi^{p}\Big(p\pi K_1^{\gamma}-2\beta\Big) (Z_N^j(t))^2}{N^{p+2}|e^{2\pi i X_N^{j}(t)}-e^{2\pi i X_N^{j-1}(t)}|^{p-1}}.
			\end{split}
		\end{equation} 
		Notice that 
		$$\frac{1}{\tan^2(\pi (x-y))}=\frac{1}{\sin^2(\pi (x-y))}-1=\frac{4}{|e^{2\pi ix}-e^{2\pi i y}|^2}-1,$$
		then it follows from \eqref{equ5.5} that
		\begin{equation}\label{est9}
			\begin{split}
				&\frac{d}{dt} \E\big[H_p^N(L_N(t))\big]\\
				&\leq (p-1)\pi \beta\E\big[H_p^N(L_N(t))\big]+\frac{(p-1)\pi^2K_1^\gamma}{N^2}\E\big[H_p^N(L_N(t))\big]\\
				&-\sum_{j=1}^{N}\E\left[\frac{(p-1)2^{p+1}\pi^{p}\Big(2\beta-p\pi K_1^{\gamma}\Big)}{N^{p+2}|e^{2\pi i X_N^{j}(t)}-e^{2\pi i X_N^{j-1}(t)}|^{p+1}}\right]+\frac{(p-1)\pi\Big(2\beta-p\pi K_1^{\gamma}\Big)}{N^2}\E\big[H_p^N(L_N(t))\big].
			\end{split}
		\end{equation} 
		By H\"older inequality, we have
		$$\left(\sum_{j=1}^{N}\frac{1}{N^{p+1}|e^{2\pi i X_N^{j}(t)}-e^{2\pi i X_N^{j-1}(t)}|^{p+1}}\right)^{\frac{p-1}{p+1}}\cdot N^{\frac{2}{p+1}}\geq \sum_{j=1}^{N}\frac{1}{N^{p-1}|e^{2\pi i X_N^{j}(t)}-e^{2\pi i X_N^{j-1}(t)}|^{p-1}}.$$
		It yields 
		\begin{equation}\label{hold1}
			\sum_{j=1}^{N}\frac{1}{N^{p+2}|e^{2\pi i X_N^{j}(t)}-e^{2\pi i X_N^{j-1}(t)}|^{p+1}} \geq \left(\sum_{j=1}^{N}\frac{1}{N^{p}|e^{2\pi i X_N^{j}(t)}-e^{2\pi i X_N^{j-1}(t)}|^{p-1}}\right)^{\frac{p+1}{p-1}}.
		\end{equation}
		Also, by H\"older's inequality, we have
		\begin{equation}\label{hold2}
			\E\big[(H_p^N(L_N(t)))^{\frac{p+1}{p-1}}\big]\geq \Big(\E\big[H_p^N(L_N(t))\big]\Big)^{\frac{p+1}{p-1}}.
		\end{equation}
		Applying \eqref{hold1} and \eqref{hold2} to \eqref{est9}, we have
		$$\frac{d}{dt} \E\big[H_p^N(L_N(t))\big]\leq  a_p\E\big[H_p^N(L_N(t))\big]-b_p\Big(\E\big[H_p^N(L_N(t))\big]\Big)^{\frac{p+1}{p-1}},$$ 
		where 
		$$a_p=(p-1)\pi \beta+1$$
		and 
		$$b_p=\frac{p-1}{4}(\frac{2\beta}{\pi}-pK_1^\gamma),$$
		when $N$ is large enough. For each $p>1$, we obtain by using comparison theorem of ODE that $\displaystyle \E\big[H_p^N(L_N(t))\big]\leq y_p(t)$, where $y_p(t)$ satisfies
		\begin{equation*}\label{odeeq}
			\Big(\ln\big(\frac{y_p}{|1-\frac{b_p}{a_p}y_p^{\frac{2}{p-1}}|}\big)\Big)'=a_p.
		\end{equation*}	
		Let $\displaystyle C^p_0=\frac{b_p}{a_p}-(H_p^N(L_N(0)))^{-\frac{2}{p-1}}$ and $\beta\geq \frac{p}{2}\pi K_1^\gamma$. We solve the equation above and find that 
		\begin{equation}\label{est11}
				\left(\E\big[H_p^N(L'_N(t))\big]\right)^{\frac{2}{p-1}}\leq \frac{1}{\frac{b_p}{a_p}-C^p_0e^{-\frac{2a_p}{p-1}t}}=\frac{a_p}{b_p}+\frac{C_0^p e^{-2\frac{b_p}{a_p}t}}{\frac{b_p}{a_p}(\frac{b_p}{a_p}-C_0^p e^{-2\frac{b_p}{a_p}t})},
		\end{equation}
		if $H_p^N(L_N(0))>\big(\frac{a_p}{b_p}\big)^{\frac{p-1}{2}}$. When $H_p^N(L_N(0))\leq\big(\frac{a_p}{b_p}\big)^{\frac{p-1}{2}}$, 
		\begin{equation}\label{est11a}
		\left(\E\big[H_p^N(L'_N(t))\big]\right)^{\frac{2}{p-1}}\leq \frac{1}{\frac{b_p}{a_p}-C_0^pe^{-2\frac{b_p}{a_p}t}}\leq\frac{a_p}{b_p}.
		\end{equation}
		

		On the other hand, according to Theorem \ref{introthm1} and Proposition \ref{ratecon}, there exists a unique fluctuating hydrodynamic limit $p^\beta_t$. By lower semicontinuity of $\textbf{H}_p$ for $p>1$ and Fatou's lemma, we have
		
		\begin{equation}\label{est5.4}
			\tilde\E [\textbf{H}_p(p^\beta_t)]\leq \liminf_{N\to\infty} \tilde\E [\textbf{H}_p(p^N_t)].
		\end{equation}
		Notice that
		$$\textbf{H}_p(p^N_t)\leq H_p^N(p^N_t),$$
		and 
		$$\tilde\E [H_p^N(p^N_t)]=\E [H^N_p(L_N(t))];$$
		therefore combing with \eqref{est11}, \eqref{est11a}, \eqref{est5.4} and Fatou's lemma, we have
		$$\tilde\E[\textbf{H}_p(p^\beta_t)]\leq \liminf_{N\to \infty} \tilde{\E} [H^N_p(p^N_t)]<  \left(\frac{a_p}{b_p}\frac{1}{1-e^{-2\frac{b_p}{a_p}t}}\right)^{\frac{p-1}{2}}, \   \ \ \  \forall \ t>0,$$
		under the condition $\beta>\frac{p \pi K_1^\gamma}{2}$. \eqref{est19} is proved. Again by the similar argument used in proving finiteness of $\textbf{S}(p_t)$, we can also prove that, for $\forall \ t>0$, it holds
		\begin{equation*}\label{est17}
			\textbf{H}_p(p^\beta_t)<\infty, \ \ \ \ a.s.
		\end{equation*}
		for any initial measure $\mu_0$ as long as $L_N(0)$ weakly converges to $\mu_0$. 
	\end{proof}
	\vskip3mm

	\vskip 5mm

	\section{Appendice}
	\begin{lemma}\label{cotineq}
		Let $\{Y_N^i\}_{i=1}^N$ be $N$ points in $\R$. Suppose that $Y_N^1<Y_N^2<\ldots <Y_N^N$ and $$|Y_N^1-Y_N^N|<1,$$
		then for any $p\geq 0$, it holds
		\begin{equation}
			\begin{split}
				&\sum_{i=1}^{N}|\sin^{-p}(\pi(Y_N^i-Y_N^{i-1}))|\sum_{j=1,j\neq i,i-1}^{N}\frac{1}{\tan(\pi(Y^i_N-Y_N^{j})\tan(\pi(Y^{i-1}_N-Y_N^{j})}\\
				&\leq 2(1-\frac{1}{N})\sum_{i=1}^{N}|\sin^{-p}(\pi(Y_N^i-Y_N^{i-1}))|\big(Z_N^i\big)^2.
			\end{split}
		\end{equation}
	\end{lemma}
	\begin{proof}
		We firstly consider the case when $\forall \, i\in\{1,2,\ldots, N\}$, 
		\begin{equation}\label{wlog}
			|Y_N^i-Y_N^{i-1}|\leq \frac{1}{2}.
		\end{equation}
		We can decompose 
		$$\sum_{j=1,j\neq i,i-1}^{N}\frac{1}{\tan(\pi(Y^i_N-Y_N^{j})\tan(\pi(Y^{i-1}_N-Y_N^{j})}$$ into three parts
		\begin{equation*}
			\begin{split}
				S^i_{+}&=\sum_{j\in B^i_{+,1}\cup B^i_{+,2}}\frac{1}{\tan(\pi(Y^i_N-Y_N^{j})\tan(\pi(Y^{i-1}_N-Y_N^{j})}\\
				S^i_{-}&=\sum_{j\in B^i_{-,1}\cup B^i_{-,2}}\frac{1}{\tan(\pi(Y^i_N-Y_N^{j})\tan(\pi(Y^{i-1}_N-Y_N^{j})}\\
				S^i_{c}&=\sum_{j\in B^i_{c} }\frac{1}{\tan(\pi(Y^i_N-Y_N^{j})\tan(\pi(Y^{i-1}_N-Y_N^{j})},
			\end{split}
		\end{equation*}
		where
		\begin{equation*}
			\begin{split}
				B^i_{+,1}&=\left\{j: 0<Y_N^i-Y_N^j\leq\frac{1}{2}, \ 0<Y_N^{i-1}-Y_N^j\leq\frac{1}{2},\  j< i-1\right\} \\
				B^i_{+,2}&=\left\{j: -1<Y_N^i-Y_N^j\leq-\frac{1}{2}, \ -1<Y_N^{i-1}-Y_N^j\leq-\frac{1}{2},\  j> i\right\} \\
				B^i_{-,1}&=\left\{j:\ -\frac{1}{2}\geq Y^{i-1}_N-Y_N^{j}<0,\ -\frac{1}{2}\geq Y^{i}_N-Y_N^{j}<0, \  j>i\right\} \\
				B^i_{-,2}&=\left\{j: \frac{1}{2}<Y_N^i-Y_N^j\leq1, \ \frac{1}{2}<Y_N^{i-1}-Y_N^j\leq1,\  j< i-1\right\}  \\
				B^i_{c}&=\left\{1,\ldots, N\right\}/B^i_{-,1}\cup B^i_{-,2}\cup B^i_{+,1}\cup B^i_{+,2}
			\end{split}
		\end{equation*}
		We define 
		$$L^i_{+,1}=\min_{j\in B^i_{+,1}} j,\ \ L^i_{+,2}=\max_{j\in B^i_{+,2}} j ,$$
		and 
		$$L^i_{-,1}=\max_{j\in B^i_{-,1}} j,\ \ L^i_{-,2}=\max_{j\in B^i_{-,2}} j .$$
		Here, if the set $B_{+,k}^i$ ($B_{-,k}^i$) is empty, we define the corresponding number $L_{+,k}^i=0$ ($L_{-,k}^i=0$). It is obvious that $S^i_{c}\leq0$. To estimate $S^i_{+}(t)$, we firstly let $y_N^i=Y_N^i+1$, if $Y_N^i-Y_N^j<0$, otherwise $y_N^i=Y_N^i$. Then 
		\begin{equation*}
			\begin{split}
				S^i_{+}&=\sum_{j\in B^i_{+,1}}\frac{1}{\tan(\pi(Y^i_N-Y_N^{j})\tan(\pi(Y^{i-1}_N-Y_N^{j})}\\
				&+ \sum_{j\in B^i_{+,2}}\frac{1}{\tan(\pi(y^i_N-y_N^{j})\tan(\pi(y^{i-1}_N-y_N^{j})}.
			\end{split}
		\end{equation*}
		Note that when $\tan A \cdot \tan B\geq0$,
		$$\tan(A+B)=\frac{\tan A+\tan B}{1-\tan A\cdot \tan B}\geq \tan A+\tan B.$$
		It follows that 
		$$\tan(\pi(Y^i_N-Y_N^{j}))\geq \sum_{k=1}^{i-j}\tan(\pi(Y^{j+k}_N-Y^{j+k-1}_N)).$$
		Therefore, if $L^i_{+,1}>0$, then $B^i_{+,2}$ is empty, and applying the inequality above gives
		\begin{equation*}
			\begin{split}
				&S_{+}^i=\sum_{j\in B^i_{+,1}}\frac{1}{\tan(\pi(Y^i_N(t)-Y_N^{j})\tan(\pi(Y^{i-1}_N(t)-Y_N^{j})}\\
				&\leq \sum_{j=L_{+,1}^i}^{i-2}\frac{1}{\Big(\sum_{k=1}^{i-j}\tan(\pi(Y^{j+k}_N-Y^{j+k-1}_N))\Big)\cdot \Big(\sum_{k=1}^{i-1-j}\tan(\pi(Y^{j+k}_N-Y^{j+k-1}_N))\Big)}.
			\end{split}
		\end{equation*}
		By Cauchy's inequality, we have
		\begin{equation*}
			\begin{split}
				&\sum_{j\in B^i_{+,1}}\frac{1}{\tan(\pi(Y^i_N-Y_N^{j})\tan(\pi(Y^{i-1}_N-Y_N^{j})}\\
				&\leq \sum_{j=L_{+,1}^i}^{i-2}\frac{1}{(i-j)^2(i-j-1)^2}\left(\sum_{m=1}^{i-j}\sum_{n=1}^{i-1-j}\cot(\pi(Y^{j+m}_N-Y^{j+m-1}_N))\cot(\pi(Y^{j+n}_N-Y^{j+n-1}_N))\right).
			\end{split}
		\end{equation*}
		Similarly, if $B^i_{+,2}$ is not empty, then for $j\in B^i_{+,2}$, we have
		$$\tan(\pi(y^i_N-y_N^{j}))\geq \sum_{k=1}^{N-j+i}\tan(\pi(y^{\left\{\frac{(j+k)}{N}\right\}N}_N-y^{\left\{\frac{(j+k-1)}{N}\right\}N}_N)).$$
		It follows that 
		\begin{equation*}
			\begin{split}
				&\sum_{j\in B^i_{+,2}}\frac{1}{\tan(\pi(y^i_N-y_N^{j})\tan(\pi(y^{i-1}_N-y_N^{j})}\\
				&\leq \sum_{j=L_{+,2}^i}^N\frac{1}{(N-j+i)^2(N-j+i-1)^2}\\
				&\cdot\left(\sum_{m=1}^{N-j+i}\sum_{n=1}^{N-j+i-1}\cot(\pi(y^{\left\{\frac{(j+m)}{N}\right\}N}_N-y^{\left\{\frac{(j+m-1)}{N}\right\}N}_N))\cot(\pi(y^{\left\{\frac{(j+n)}{N}\right\}N}_N-y^{\left\{\frac{(j+n-1)}{N}\right\}N}_N))\right).
			\end{split}
		\end{equation*}
		
		Now we conclude that 
		\begin{equation}\label{S+}
			S_+^i\leq\sum_{k=1}^{T^{i}_{k,+}}\frac{1}{(k+1)^2k^2}\left(\sum_{m=1}^{k+1}\sum_{n=1}^{k}Z_N^{N\left\{\frac{i-k-1+m}{N}+1\right\}}Z_N^{N\left\{\frac{i-k-1+n}{N}+1\right\}}\right),
		\end{equation}
		where 
		\begin{equation*}
			T^{i}_{k,+}=\left\{
			\begin{aligned}
				&i-L^i_{+,1}-1, \  \ \ \textup{if } \ L^i_{+,2}= 0;\\
				&N-L^i_{+,2}+i-1, \ \ \  \textup{if} \ L^i_{+,2}\neq 0.
			\end{aligned}
			\right.
		\end{equation*}
		
		We can estimate $S^i_{-}$ in a similar way. Without showing the details, we get 
		
		\begin{equation*}
			\begin{split}
				&S_-^i\leq\sum_{k=1}^{T^{i}_{k,-}}\frac{1}{(k+1)^2k^2}\left(\sum_{m=1}^{k+1}\sum_{n=1}^{k}Z_N^{N\left\{\frac{i+m}{N}-1\right\}}Z_N^{N\left\{\frac{i+n-1}{N}-1\right\}}\right),
			\end{split}
		\end{equation*}
		where 
		\begin{equation*}
			T^{i}_{k,-}=\left\{
			\begin{aligned}
				&L^i_{-,1}-i, \  \ \ \textup{if } \ L^i_{+,2}= 0;\\
				&L^i_{-,2}+N-i, \ \ \  \textup{if} \ L^i_{+,2}\neq 0.
			\end{aligned}
			\right.
		\end{equation*}

		Because of the rearrangement inequality, we have, for fixed $k,m,n$, 
		$$\sum_{i=1}^{N}Z_N^{N\left\{\frac{i-k-1+m}{N}+1\right\}}Z_N^{N\left\{\frac{i-k-1+n}{N}+1\right\}}(t)\leq \sum_{i=1}^{N}\big(Z_N^i\big)^2.$$
		Thus,
		\begin{equation*}
			\begin{split}
				&\sum_{i=1}^{N}|\sin^{-p}(\pi(Y_N^i-Y_N^{i-1}))|\sum_{j=1,j\neq i,i-1}^{N}\frac{1}{\tan(\pi(Y^i_N-Y_N^{j})\tan(\pi(Y^{i-1}_N-Y_N^{j})}\\
				&\leq \sum_{i=1}^{N}|\sin^{-p}(\pi(Y_N^i-Y_N^{i-1}))|(S_+^i+S_-^i)\\
				&<\sum_{k=1}^{N-1}\frac{1 }{k^2(k+1)^2}\sum_{m=1}^{k+1}\sum_{n=1}^{k}\left(\sum_{i=1}^{N}|\sin^{-p}(\pi(Y_N^i-Y_N^{i-1}))|Z_N^{N\left\{\frac{i-k-1+m}{N}+1\right\}}Z_N^{N\left\{\frac{i-k-1+n}{N}+1\right\}}\right)\\
				&+\sum_{k=1}^{N-1}\frac{1 }{k^2(k+1)^2}\sum_{m=1}^{k+1}\sum_{n=1}^{k}\left(\sum_{i=1}^{N}|\sin^{-p}(\pi(Y_N^i-Y_N^{i-1}))|Z_N^{N\left\{\frac{i+m}{N}-1\right\}}Z_N^{N\left\{\frac{i+n-1}{N}-1\right\}}\right)\\
				&\leq 2\sum_{k=1}^{N-1}\frac{1 }{k^2(k+1)^2}\sum_{m=1}^{k+1}\sum_{n=1}^{k}\sum_{i=1}^{N}|\sin^{-p}(\pi(Y_N^i-Y_N^{i-1}))|\big(Z_N^i\big)^2\\
				&=\sum_{k=1}^{N-1}\frac{2}{k(k+1)}\sum_{i=1}^{N}|\sin^{-p}(\pi(Y_N^i-Y_N^{i-1}))|\big(Z_N^i\big)^2\\
				&< 2(1-\frac{1}{N})\sum_{i=1}^{N}|\sin^{-p}(\pi(Y_N^i-Y_N^{i-1}))|\big(Z_N^i\big)^2.
			\end{split}
		\end{equation*}
		
		\quad Next, when the assumption \eqref{wlog} is not satisfied, there exists exactly one $i^*$ such that $\frac{1}{\tan(\pi(Y^{i^*}_N(t)-Y_N^{i^*-1}(t))}<0$. In such case, it holds
		$$\frac{1}{\tan(\pi(Y^{i^*}_N(t)-Y_N^{j}(t))\tan(\pi(Y^{i^*-1}_N(t)-Y_N^{j}(t))}<0, \ \ \ \forall \, j\neq i^*,i^*-1,$$
		which implies that
		\begin{equation*}
			\begin{split}
				&\sum_{i=1}^{N}|\sin^{-p}(\pi(Y_N^i-Y_N^{i-1}))|\sum_{j=1,j\neq i,i-1}^{N}\frac{1}{\tan(\pi(Y^{i}_N(t)-Y_N^{j}(t))\tan(\pi(Y^{i-1}_N(t)-Y_N^{j}(t))}\\
				&\leq \sum_{i=1, i\neq i^*}^{N}|\sin^{-p}(\pi(Y_N^i-Y_N^{i-1}))|\sum_{j=1,j\neq i,i-1}^{N}\frac{1}{\tan(\pi(Y^{i}_N(t)-Y_N^{j}(t))\tan(\pi(Y^{i-1}_N(t)-Y_N^{j}(t))}\\
				&=\sum_{i=1, i\neq i^*}^{N}|\sin^{-p}(\pi(Y_N^i-Y_N^{i-1}))|( S_+^i+S_-^i).
			\end{split}
		\end{equation*}
		Note that in such case, we have 
		$$j\in B^i_{+,1}\cup B^i_{+,2}  \ \ \ \ \forall \, j\neq i^*,$$
		which means 
		$$\sum_{i=1, i\neq i^*}^{N}  S_-^i=0.$$
		Subsequently, following the estimate \eqref{S+} which was derived based on the assumption \eqref{wlog}, we obtain
		\begin{equation*}
			\begin{split}
				&\sum_{i=1, i\neq i^*}^{N}|\sin^{-p}(\pi(Y_N^i-Y_N^{i-1}))|\sum_{j=1,j\neq i,i-1}^{N}\frac{1}{\tan(\pi(Y^{i}_N(t)-Y_N^{j}(t))\tan(\pi(Y^{i-1}_N(t)-Y_N^{j}(t))}\\
				&<(1-\frac{1}{N})\sum_{i=1}^{N}|\sin^{-p}(\pi(Y_N^i-Y_N^{i-1}))|\big(Z_N^i\big)^2.
			\end{split}
		\end{equation*}
		Now we have completed the proof.
	\end{proof}


	\bibliographystyle{alpha}
	\bibliography{DW}
\end{document}